\documentclass[12pt,a4paper]{article}
\usepackage[dvipdfmx]{graphicx,color}
\usepackage{amsmath,amsfonts,amssymb}
\usepackage{hyperref}
\usepackage{ulem}

\usepackage[numbers]{natbib}
\newcommand{\vc}[1]{{\boldsymbol #1}} 
\newcommand{\vcn}[1]{{\bf #1}}

\newcommand{\sr}[1]{{\cal #1}}
\newcommand{\dd}[1]{\mathbb{#1}}

\newcommand{\rmn}[1]{\if#11I\else {\if#12I\hspace{-0.12ex}I\hspace{-0.85ex}\else {\if #13I\hspace{-0.16ex}I\hspace{-0.16ex}I\hspace{-1.6ex}\else I\hspace{-1.2ex}V \fi} \fi} \fi}

\newcommand{\eqn}[1]{(\ref{eqn:#1})}
\newcommand{\lem}[1]{Lemma~\ref{lem:#1}}
\newcommand{\cor}[1]{Corollary~\ref{cor:#1}}
\newcommand{\thr}[1]{Theorem~\ref{thr:#1}}

\newcommand{\supp}[1]{Supplement~\ref{supp:#1}}
\newcommand{\rem}[1]{Remark~\ref{rem:#1}}

\newcommand{\app}[1]{Appendix~\ref{app:#1}}
\newcommand{\sectn}[1]{Section~\ref{sect:#1}}

\newcommand{\lemt}[1]{\ref{lem:#1}}

\newcommand{\thrt}[1]{\ref{thr:#1}}

\newcommand{\figt}[1]{\ref{fig:#1}}

\newcommand{\sect}[1]{\ref{sect:#1}}

\newcommand{\ol}{\overline}
\newcommand{\ul}{\underline}
\newcommand{\pend}{\hfill \thicklines \framebox(6.6,6.6)[l]{}}

\newenvironment{proof*}[1]{\noindent {\sc  #1} \rm}{\pend}

\newtheorem{theorem}{Theorem}[section]
\newtheorem{lemma}{Lemma}[section]

\newtheorem{remark}{Remark}[section]

\newtheorem{corollary}{Corollary}[section]

\newcommand{\setsection}[2] {
\setcounter{section}{#1}
\setcounter{subsection}{0}
\setcounter{equation}{0}
\setcounter{conjecture}{0}
\setcounter{assumption}{0}
\setcounter{question}{0}
\setcounter{definition}{0}
\setcounter{theorem}{0}
\setcounter{corollary}{0}
\setcounter{lemma}{0}
\setcounter{proposition}{0}
\setcounter{remark}{0}
\setcounter{appen}{0}
\setsection*{\large \bf \thesection. #2}}

\newcommand{\setnewcounter} {
\setcounter{subsection}{0}
\setcounter{equation}{0}
\setcounter{conjecture}{0}
\setcounter{assumption}{0}
\setcounter{question}{0}
\setcounter{definition}{0}
\setcounter{theorem}{0}
\setcounter{corollary}{0}
\setcounter{lemma}{0}
\setcounter{proposition}{0}
\setcounter{remark}{0}
}

\begin{document}
\title{\bf \Large A unified approach for large queue asymptotics\\ in a heterogeneous multiserver queue}

\author{Masakiyo Miyazawa\\ Tokyo University of Science}
\date{\today, supplemented version 2}

\maketitle

\begin{abstract}
We are interested in a large queue in a $GI/G/k$ queue with heterogeneous servers. For this, we consider tail asymptotics and weak limit approximations for the stationary distribution of its queue length process in continuous time under a stability condition. Here, two weak limit approximations are considered. One is when the variances of the inter-arrival and/or service times are bounded, and the other is when they get large. Both require a heavy traffic condition. Tail asymptotics and heavy traffic approximations have been separately studied in the literature. We develop a unified approach based on a martingale produced by a good test function for a Markov process to answer both problems. 
\end{abstract}
  
\begin{quotation}

\medskip

\end{quotation}

\section{Introduction}
\label{sect:introduction}

We are interested in a large queue in a $GI/G/k$ queue with FCFS (First-Come First-Served) service discipline. This queueing model has a single queue and $k$ servers who may have different service time distributions. Throughout the paper, $k$ is an arbitrary positive integer, but fixed. We refer this model as a heterogeneous multiserver queue. As a queueing model, this queue is basic and classical, but it is known to be difficult to get analytical results for this model. Because of this, theoretical studies on it have been directed to asymptotic analysis. Among them, behavior of a large queue is one of their major interests.

In this paper, we study the large queue behavior through tail asymptotics and weak limit approximations for the stationary distribution of the queue length process of the heterogeneous multiserver queue. We describe this queueing model by a continuous time Markov process, and assume that it has a stationary distribution, that is, it is stable. Here, the weak limit approximation is meant to derive a limiting distribution for approximation from a sequence of the queueing models under appropriate scaling of their queue length. Such an approximation traditionally assumes that the variances of the inter-arrival and/or service times are bounded, and is obtained by increasing the traffic intensity to one from below, which is called a heavy traffic condition. We also consider another approximation under the heavy traffic condition by increasing the variances of the inter-arrival and/or service times, which is referred to as large variance approximation.

The tail asymptotics and heavy traffic approximations have been separately studied in the literature. For example, the tail asymptotics are obtained in Neuts and Takahashi \cite{NeutTaka1981} and Sadowsky and Szpankowski \cite{SadoSzpa1995} for the stationary waiting time distributions and some other characteristics at arrival instants of customers, while the heavy traffic approximation is obtained for the stationary waiting time distribution in a $GI/G/k$ queue with homogeneous servers by K{\"o}llerstr{\"o}m \cite{Koll1974}. This heavy traffic approximation is related to a diffusion approximation for a process limit, which has been widely studied. In particular for the heterogeneous multiserver queue, a diffusion approximation is obtained for the queue length process in Chen and Ye \cite{ChenYe2011a}. Unlike those studies, we aim to give a unified approach to simultaneously answer both of the tail asymptotic problem and the heavy traffic approximation on the stationary distribution of the queue length process in continuous time.

For this unified approach, we describe the heterogeneous multiserver queue by a piecewise deterministic Markov process due to Davis \cite{Davi1984,Davi1993}, which is called a PDMP for short. A major difficulty in studying its stationary distribution comes from state transitions at arriving and departing instants of customers, which are particularly complicated when some of the servers are idle. In the PDMP, the effect of those transitions can be removed from system equations for its dynamics by the so called boundary condition. Since the boundary in this terminology is different from boundary states for the queue length process, we refer to the boundary condition of \cite{Davi1984} as a terminal condition in this paper. A key step of our approach is to find a good test function satisfying this terminal condition, which produces a martingale from the PDMP. This idea will be shown to be useful not only for the tail asymptotic problem but also for the weak limit approximations.
 
Obviously, the tail asymptotics of the stationary distribution is crucially influenced by the heaviness of the service time distributions. Because of this, different techniques have been employed for the light and heavy tail cases in the literature (e.g., see \cite{FossKors2012,Miya2011}). However, it may be also interesting to see what occurs when they are mixed. This question was answered for the stationary waiting time distribution by Sadowsky \cite{Sado1995}. We revisit this question for the stationary queue length distribution in continuous time, and get a similar answer to it (see \thr{tail asymptotic 1} and \cor{slow server}).

Our weak limit approximations for a large queue have the same spirit as Kingman \cite{King1962}'s two moment approximation for the stationary waiting time distribution of the $GI/G/1$ queue in heavy traffic, which is extended to the homogeneous multiserver queue by K{\"o}llerstr{\"o}m \cite{Koll1974}. Motivated by them, the author \cite{Miya2015} recently studied heavy traffic approximation of the stationary queue length distribution in continuous time in the $GI/G/1$ queue. Its basic idea is close to use the aforementioned terminal condition. In this paper, we refine it using the PDMP, and attack various asymptotic problems for the heterogeneous multiserver queue, which have been not studied in the line of Kingman \cite{King1962}'s two moment approximation.

The present approach is flexible in modeling assumptions. For example, it allows arrivals and service completions of different servers to occur simultaneously. This considerably complicates analysis on state transitions, and is often excluded by assuming additional assumptions such that the service time distributions have densities. We do not need such assumptions. We also demonstrate it to be applicable to the case that the arrival process is a superposition of the renewal processes. Furthermore, the approach may be amendable for some other queueing models. In particular, a queueing network is apparently more interesting for application. Braverman et al.\ \cite{BravDaiMiya2015} shows that its heavy traffic approximation can be studied in the same spirit (see \sectn{some} for more discussion).

This paper is made up by five sections. In \sectn{hetero}, we formally introduce the heterogeneous multiserver queue and the PDMP for describing it. We also discuss a martingale obtained from the PDMP. A representation of the PDMP using this martingale is used to derive a stationary equation with a certain good test function. Main results on the tail asymptotics and weak limit approximations are presented in \sectn{answers}. In this section, we also consider an extension of the heterogeneous multiserver queue so that it has multiple renewal arrivals. Theorems are proved in \sectn{proofs}, where we detail change of measure for the tail asymptotic because they may have their own interest. In \sectn{concluding}, we discuss possible extensions of the present approach. Technical lemmas are proved in the appendices. In the supplement, Lemmas 2.5 and 4.8 are proved, and the proofs of Lemmas 3.1 and 4.2 are fully detailed. Those proofs are omitted in a published version \cite{Miya2017}.

\section{Heterogeneous multiserver queue}
\label{sect:hetero}

We introduce the heterogeneous multiserver queue in \sectn{model}. We then introduce the Markov process and the related tools in Sections \sect{extended}--\sect{stationary}. These tools will be used to prove our asymptotic results (Theorems \thrt{tail asymptotic 1}--\thrt{variance 1}) that will be stated in \sectn{answers}. Developing the tools is an important part of the contribution of this paper, which may have independent interest. Readers who wish to read the statement of the theorems without
focusing on their proofs can go to \sectn{answers} directly after reading \sectn{model}, picking up the definitions of functions $\eta(v,\cdot)$ and $\xi(v, \cdot)$ in \sectn{good}.

\subsection{Model description and PDMP}
\label{sect:model}

By a heterogeneous $k$-server queue for $k \ge 1$, we mean that customers arrive at a system according to a renewal process and are served in FCFS (Fist-Come First-Served) manner by $k$ servers, numbered as $1,2,\ldots,k$, who have their own service time distributions, which may be different. It is assumed that service times at each server are $i.i.d.$ and independent of everything else. We denote the inter-arrival time distribution by $F_{0}$, which is assumed to have finite positive mean $1/\lambda_{0}$. Thus, $\lambda_{0}$ is the arrival rate. Denote the service time distribution for the $i$-th server by $F_{i}$, which is assumed to have a finite mean $m_{i}$. Let $\lambda_{i} = 1/m_{i}$, which is the service rate for $i \in K$, and let $\rho = \lambda_{0} /\sum_{i \in K} \lambda_{i}$, where
\begin{align*}
  K = \{1,2,\ldots,k\}.
\end{align*}
For convenience, we let $\ol{K} = \{0\} \cup K$, where $0$ represents arrivals of customers.

Let $T_{i}$ be a random variable subject to the distribution $F_{i}$ for $i \in \ol{K}$. To exclude trivial cases, we assume throughout the paper that 
\begin{itemize}
\item [(i)] $\dd{P}(T_{i} > 0) > 0$ for all $i \in \ol{K}$.
\end{itemize}
This condition is equivalent to that $\dd{E}(T_{i}) > 0$ for all $i \in \ol{K}$. It is notable that $F_{0}(0) (=\dd{P}(T_{i} = 0)) > 0$ is possible under (i).

Since servers are heterogeneous, we need to specify which idle server is chosen by the next customer to be served in FCFS manner. Furthermore, if more than one servers simultaneously complete service and if the number of waiting customers is less than the number of those servers, then we also need a rule to assign the waiting customers to available servers. This may not occur with a positive probability if service time distributions have densities. However, we do not make such an assumption, and give a server assigning rule to fully specify the model.

For this, we first put an arriving customer in the waiting line, and put the servers who complete service at the same moment into the set of available servers. We then select customers in the waiting line in FCFS manner, and exclusively assign them to servers from this set. Let $j$ be the number of waiting customers, and let $A$ be the index set of available servers. If $j \ge |A|$, then all servers becomes busy, so we do not need any rule, where $|A|$ is the cardinality of $A$. If $1 \le j < |A|$, we denote the index set of the servers who start new services by $D(A,j)$. For convenience, we let $D(A,j) = A$ if $j \ge |A|$. Thus, $D(A,j)$ is a set-valued random variable for each $A \in 2^{K}$ and $j \ge 1$. We assume that $D(A,j)$'s have the same distribution for each fixed $A$ and $j$ and are independent of everything else. We refer to this server selection according to $D(A,j)$ as a server selection rule. Since we consider the large queue asymptotics, the majority of our results are independent of this rule.

Since $\dd{P}(T_{0} = 0) > 0$ is possible, there is some arbitrariness when we apply the server selection rule. In this paper, we number arriving customers one by one, and independently apply the server selection rule for them even when they arrive simultaneously. Similarly, if some servers have zero service times, then the server selection rule is repeatedly applied every time when they occur. This dynamics is detailed in the proof of \lem{terminal 1} in  \app{terminal 1}.

We describe this queueing model by a continuous time Markov process. For this, we define the following characteristics at each time $t$, all of which are assumed to be right continuous.
\begin{align*}
 & L(t) \mbox{: the number of customers in system.}\\
 & I(t) \mbox{: the index set of idle servers.}\\
 & J(t) \equiv K \setminus I(t) \mbox{: the index set of busy servers.}\\
 & \ol{J(t)} \equiv \ol{K} \setminus I(t) \mbox{: the index $0$ is included}.\\
 & R_{0}(t) \mbox{: the residual exogenous arrival time.}\\
 & R_{i}(t) \mbox{: the residual service time at server $i \in J(t)$, and $R_{i}(t) = 0$ for $i \in I(t)$, }
\end{align*}
where $T_{i}$ is independently sampled subject to $F_{i}$ when $i \in \overline{J(t)}$ and $R_i(t-) = 0$. Here,
\begin{align*}
  R(t-) = \lim_{\epsilon \downarrow 0} R(t-\epsilon), \qquad t > 0.
\end{align*}

We denote the $(k+1)$-dimensional vector whose $i$-th entry is $R_{i}(t)$ by $\vc{R}(t)$. Let
\begin{align*}
  X(t) = (L(t), J(t), \vc{R}(t)), \qquad t \ge 0,
\end{align*}
then $X(\cdot) \equiv \{X(t); t \ge 0\}$ is a Markov process, which describes the heterogeneous $k$-server queue. By the right continuous assumption, $X(t)$ is right-continuous and has a left-hand limit for $t \ge 0$. Its state space is a subset of $\dd{Z}_{+} \times 2^{K} \times \dd{R}_{+}^{k+1}$, where $\dd{Z}_{+}$ and $\dd{R}_{+}$ are the sets of nonnegative integers and numbers, respectively, and $2^{K}$ is the set of all subsets of $K$.

This Markov process can be considered as a piecewise deterministic Markov process, PDMP for short, introduced by Davis \cite{Davi1984}. To have an appropriate state space for the Markov process, we first define faces which are specified by idle servers. For $A \in 2^{K}$, let
\begin{align*}
  S_{A} = \left\{
\begin{array}{ll}
 \{(\ell, K, \vc{y}) \in \dd{Z}_{+} \times 2^{K} \times \dd{R}_{+}^{k+1}; \ell \ge k\}, &  A = \emptyset,  \\
 \{(k - |A|, K \setminus A, \vc{y}) \in \dd{Z}_{+} \times 2^{K} \times \dd{R}_{+}^{k+1} \}, &  A \ne \emptyset.
\end{array}
\right.
\end{align*}
These $S_{A}$ are referred to as faces. Thus, the state space of PDMP $\{X(t)\}$ is given by $S = \cup_{A \in 2^{K}} S_{A}$. In some cases, it will be convenient to use other faces defined as
\begin{align*}
  \Gamma_{i} = \{(\ell, U, \vc{y}) \in S; i \in K \setminus U \}, \qquad i \in K.
\end{align*}
It is notable that $\Gamma_{i}$ is different from $S_{\{i\}}$. Yet another useful notation is the boundary of $S$ which represents the points where some residual times vanish. Similar to \cite{Davi1984}, we denote it by $\Gamma$. Namely,
\begin{align*}
  \Gamma = \big\{(\ell, U, \vc{y}) \in S; y_{0} = 0 \mbox{ or } K \setminus U \ne \emptyset \big\}.
\end{align*}
We will refer to $\Gamma$ as a terminal set of $S$ rather than a boundary.

In a PDMP, $\Gamma$ is the set of states just before discontinuous state changes. However, $\Gamma$ can not be identified by $L(t)$ and $\vc{R}(t)$. We here need their left-hand limits. This is the reason why we include $J(t)$ in $X(t)$. 

Throughout the paper, we assume the stability condition that
\begin{itemize}
\item [(ii)] The PDMP process $X(\cdot)$ has a stationary distribution, denoted by $\pi$. 
\end{itemize}
Recall $\rho = \lambda_{0}/(\sum_{i \in K} \lambda_{i})$. It is easy to see that (ii) implies $\rho \le 1$. On the other hand,
\begin{align}
\label{eqn:stability 1}
  \lambda_{0} < \sum_{i \in K} \lambda_{i}, \quad \mbox{ equivalently, } \quad \rho < 1,
\end{align}
is necessary for the stability of the fluid limit in Theorem 3.1 of \cite{Dai1995}, and therefore implies (ii) under the spread-out condition on either the inter-arrival time distribution or the service time distributions, where a distribution on $\dd{R}_{+}$ is said to be spread-out if its $n$ times convolution for some $n \ge 1$ has a absolute continuous component with respect to Lebesgue measure (see Section VII.1 of \cite{Asmu2003}). (ii) is unlikely satisfied for $\rho=1$ except for the case that the inter-arrival and service times are deterministic, but its verification is not known. So, we will just assume \eqn{stability 1} in addition to (i) and (ii).

\subsection{Extended generator of PDMP}
\label{sect:extended}

We first consider the dynamics of the Markov process $X(\cdot)$. Let $\sr{F}_{t} = \sigma(X(s); s \le t)$, where $\sigma(\cdot)$ stands for the minimal $\sigma$-field. $\{\sr{F}_{t}; t \ge 0\}$ is called a filtration. Then, $X(\cdot)$ is a strong Markov process with respect to $\{\sr{F}_{t}; t \ge 0\}$. We describe the dynamics of this Markov process using real valued test functions. For this, we introduce useful classes for them. By $C^{(1,p)}_{k+1}(S)$, we denote the set of measurable functions $f$ from $S$ to $\dd{R}$ such that, for each fixed $(\ell,U) \in \dd{Z}_{+} \times 2^{K}$,  $f(\ell, U, \vc{y})$ is continuous for $\vc{y} \in \dd{R}_{+}^{k+1}$, and has continuous partial derivatives for $\vc{y} \in \dd{R}_{+}^{k+1}$ except for $y_{i} \ne 0$. The set of bounded functions in $C^{(1,p)}_{k+1}(S)$ is denoted by $C^{(1,p)}_{b,k+1}(S)$.

For $f \in C^{(1,p)}_{k+1}(S)$, define a partial differential operator $\sr{A}$ as
\begin{align}
\label{eqn:generator 1}
  \sr{A} f(\ell, U, \vc{y}) = - \sum_{i \in \ol{U}} \frac {\partial f} {\partial y_{i}} (\ell, U, \vc{y}), \qquad (\ell, U, \vc{y}) \in S, y_{j} > 0, \forall j \in \ol{U},
\end{align}
where $\ol{U} = \{0\} \cup U$, and recall that $D_{i}$ is a finite subset of $\dd{R}_{+}$ such that the partial derivative of $f(\ell,U, \vc{y})$ with respect to $y_{i} > 0$.

For each $i \in \ol{K}$, let $T_i(\ell)$ be the $\ell$-th interarrival time for $i=0$ and service time for $i \in K$, let $t_{i,0} = 0$ and define $t_{i,n}$ as
\begin{align*}
  t_{i,n} = T_i(1) + T_i(2) + ... + T_i(n), \qquad n=1,2,\ldots,
\end{align*}
and define the counting process $N_{i}(\cdot) \equiv \{N_{i}(t); t \ge 0\}$ as
\begin{align*}
  N_{i}(t) = \sum_{n=1}^{\infty} B_{i,n} 1(t_{i,n} \le t), \qquad t \ge 0,
\end{align*}
where $1(\cdot)$ is the indicator function of proposition ``$\cdot$''. By superposing them, define the counting process $N(\cdot) \equiv \{N(t); t \ge 0\}$ as
\begin{align*}
  N(t) = N_{0}(t) + \sum_{i \in K} N_{i}(t), \qquad t \ge 0.
\end{align*}
By the definitions, $N(0) = N_{i}(0) = 0$. It also is easy to see that
\begin{align*}
  L(t) = L(0) + N_{0}(t) - \sum_{i \in K} N_{i}(t), \qquad t \ge 0.
\end{align*}

Note that $N(\cdot)$ may have jumps with size greater than 1. If this is the case, an arrival and/or multiple departures from different servers occur simultaneously. This situation is already argued as for the server selection rule. For our arguments, it is convenient to define the simplified version of $N(\cdot)$, which counts the number of atoms of $N(\cdot)$ but ignores their magnitude. To do so, we let $\{t_{n}; n \geq 1\}$ be a sequence of time instants defined as
\begin{align}
\label{eqn:jump instants 1}
  t_{n} = \inf\{s > t_{n-1}; N(s) - N(t_{n-1}) > 0 \}, \qquad n=1,2,\ldots,
\end{align}
where $t_{0} = 0$. We then define $N^{*}(\cdot) \equiv \{N^{*}(t); t \ge 0\}$ as
\begin{align}
\label{eqn:simplified N}
  N^{*}(t) = \sum_{n=1}^{\infty} 1(t_{n} \le t), \qquad t \ge 0.
\end{align}

Since $X(t)$ is absolutely continuous in $t \in \dd{R}_{+}$ with respect to the sum of the Lebesgue measure and the random measure generated by the counting process $N^{*}(\cdot)$, the integration formula yields
\begin{align}
\label{eqn:dynamics 1}
  f(X(t)) = f(X(0)) + \int_{0}^{t} \sr{A} f(X(s)) ds + \int_{0}^{t} \Delta f(X(s)) dN^{*}(s), \quad f \in C^{(1,p)}_{k+1}(S),
\end{align}
where $\Delta f(X(s)) = f(X(s)) - f(X(s-))$. Throughout the paper, operator $\Delta$ stands for such a difference. Note that $\Delta N^{*}(s) > 0$ if and only if $X(s-) \in \Gamma$, which causes a discontinuous state change, that is, $X(s-) \ne X(s)$. We refer to it as a jump. To describe this jump, we define the transition operator $Q$ on $C^{(1,p)}_{k+1}(S)$ as
\begin{align}
\label{eqn:Q 1}
  Q f \big(X(s-) \big) = \dd{E}\big( f(X(s)) | X(s-)\big), \qquad X(s-) \in \Gamma, f \in C^{(1,p)}_{k+1}(S),
\end{align}
where $X(0) , X(0-)$ should be considered as $X(0+)$ and $X(0)$, respectively. We refer to this $Q$ as a jump kernel. 

Our arguments will be based on the following fact. 
\begin{lemma}[A version of Proposition 4.3 of Davis \cite{Davi1984}]
\label{lem:martingale 1}
For a Borel-measurable function $f$ from $S$ to $\dd{R}$, let
\begin{align}
\label{eqn:martingale 1}
  M_{0}(t) = \int_{0}^{t} \big( f(X(s)) - Qf(X(s-)) \big) dN^{*}(s), \qquad t \ge 0,
\end{align}
then $M_{0}(\cdot) \equiv \{M_{0}(t); t \ge 0\}$ is  an $\sr{F}_{t}$-martingale if
\begin{align}
\label{eqn:finite 1}
  \dd{E}\Big(\int_{0}^{t} \big| f(X(s)) - Q f(X(s-)) \big| dN^{*}(s) \Big) < \infty, \qquad \mbox{for each } t > 0.
\end{align}

\end{lemma}

This fact is proved in \cite{Davi1984} using the compensator of the point process $N^{*}(\cdot)$, which is derived in \cite{Davi1976}. However, its proof can be much simpler, which is given in \app{martingale 1}. Let $M(t) = M_{0}(t) + f(X(0))$, then \eqn{dynamics 1} implies that, for $f \in C^{(1,p)}_{k+1}(S)$, 
\begin{align}
\label{eqn:martingale 2}
 M(t) = f(X(t)) - \Big(\int_{0}^{t} \sr{A} f(X(s)) ds + \int_{0}^{t} (Qf(X(s-)) - f(X(s-))) dN^{*}(s)\Big).
\end{align}
The following fact is immediate from \lem{martingale 1}.

\begin{lemma}[A special case of Theorem 5.5 of \cite{Davi1984}]
\label{lem:martingale 2}
For $f \in C^{(1,p)}_{k+1}(S)$, if either \eqn{finite 1} with $\dd{E}(|f(X(0))|) < \infty$ or $\dd{E}(|M(t)|) < \infty$ is satisfied, then $M(\cdot) \equiv \{M(t); t \ge 0\}$ is an $\sr{F}_{t}$-martingale. In particular, if $f$ satisfies that
\begin{align}
\label{eqn:terminal 1}
  Qf(\vc{x}) = f(\vc{x}), \qquad \forall \vc{x} \in \Gamma,
\end{align}
then the $\sr{F}_{t}$-martingale $M(\cdot)$ is simplified to
\begin{align}
\label{eqn:martingale 3}
  M(t) = f(X(t)) - \int_{0}^{t} \sr{A} f(X(s)) ds.
\end{align}
\end{lemma}

Davis \cite{Davi1984} refers to \eqn{terminal 1} as a boundary condition. Since the term boundary is used in a different sense in the queueing process context, we refer to \eqn{terminal 1} as a terminal condition. Recall that $\Gamma$ is named as the terminal set, which corresponds to this terminology. Ideally, one would want to find a class of test functions such that its elements satisfy \eqn{terminal 1} and characterize the Markov process $X(\cdot)$ (e.g., see Chapter 4 of \cite{EthiKurt1986}). However, no recipe is available for it. We attack this problem from a different viewpoint. Namely, we work with a smaller class of test functions whose elements satisfy \eqn{terminal 1}, while it does not characterize $X(\cdot)$ but is still useful for asymptotic analysis.

\subsection{Good test functions for the terminal condition}
\label{sect:good}

Let us consider the following form of a test function for each fixed $\theta, \eta, \zeta_{i} \in \dd{R}$.
\begin{align}
\label{eqn:test function 1}
  f(\ell,U,\vc{y}) = e^{\theta (\ell \vee k)} e^{\eta g_{0}(y_{0})} \prod_{i \in U} e^{\zeta_{i} g_{i}(y_{i})}, \qquad (\ell,U,\vc{y}) \in S,
\end{align}
where $g_{i}$ is a function from $\dd{R}_{+}$ to $\dd{R}_{+}$ such that $g_{i}(0) =0$, $g_{i}(y)$ is continuous and has a continuous derivative in $y$ except for finitely many values, and $a \vee b = \max(a,b)$ for real numbers $a$ and $b$. We also used the notation $\wedge$ as $a \wedge b = \min(a,b)$. 

We will use the following two functions for these $g_{i}$ for $i \in \ol{K}$.
\begin{align}
\label{eqn:test function g}
  g_{i}(y) = y, \qquad \mbox{or} \qquad g_{i}(y) = y \wedge v_{i}, \quad v_{i} > 0.
\end{align}
Obviously, $g_{i}(y) = y$ is a natural choice for $g_{i}$, but this choice may not work well when the tail probability of $T_{i}$ is heavy. This is the reason why we introduced the second choice in \eqn{test function g}, where $v_{i}$ is referred to as a truncation level. In what follows, a common $v$ is used for the truncation level for $i$ in a subset of $K$, while $u$ is used for $i =0$.

Our next task is to find $\eta, \zeta_{i}$ as functions of $\theta$ for the $f$ of \eqn{test function 1} to satisfy \eqn{terminal 1}. The idea behind this is to make the expected values of $\Delta f(X(t))$ equal zero at the jump instants of $N^{*}(\cdot)$ to vanish. For this, we choose them as
\begin{align}
\label{eqn:eta 1}
 & \mbox{$\eta$ is the solution of the equation: } e^{\theta} \dd{E}(e^{\eta (T_{0} \wedge u)}) = 1, \qquad u > 0,\\
\label{eqn:zeta 1}
 & \mbox{$\zeta_{i}$ is the solution of the equation: } e^{-\theta} \dd{E}(e^{\zeta_{i} (T_{i} \wedge v)}) = 1, \qquad v > 0,
\end{align}
and denote these $\eta$ and $\zeta_{i}$ by $\eta(u, \theta)$ and $\zeta_{i}(v, \theta)$, respectively. Clearly, these functions are well defined and finite for all $\theta$ in appropriate half lines or in $\dd{R}$ for finite $u, v$ because $\dd{E}(e^{x (T_{i} \wedge u)})$ is strictly increasing in $x \in \dd{R}$ for $i \in \ol{K}$. For simplicity, we denote $\eta(\infty, \theta)$ and $\zeta_{i}(\infty, \theta)$, respectively, by $\eta(\theta)$ and $\zeta_{i}(\theta)$, which may not exists for some $\theta \in \dd{R}$.

For $A \in 2^{K}$, we denote $f$ of \eqn{test function 1} with $g_{0}(y) = y \wedge u$, $g_{i}(y) = y$ for $i \in K \setminus A$, $g_{i}(y) = y \wedge v$ for $i \in A$ for each $u, v > 0$ and $\theta \in \dd{R}$ by $f_{u,A(v),\theta}$. Namely,
\begin{align}
\label{eqn:test function 3}
  f_{u,A(v),\theta}(\ell,U,\vc{y}) = e^{\theta (\ell \vee k) + \eta(u,\theta) (y_0 \wedge u) + \sum_{i \in U \setminus A}\zeta_{i}(\theta) y_{i} + \sum_{i \in A \cap U}\zeta_{i}(v,\theta) (y_{i} \wedge v)}, 
\end{align}
for $(\ell,U,\vc{y}) \in S$. If $u=\infty$, this test function is denoted by $f_{A(v),\theta}$. If $A = \emptyset$, then it is denoted by $f_{u,\theta}$. If $u = \infty$ and $A = \emptyset$, then it is denoted by $f_{\theta}$. For example,
\begin{align*}
  f_{K(v),\theta}(\ell,U,\vc{y}) = e^{\theta (\ell \vee k) + \eta(\theta) y_0 + \sum_{i \in U}\zeta_{i}(v,\theta) (y_{i} \wedge v)}.
\end{align*}
This notation system will be used throughout the paper. The next lemma  will be proved in \app{terminal 1}.

\begin{lemma}
\label{lem:terminal 1}
Under the setting \eqn{eta 1} and \eqn{zeta 1} and assumption (i), the terminal condition \eqn{terminal 1} holds for $f = f_{u,A(v),\theta}$ for finite $u, v > 0$ and $A \subset K$ as long as $\zeta_{i}(\theta)$ are well defined for $i \in K \setminus A$ and $\eta(u,\theta)$ is well defined when $u=\infty$.
\end{lemma}

We will note that $\eta(u,\theta)$ and $\zeta_{i}(v,\theta)$ have nice probabilistic interpretations. Because they are similar, we present their general properties using some new notation. Let $T$ be a nonnegative random variable with a finite and positive expectation. Denote its distribution by $F$, and let $\widehat{F}(\theta) = \dd{E}(e^{\theta T})$ for $\theta \in \dd{R}$. We define $\beta_{*}$ and $\theta_{*}$ as
\begin{align}
\label{eqn:limit point 1}
  \beta_{*} = \sup \{\theta \in \dd{R}; \widehat{F}(\theta) < \infty\}, \qquad \theta_{*} = - \sup\{ \theta \in \dd{R}; e^{\theta} <  \widehat{F}(\beta_{*}) \}.
\end{align}
Note that $\beta_{*} \ge 0 \ge \theta_{*}$. If $\beta_{*} = 0$, then the distribution $F$ is said to have a heavy tail. Otherwise, it is said to have a light tail. For $v > 0$, we consider $T\wedge v$ instead of $T$, and define $F(v,\cdot)$, $\widehat{F}(v,\cdot)$, $\theta_{*}(v)$ and $\beta_{*}(v)$ in a similar manner. It is notable that
\begin{align*}
  \beta_{*}(v) = +\infty, \qquad \theta_{*}(v) = - \infty,
\end{align*}
because $T \wedge v$ is bounded by $v$. Let $\ol{\theta} = - \log F(0)$, where $\ol{\theta} = \infty$ if $F(0)=0$. Similar to \eqn{eta 1}, we define the functions $\xi(v,\cdot)$ as the solution to
\begin{align}
\label{eqn:xi 1}
  e^{\theta} \widehat{F}(v,\xi(v,\theta)) = 1.
\end{align}
Then, for finite $v > 0$, $\xi(v,\theta)$ is well defined and finite for all $\theta < \ol{\theta}$, but this may not be the case for $\xi(\theta) \equiv \xi(\infty,\theta)$. 

We note that $-\xi(v,\cdot)$ can be considered as a rate function for large deviations. To see this, let $\{N(v,t); t \ge 0 \}$ be the renewal counting process with interval time distribution $F(v,\cdot)$. From \eqn{xi 1} and Theorem 1 of Glynn and Whitt \cite{GlynWhit1994}, we have
\begin{align}
\label{eqn:zeta 3}
 & \xi(v,\theta) = - \lim_{t \to \infty} \frac 1t \log \dd{E}\big( e^{\theta N(v,t)} \big), \qquad \theta < \ol{\theta}.
\end{align}
Intuitively, this fact is helpful to understand what the function $\xi(v,\cdot)$ is, but technically it is not so much helpful. Because of this, we will not use the fact that \eqn{zeta 3} holds in proofs. We now present important properties of $\xi(v,\cdot)$.

\begin{lemma}
\label{lem:concave 1}
(a) For each finite $v > 0$, $\xi(v,\theta)$ is decreasing, concave and infinitely differentiable in $\theta < \ol{\theta}$. (b) For each fixed $\theta \in (0,\ol{\theta})$ ($\theta < 0$, respectively), $\xi(v,\theta)$ is negative (positive) and increasing (decreasing) in $v$. (c) For each $\delta \in (0, \ol{\theta})$ and $v > 0$, we have
\begin{align}
\label{eqn:xi bound 1}
  |\xi(v,\theta)| \le \max\Big(\frac {1} {\dd{E}(T \wedge v)}, \frac 1\delta |\xi(v,\delta)|\Big) |\theta|, \qquad \theta < \delta.
\end{align}
\end{lemma}

Although (a) of this lemma is immediate from \eqn{zeta 3}, we prove it in \app{concave 1} as mentioned above. Using \lem{concave 1}, we define ${\xi}(\vartriangle,\cdot)$ as
\begin{align}
\label{eqn:xi infty}
  {\xi}(\vartriangle,\theta) = \lim_{v \uparrow \infty} \xi(v,\theta), \qquad \theta < \ol{\theta}.
\end{align}
It is notable that ${\xi}(\vartriangle,\theta)$ may be different from $\xi(\theta) \equiv \xi(\infty,\theta)$ as shown in the lemma below. From \lem{concave 1}, we immediately have the following facts.  See \supp{S1} on its proof.

\begin{lemma}
\label{lem:concave 2}
(a) ${\xi}(\vartriangle,\theta)$ is finite and concave for all $\theta < \ol{\theta}$. (b) ${\xi}(\vartriangle,\theta) \le \beta_{*}$ for all $\theta < \ol{\theta}$ with equality for $\theta \le \theta_{*}$, and ${\xi}(\vartriangle,\theta) = \xi(\theta)$ for $\theta \in (\theta_{*},\ol{\theta})$. (c) $\widehat{F}({\xi}(\vartriangle,\theta)) < \infty$ for all $\theta < \ol{\theta}$.
\end{lemma}

\begin{remark}
\label{rem:concave 2}
Despite (c), the derivative $\widehat{F}'({\xi}(\vartriangle,\theta)) = \dd{E}(T_{i} e^{{\xi}(\vartriangle,\theta) T_{i}})$ may not exist at $\theta = \theta_{*} > -\infty$, namely at $\theta < 0$ satisfying ${\xi}(\vartriangle,\theta) = \beta_{*}$.  A typical example for this is that $F(x) = 1 - (1+x)^{-r} e^{- \delta x}$ with $1 < r \le 2$ and $\delta \ge 0$, for which $\beta_{*} = \delta$, $\widehat{F}(\delta) < \infty$ and $\theta_{*} = - \log \widehat{F}(\delta) < 0$, and therefore $\xi(\vartriangle,\theta_{*}) = \delta$, but $\widehat{F}'(\delta) = \infty$. \end{remark}

\begin{figure}[h] 
   \centering
   \includegraphics[height=3cm]{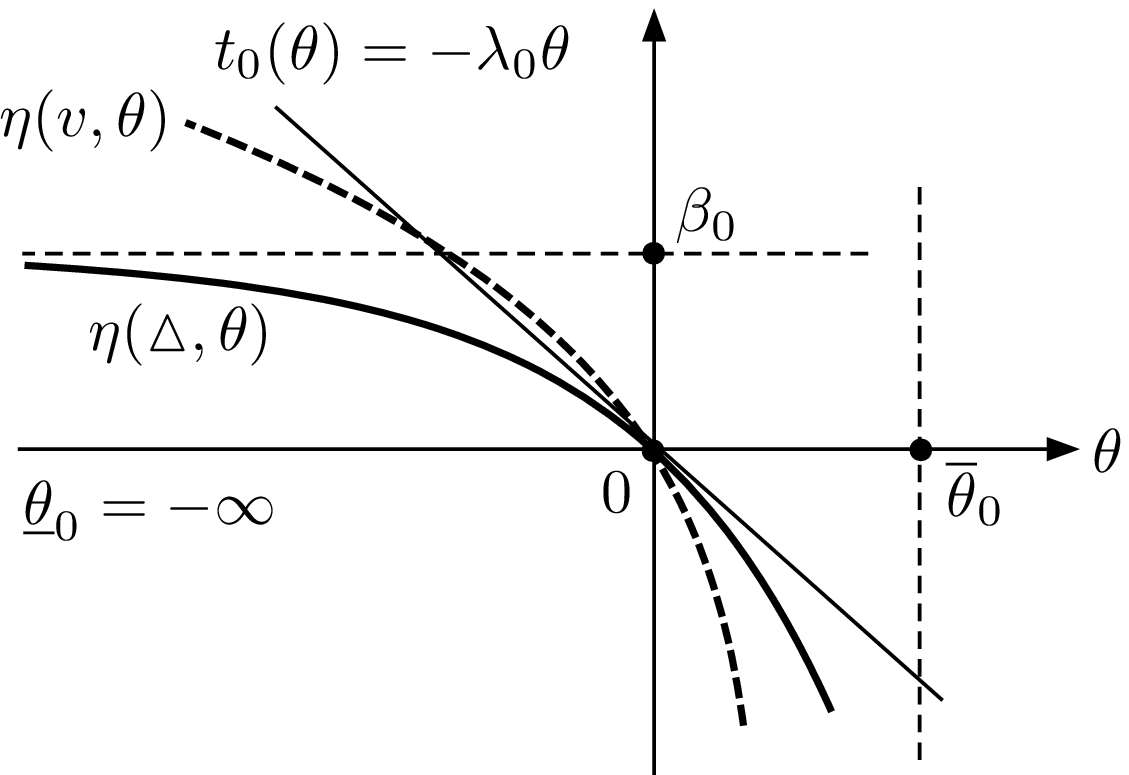} \hspace{2ex}
   \includegraphics[height=3cm]{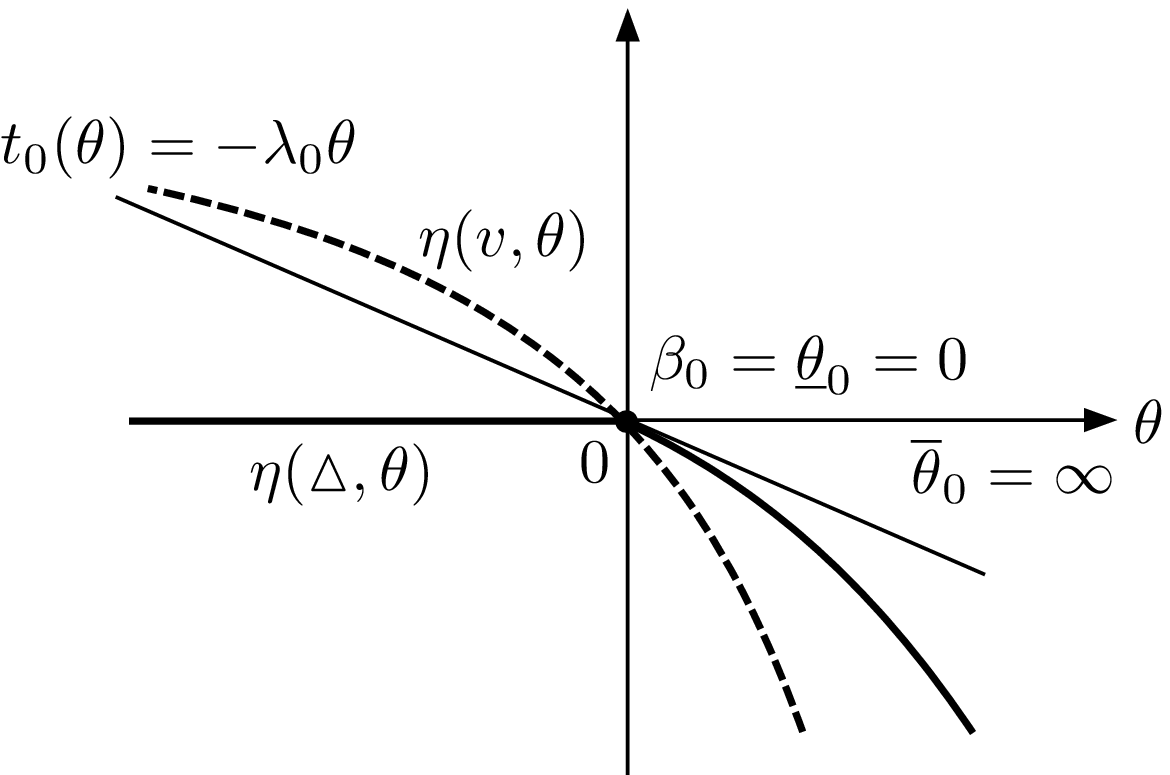} 
   \caption{$\eta(v,\theta)$ and ${\eta}(\vartriangle,\theta)$ for light and heavy (left and right panels) tailed $F_{0}$}
   \label{fig:shape of eta}
\end{figure}

\begin{figure}[h] 
   \centering
   \includegraphics[height=3cm]{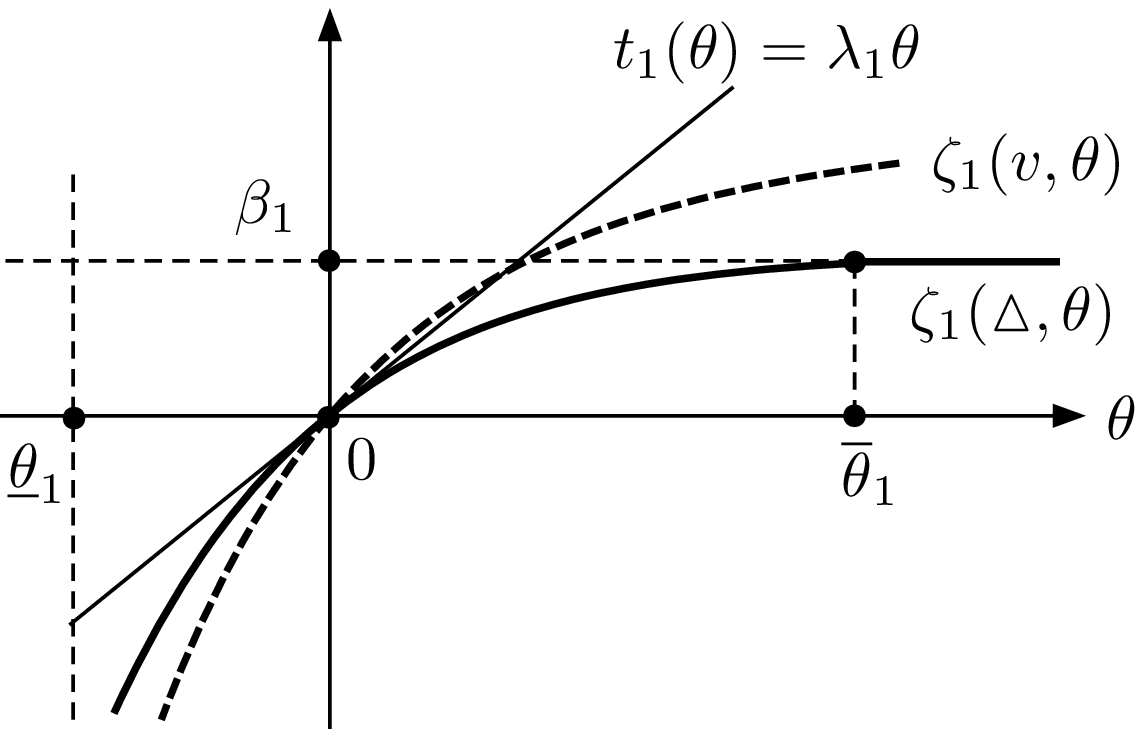} \hspace{2ex}
   \includegraphics[height=3cm]{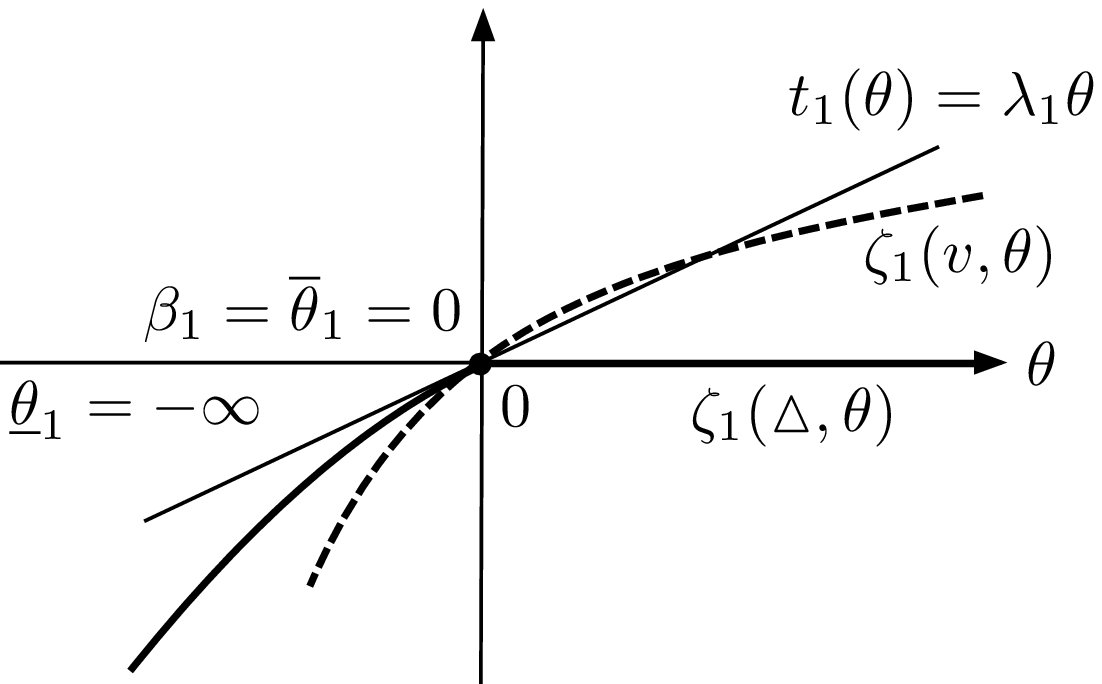} 
   \caption{$\zeta_{1}(v,\theta)$ and $\zeta_{1}(\vartriangle,\theta)$ for light and heavy (left and right panels) tailed $F_{1}$}
   \label{fig:shape of xi}
\end{figure}

We next apply these properties for $T_{0}$ and $T_{i}$ for $i \in K$. For $T = T_{0}$, we obviously have that $\eta = \xi$ and $\eta(v,\cdot) = \xi(v,\cdot)$, and therefore ${\eta}(\vartriangle,\cdot)$ is defined as ${\xi}(\vartriangle,\cdot)$. Similarly, for $T = T_{i}$ with $i \in K$, we have $\zeta_{i}(\theta) = \xi(-\theta)$ and $\zeta_{i}(v,\theta) = \xi(v,-\theta)$ by \eqn{zeta 1}, and ${\zeta}_{i}(\vartriangle,\theta)$ is defined as ${\xi}(\vartriangle,-\theta)$. Similar to $\theta_{*}$ and $\beta_{*}$ of \eqn{limit point 1}, define, for $i \in K$,
\begin{align}
\label{eqn:limit point 2}
  \ul{\theta}_{i} = \log F_{i}(0), \quad \beta_{i} = \sup \{ \theta \in \dd{R}; \widehat{F}_{i}(\theta) < \infty\}, \quad \theta_{i} = \sup \{\theta \in \dd{R}; e^{\theta} < \widehat{F}_{i}(\beta_{i})\},
\end{align}
while $\ol{\theta}_{0} = -\log F_{0}(0)$, $\beta_{0} = \beta_{*}$ and $\theta_{0} = \theta_{*}$ for $T=T_{0}$, where $\widehat{F}_{i}$ is the moment generating function of the distribution $F_{i}$, namely, $\widehat{F}_{i}(\theta) = \dd{E}(e^{\theta T_{i}})$, and, for $i \in \ol{K}$, $\log F_{i}(0) = -\infty$ if $F_{i}(0)=0$. Based on Lemmas \lemt{concave 1} and \lemt{concave 2}, typical shapes of the functions $\eta(v,\cdot)$, $\eta(\vartriangle,\cdot)$, $\zeta(v,\cdot)$ and $\zeta(\vartriangle,\cdot)$ are drawn in Figures \figt{shape of eta} and \figt{shape of xi}.

\subsection{Stationary equations}
\label{sect:stationary}

Let $u=v$ and $A = K$ for $f_{u,A(v),\theta}$ of \eqn{test function 3}. We will use this test function $f_{v,K(v),\theta}$ for deriving the stationary equations. Let us introduce the following notation.
\begin{align*}
 & \Phi_{i}(v,\theta) = \dd{E}_{\pi}(f_{v,K(v),\theta}(X(0)) 1(R_{i}(0) < v)), \qquad i \in \ol{K},\\
 & \Phi_{i,0}(v,\theta) = \dd{E}_{\pi}(f_{v,K(v),\theta}(X(0))1(i \not\in J(0))), \qquad i \in K,
\end{align*}
where $\dd{E}_{\pi}$ stands for the expectation conditioned on $X(0)$ subject to $\pi$. We define these functions also for $v \uparrow \infty$, and denote the limit of $\Phi_{i}(v,\theta)$ by $\Phi(\vartriangle,\theta)$, because it does not depend on $i$, while $\Phi_{i,0}(\vartriangle,\theta)$ is used as is. We wish to apply \lem{martingale 2} to the function $f = f_{v,K(v),\theta}$ under the stationary distribution $\pi$. From the definition of $\sr{A}$ in \eqn{generator 1}, we have
\begin{align}
\label{eqn:Af 1}
 \sr{A} f_{v,K(v),\theta} & (\ell,U,\vc{y}) = - \Big(\eta(v,\theta) 1(y_{0} < v) + \sum_{i \in U} \zeta_{i}(v,\theta) 1(y_{i} < v) \Big) f_{v,K(v),\theta}(\ell,U,\vc{y}) \nonumber\\
  &= - \Big( \eta(v,\theta) 1(y_{0} < v) + \sum_{i \in K} \zeta_{i}(v,\theta) 1(y_{i} < v) \nonumber\\
  & \hspace{22ex} - \sum_{i \in K \setminus U} \zeta_{i}(v,\theta) 1(y_{i} < v) \Big) f_{v,K(v),\theta}(\ell,U,\vc{y}).
\end{align}
Taking the expectation of \eqn{martingale 3} in \lem{martingale 2} for $f = f_{v,K(v),\theta}$ and assuming that $\dd{E}_{\pi}(f(X(0)))$ is finite, we have the following lemma since $\Phi_{i}(v,\theta) < \infty$ implies that $\Phi_{i,0}(v,\theta)$ are finite for all $i \in \ol{K}$.

\begin{lemma}
\label{lem:SE 1}
Under assumptions (i) and (ii), for $v > 0$ and $\theta \in (\ul{\theta}_{K}, \ol{\theta}_{0})$, where $\ul{\theta}_{K} = \vee_{i \in K} \ul{\theta}_{i}$,
\begin{align}
\label{eqn:SE 1}
  \eta(v,\theta) \Phi_{0}(v,\theta) + \sum_{i \in K} \zeta_{i}(v,\theta) \Phi_{i}(v,\theta) - \sum_{i \in K} \zeta_{i}(v,\theta) \Phi_{i,0}(v,\theta) = 0,
\end{align}
as long as $\Phi_{i}(v,\theta)$ is finite for all $i \in \ol{K}$. In particular, \eqn{SE 1} always holds for $\theta \le 0$.
\end{lemma}

Let $v$ go to infinity in \eqn{SE 1}, then Lemmas \lemt{concave 1} and \lemt{concave 2} yield the following lemma, recalling that $\Phi(\vartriangle,\theta) = \Phi_{i}(\vartriangle,\theta)$.

\begin{lemma}
\label{lem:SE 2}
Under the same assumptions as \lem{SE 1},
\begin{align}
\label{eqn:SE 2}
  \Big({\eta}(\vartriangle,\theta) + \sum_{i \in K} {\zeta}_{i}(\vartriangle,\theta) \Big) \Phi(\vartriangle \theta) - \sum_{i \in K} {\zeta}_{i}(\vartriangle,\theta) \Phi_{i,0}(\vartriangle,\theta) = 0,
\end{align}
as long as $\Phi(\vartriangle,\theta)$ is finite.
\end{lemma}

We may consider \eqn{SE 1} and \eqn{SE 2} as stationary equations. For studying the asymptotic problems related to the stationary distribution of $L(\cdot) \equiv \{L(t); t \ge 0\}$, one may consider to apply the techniques developed in \cite{KobaMiya2014,Miya2015a}, which are based on the stationary equation similar to \eqn{SE 2} and its variant called a stationary inequality. However, one must be careful about it. First, $\Phi(\vartriangle,\theta)$ is not the moment generating function of $L$. Second, ${\eta}(\vartriangle, z) + \sum_{i \in K} {\zeta}_{i}(\vartriangle, z)$ may not be analytically extendable from the complex region $\{ z \in \dd{C}; \Re z < 0\}$, since it is difficult to see their analytic properties. Therefore it is hard to apply the method of analytic extension in general, which is typically used for tail asymptotic problems (see, e.g., \cite{KobaMiya2012}).

\section{Large queue asymptotics}
\label{sect:answers}
\setnewcounter

In this section, we present main results on tail asymptotics and weak limit approximations for the stationary distribution of $L(\cdot)$.

\subsection{Tail asymptotics}
\label{sect:tail}

We first consider the tail asymptotics. For simplicity, we assume that $F_{0}(0)=0$ in this subsection, so $\ol{\theta}_{0} = \infty$. As we argued at the end of the previous section, \lem{SE 2} may not be helpful for this purpose. Nevertheless, the coefficient of $\Phi(\theta)$ in it suggests that the tail decay rate would be given by
\begin{align}
\label{eqn:alpha 1}
  \alpha \equiv \sup\Big\{\theta \ge 0; - \Big(\eta(\theta) + \sum_{i \in K} {\zeta}_{i}(\vartriangle,\theta) \Big) \le 0 \Big\},
\end{align}
where we have replaced ${\eta}(\vartriangle,\theta)$ by $\eta(\theta)$ using \lem{concave 2}.  Recall the definitions of $\beta_{i}$ and $\theta_{i}$ for the distribution $F_{i}$ (see \eqn{limit point 2}). If $\alpha < \theta_{i}$, \lem{concave 2} implies that we can also replace ${\zeta}_{i}(\vartriangle,\theta)$ in \eqn{alpha 1} by $\zeta_{i}(\theta)$. Let
\begin{align}
\label{eqn:gamma 1}
  \gamma_{K(\vartriangle)}(\theta) = - \Big(\eta(\theta) + \sum_{i \in K} {\zeta}_{i}(\vartriangle,\theta) \Big),
\end{align}
then $\gamma_{K(\vartriangle)}(\theta)$ is finite for $\theta \ge 0$ and a convex function by \lem{concave 2}. This $\gamma_{K(\vartriangle)}$ corresponds to the rate function corresponding to $-\xi(v,\theta)$ of \eqn{zeta 3}, which is the reason why minus signs are in \eqn{alpha 1} and \eqn{gamma 1}. If $\alpha > 0$, then $\gamma_{K(\vartriangle)}(\theta) < 0$ for $\theta \in (0,\alpha)$ since $\gamma'_{K(\vartriangle)}(0) = \lambda_{0} - \sum_{i \in K} \lambda_{i} < 0$ by \eqn{stability 1}. Furthermore, if ${\zeta}_{i}(\vartriangle,\alpha) < \beta_{i}$, equivalently, $\alpha < \theta_{i}$, for all $i \in K$, then $0 < \gamma'_{K(\vartriangle)}(\alpha) < \infty$. Otherwise, $\gamma'_{K(\vartriangle)}(\alpha)$ may diverge by \rem{concave 2}.

We prove the following theorem in \sectn{proofs}. It shows that $\alpha$ is indeed the decay rate in the logarithmic sense. For this, we need a further condition. 
\begin{itemize}
\item [(iii)] For all $i \in K$, the distribution of $T_{i}$ is spread-out.
\end{itemize}

\begin{theorem}
\label{thr:tail asymptotic 1}
Assume that the PDMP $\{(L(t), J(t), \vc{R}(t))\}$ for the heterogeneous $k$ server queue satisfies (i) with $F_{0}(0)=0$, (ii), (iii) and \eqn{stability 1}, and let $L$ be a random variable subject to the stationary distribution of $L(\cdot)$. Then, (a) If $\alpha < \theta_{i}$ for all $i \in K$, where $\theta_{i}$ is defined by \eqn{limit point 2}, then $\alpha > 0$ and, for some $c > 0$,
\begin{align}
\label{eqn:exact asymptotic 1}
  \lim_{x \to \infty} e^{-\alpha x} \dd{P}(L > x) = c.
\end{align}
(b) For each $\alpha \ge 0$, we have
\begin{align}
\label{eqn:decay rate 1}
  \lim_{x \to \infty} \frac 1x \log \dd{P}(L > x) = - \alpha.
\end{align}
\end{theorem}

\begin{remark}
\label{rem:tail asymptotic 2}
For a heterogenous multiserver queue with batch arrivals, similar tail asymptotics are obtained for the stationary distributions of the waiting time, the queue length just after arrival instants and some related characteristics embedded at arrival instants in \cite{Sado1995,SadoSzpa1995}. Thus, \thr{tail asymptotic 1} may be considered as the continuous time counterparts of those discrete time results for single arrivals. Its batch arrival version will be discussed in \sectn{on arrival}, which includes the case that $F_{0}(0) > 0$. 
\end{remark}

From Figures \figt{shape of eta} and \figt{shape of xi}, we can see that the heaviness of the inter-arrival time distribution $F_{0}$ has no significant influence on the  tail decay rate $\alpha$. This is intuitively clear. Let us consider how the heaviness of the service time distributions influences the tail decay rate. Sadowsky \cite{Sado1995} calls a server with a heavy tailed service time distribution lazy, and remarks that they have no contribution for the decay rate $\alpha$ (see the second last paragraph in Section 1 in \cite{Sado1995}). In this paper, such a server is called slow. To see the influence of slow servers, define, for $A \in 2^{K}$, 
\begin{align*}
 & \rho_{A} = \frac {\lambda_{0}} {\sum_{i \in A} \lambda_{i}},\\
 & \alpha_{A} = \sup\Big\{\theta \ge 0; - \Big(\eta(\theta) + \sum_{i \in A} {\zeta}_{i}(\vartriangle,\theta) \Big) \le 0 \Big\}.
\end{align*}

Let $K_{0} = \{i \in K; \beta_{i} = 0\}$, which is the index set of the heavy tailed service time distributions. Then, the following fact is immediate from \thr{tail asymptotic 1} since $\zeta_{i}(\vartriangle,\alpha) = 0$ for $i \in K_{0}$ and $\gamma_{K(\vartriangle)}(\theta)$ is convex in $\theta \in \dd{R}$.

\begin{corollary}
\label{cor:slow server}
 Under the assumptions of \thr{tail asymptotic 1}, $\alpha = \alpha_{K \setminus K_{0}}$, and $\alpha > 0$ if and only if $\rho_{K \setminus K_{0}} < 1$.
\end{corollary}

As we mentioned in \sectn{introduction}, our approach is flexible for modeling extensions. As an example, we here modify the arrival process. The arrival process has been assumed to be a renewal process for the the heterogeneous multiserver queue. We replace it by the superposition of renewal processes.  Assume that there are $m$ arrival streams, indexed by $1,2,\ldots, m$, which are renewal processes with the interval time distributions $F_{e,1}, F_{e,2}, \ldots, F_{e,m}$, respectively. We further assume that they are mutually independent and independent of everything else. Denote those renewal processes as counting processes by $N_{e,1}(\cdot), N_{e,2}(\cdot), \ldots, N_{e,m}(\cdot)$, then their superposition is given by
\begin{align*}
  N_{e}(t) \equiv \sum_{i=1}^{m} N_{e.i}(t), \qquad t \ge 0.
\end{align*}
We take this $N_{e}(\cdot) \equiv \{N_{e}(t); t \ge 0\}$ as an arrival process for the heterogeneous multiserver queue. We assume the same servers as in \sectn{hetero}.

In general, superposing renewal processes is analytically not tractable for an arriving process of a queueing model. However, this is not the case in the present formulation because we only need to redefine $\eta(\theta)$ as
\begin{align}
\label{eqn:eta e1}
  \eta(\theta) = \sum_{i=1}^{m} \eta_{i}(\theta), \qquad \theta \in \dd{R},
\end{align}
as long as all $\eta_{i}(\theta)$ exists, where $\eta_{i}(\theta)$ is obtained as the solution to
\begin{align}
\label{eqn:eta e2}
  e^{\theta} \widehat{F}_{e,i} (\eta_{i}(\theta)) = 1, \qquad i =1,2,\ldots,m.
\end{align}
Note that $\eta_{i}(\theta) \le 0$ for $\theta \ge 0$, and therefore it exists and is finite for $\theta \ge 0$.

\thr{tail asymptotic 1} can be extended for the superposed renewal arrivals. This extension is routine, so we here only present results without proof. Let $R_{e,i}(t)$ for $i \in J$ be the residual arrival time at time $t \ge 0$ for the $i$-th renewal arrival process, and let $R_{s,j}(t)$ be the residual service time at time $t \ge 0$ for the $j$-th server. Denote $\vc{R}_{e}(t)$ and $\vc{R}(t)$ be the vectors whose $i$-th and $j$-th entries are $R_{e,i}(t)$ and $R_{s,i}(t)$, respectively.

\begin{theorem}
\label{thr:tail asymptotic 2}
Assume that the PDMP $\{(L(t), J(t), \vc{R}_{e}(t), \vc{R}(t))\}$ for the heterogeneous $k$ server queue with superposed renewal arrivals satisfies (i) with $F_{e,i}(0)=0$ for $i=1,2,\ldots,m$, (ii), (iii) and \eqn{stability 1}, where $\lambda$ is the total arrival rate, and let $L$ be a random variable subject to the stationary distribution of $L(\cdot)$. Then, (a) and (b) of \thr{tail asymptotic 1} hold true for the $\eta$ of \eqn{eta e1}.
\end{theorem}

\subsection{Weak limit approximations for a large queue}
\label{sect:weak limit}

We next consider weak limit approximations for a sequence of the heterogeneous multiserver queues, where $F_{0}(0)$ may not necessarily vanish. For this queue, Chen and Yao \cite{ChenYe2011a} study a diffusion approximation, that is, a process limit of those queues under diffusion scaling in heavy traffic. They use a bounding method, called sandwich. As mentioned in \sectn{introduction}, we do not deal with process limits, and directly consider the stationary queue length distribution. In this situation, we do not need to scale time. We consider two weak limit approximations for bounded and large variances under the heavy traffic condition.

We consider a sequence of systems indexed by $n = 1,2,\ldots$, and denote the characteristics of the $n$-th system by superscript $^{(n)}$ (for example, $X^{(n)}(t)$, $R_{i}^{(n)}(t)$ and $T^{(n)}_{i}$ for $i \in \ol{K}$). Throughout the indexed systems, we assume that the number $k$ of servers is unchanged, while the server selection rule may change. The mean arrival rate and the mean service rate of server $i$ are denoted by $\lambda_{0}^{(n)}$ and $\lambda_{i}^{(n)}$, respectively. Denote the variances of $T^{(n)}_{i}$ by $(\sigma^{(n)}_{i})^{2}$, which are assumed to be finite. Let
\begin{align*}
  \rho^{(n)} = \lambda_{0}^{(n)} \Big(\sum_{i \in K} \lambda_{i}^{(n)} \Big)^{-1},
\end{align*}
and let $L^{(n)}$ be a random variable subject to the stationary distribution of $L^{(n)}(\cdot) \equiv \{L^{(n)}(t); t \ge 0 \}$. 

We will use the stationary equations \eqn{SE 1} and \eqn{SE 2} for $\theta \le 0$ to study the weak limit of the stationary distribution of $L^{(n)}(\cdot)$ under a suitable scaling. Let $\eta^{(n)}(v,\theta)$ and $\zeta^{(n)}_{i}(v,\theta)$ be the truncated versions of $\eta^{(n)}(\theta)$ and $\zeta^{(n)}_{i}(\theta)$, respectively, for the $n$-th system, where $v$ is a truncation level. Of importance to us will be their Taylor expansions around the origin. Let $\lambda^{(n)}_{i}(v) = 1/\dd{E}(T^{(n)}_{i} \wedge v)$, and let $(\sigma_{i}^{(n)}(v))^{2}$ be the variance of $T^{(n)}_{i} \wedge v$ for $i \in \ol{K}$. For $v > 0$, we define the functions,
\begin{align}
\label{eqn:Taylor 1a}
 & \epsilon^{(n)}_{0}(v,\theta) = \eta^{(n)}(v,\theta) + \lambda^{(n)}_{0}(v) \theta + \frac 12 (\lambda^{(n)}_{0}(v))^{3} (\sigma_{0}^{(n)}(v))^{2} \theta^{2}, \quad \theta < \ol{\theta}^{(n)}_{0},\\
\label{eqn:Taylor 1b}
 & \epsilon^{(n)}_{i}(v,\theta) = \zeta_{i}^{(n)}(v,\theta) - \lambda^{(n)}_{i}(v) \theta + \frac 12 (\lambda^{(n)}_{i}(v))^{3} (\sigma_{i}^{(n)}(v))^{2} \theta^{2}, \quad \theta > \ul{\theta}^{(n)}_{i}, i \in K,
\end{align}
which correspond to the error terms for the 2nd order Taylor expansions of $\eta^{(n)}(v,\theta)$ and $\zeta^{(n)}_{i}(v,\theta)$. We define these error functions also for $v=\vartriangle$ similarly to the definition \eqn{xi infty}, in which $\eta^{(n)}(\vartriangle,\theta) = \lim_{v \uparrow \infty} \eta^{(n)}(v,\theta)$, $\zeta^{(n)}_{i}(\vartriangle,\theta) = \lim_{v \uparrow \infty} \zeta^{(n)}_{i}(v,\theta)$, $\lambda^{(n)}_{i}(\vartriangle) = \lambda^{(n)}_{i}(\infty) = \lambda^{(n)}_{i}$ and $\sigma^{(n)}_{i}(\vartriangle) = \sigma^{(n)}_{i}(\infty) = \sigma^{(n)}_{i}$.

We will study the weak convergence of $q_{n} L^{(n)}$ for a sequence of positive numbers $\{q_{n}; n=1,2,\ldots\}$ which vanishes as $n$ goes to infinity. To this end, we will evaluate $\epsilon^{(n)}_{i}(1/q_{n}, q_{n} \theta)$. We first assume the following mean rate condition:
\begin{align}
\label{eqn:moment 1}
  \lim_{n \to \infty} \lambda^{(n)}_{i} = \lambda_{i}, \qquad i \in \ol{K},
\end{align}
for some finite $\lambda_{i} > 0$. Clearly, this condition is equivalent to the convergence of the first moments of $T^{(n)}_{i}$ to finite and positive constants. We need the following fact, which is essentially a special case of Lemma 4.4 of \cite{BravDaiMiya2015}, but an extra condition is assumed in \cite{BravDaiMiya2015}. We show that it is not needed in \app{moment 1}.
\begin{lemma}
\label{lem:moment 1}
If \eqn{moment 1} and $\{T_{i}^{(n)}; n \ge 1\}$ is uniformly integrable, then, for a sequence $q_{n} > 0$ converging to $0$, we have, for some positive constants $c_{i}(\delta)$ for each $\delta > 0$,
\begin{align}
\label{eqn:bound 2}
 & \limsup_{n \to \infty} \Big| \max\Big( \frac {\eta^{(n)}(\vartriangle, q_{n} \theta)} {q_{n}}, \frac {\eta^{(n)}(1/q_{n}, q_{n} \theta)} {q_{n}} \Big) \Big| \le c_{0}(\delta) |\theta|, \qquad |\theta| < \delta,\\
\label{eqn:bound 3}
 & \limsup_{n \to \infty} \Big| \max\Big( \frac {\zeta^{(n)}_{i}(\vartriangle, q_{n} \theta)} {q_{n}}, \frac {\zeta^{(n)}_{i}(1/q_{n}, q_{n} \theta)} {q_{n}} \Big) \Big| \le c_{i}(\delta) |\theta|, \qquad i \in K, |\theta| \le \delta.
\end{align}
\end{lemma}
Up to now, we have only assumed the convergence of the first moments, and $q_{n} > 0$ can be arbitrary as long as it converges to $0$. We are ready to consider two scenarios on the variances of inter-arrival and service time distributions for this sequence. 

We first consider the heavy traffic approximation for bounded variances. That is, in addition to \eqn{moment 1}, we assume that there are nonnegative constants $\sigma_{i}$ for $i \in \ol{K}$ such that
\begin{align}
\label{eqn:heavy 1}
 & 0 < \rho^{(n)} < 1 \; \mbox{ for all } n \ge 1, \qquad \lim_{n \to \infty} \rho^{(n)} = 1,\\
\label{eqn:heavy 2}
 & \lim_{n \to \infty} \sigma^{(n)}_{i} = \sigma_{i} \ge 0, \quad i \in \ol{K}, \qquad \sum_{i \in \ol{K}} \sigma_{i} > 0,\\
\label{eqn:heavy 3}
 & \lim_{a \to \infty} \sup_{n \ge 1} \dd{E}\big((T_{i}^{(n)})^{2} 1(T_{i}^{(n)} > a)\big) = 0, \qquad i \in \ol{K},
\end{align}
where $(\sigma_{i}^{(n)})^{2}$ is the variance of $T^{(n)}_{i}$. The conditions \eqn{heavy 1} and \eqn{heavy 3} are so called heavy traffic and uniform integrability conditions, respectively. Note that $\lambda_{i}$'s in \eqn{moment 1} must be positive under \eqn{heavy 1}. Let $r_{n} = 1 - \rho^{(n)}$ and choose $q_{n} = r_{n}$. The next lemma is a key for our proof, which is proved in \app{Taylor 1}.

\begin{lemma}
\label{lem:Taylor 1}
Assume \eqn{moment 1} and the heavy traffic conditions \eqn{heavy 1}--\eqn{heavy 3}, then, for each $i \in \ol{K}$ and each $a > 0$, 
\begin{align}
\label{eqn:Taylor 1c}
 \limsup_{n \to \infty} \sup_{0 < |\theta| < a} \frac 1{r_{n}^{2} \theta^{2}} |\max(\epsilon^{(n)}_{i}(\vartriangle,r_{n} \theta), \epsilon^{(n)}_{i}(1/r_{n},r_{n} \theta))| = 0.
\end{align}
\end{lemma}

This lemma is essentially the same as Lemma 4.3 of \cite{BravDaiMiya2015}. We now present the heavy traffic approximation. The following theorem is proved in \sectn{heavy}.

\begin{theorem}
\label{thr:heavy 1}
  For a sequence of the stable heterogeneous multiserver queues satisfying assumptions (i), (ii) and \eqn{stability 1}, assume \eqn{moment 1} and the heavy traffic conditions \eqn{heavy 1}--\eqn{heavy 3}, and let $L^{(n)}$ be a random variable subject to the stationary distribution of $L^{(n)}(\cdot)$, then we have
\begin{align}
\label{eqn:L moment 1}
  \lim_{n \to \infty} \dd{P}\Big( (1-\rho^{(n)}) L^{(n)} \le x\Big) = 1 - \exp\Big(-\frac {2\lambda_{0}} {\sum_{i \in \ol{K}} \lambda_{i}^{3} \sigma_{i}^{2}} \, x\Big), \qquad x \ge 0.
\end{align}
\end{theorem}

\begin{remark}
\label{rem:heavy 1}
\thr{heavy 1} is a direct extension of Theorem 2.1 of \cite{Miya2015} for the $GI/G/1$ queue. For a multiserver queue, this theorem generalizes Theorem 2 of \cite{Koll1974} which assumes homogeneous servers.
\end{remark}

The second scheme assumes large variances. In addition to \eqn{moment 1} and \eqn{heavy 1}, we assume that there is constant $b_{i} \ge 0$ for $i \in \ol{K}$ and a sequence of positive numbers ${s}_{n}$ for $n \ge 1$ satisfying the following conditions.
\begin{align}
\label{eqn:variance 1}
 & \lim_{n \to \infty} {s}_{n} = 0, \qquad \lim_{n \to \infty} {s}_{n} \big(\sigma^{(n)}_{i}\big)^{2} = b_{i}^{2}, \quad i \in \ol{K}, \qquad \sum_{i \in \ol{K}} b_{i}^{2} > 0,
\end{align}
which is referred to as a large variance condition. Let us consider the condition that
\begin{align}
\label{eqn:variance 3}
  \lim_{a \to \infty} \sup_{n \ge 1} {s}_{n} \dd{E}\big((T_{i}^{(n)})^{2} 1(T_{i}^{(n)} > a)\big) = 0, \qquad i \in \ol{K}.
\end{align}
If this were true, then we can get a weak limit approximation for $q_{n} = s_{n}$, replacing \eqn{heavy 1} by
\begin{align}
\label{eqn:rho fixed}
  \rho^{(n)} > 0, \qquad \lim_{n \to \infty} \rho^{(n)} = \rho < 1.
\end{align}
However, \eqn{variance 3} can not be true because, for each $a > 0$, \eqn{variance 1} implies
\begin{align*}
  \lim_{n \to \infty} {s}_{n} \dd{E}\big((T_{i}^{(n)})^{2} 1(T_{i}^{(n)} > a)\big) = b^{2}_{i}.
\end{align*}
Hence, the scaling factor must be the small order of $s_{n}$. From the stationary equation \eqn{SE 1}, we can see that the heavy traffic condition \eqn{heavy 1} is still needed and the scaling factor must be proportional to $r_{n} s_{n} \equiv (1-\rho^{(n)}) s_{n}$. In this case, \lem{Taylor 1} is changed to

\begin{lemma}
\label{lem:Taylor 2}
Assume \eqn{moment 1}, \eqn{heavy 1} and the large variance condition \eqn{variance 1}, then, for each $i \in \ol{K}$ and each $a > 0$, 
\begin{align}
\label{eqn:Taylor 2c}
  \limsup_{n \to \infty} \sup_{0 < |\theta| < a} \frac 1{s_{n} r_{n}^{2} \theta^{2}}  |\max(\epsilon^{(n)}_{i}(\vartriangle, r_{n} s_{n} \theta), \epsilon^{(n)}_{i}(1/(r_{n} s_{n}), r_{n} s_{n} \theta))| = 0.
\end{align}
\end{lemma}

This lemma is proved in \app{Taylor 2}. Using this lemma, we prove the following theorem in \sectn{variance}.

\begin{theorem}
\label{thr:variance 1}
  For a sequence of stable heterogeneous multiserver queues satisfying (i), (ii) and \eqn{stability 1}, assume \eqn{moment 1} and the large variance condition \eqn{variance 1} in addition to \eqn{heavy 1}, then we have
\begin{align}
\label{eqn:variance 5}
  \lim_{n \to \infty} \dd{P}\big( (1-\rho^{(n)}) {s}_{n} L^{(n)} \le x\big) = 1 - \exp\Big(-\frac {2 \lambda_{0}} {\sum_{i \in \ol{K}} \lambda_{i}^{3} b_{i}^{2}} \, x\Big), \qquad x \ge 0,
\end{align}
where $L^{(n)}$ is a random variable subject to the stationary distribution of $L^{(n)}(\cdot)$.
\end{theorem}

\section{Proofs of theorems}
\label{sect:proofs}
\setnewcounter

We first prove \thr{tail asymptotic 1}. Throughout Sections \sect{martingale} and \sect{tail asymptotic}, it is assumed that $F_{0}(0) = 0$. We will use a change of measure using a martingale obtained from \lem{martingale 2}, which starts from $L(0) = k$. However, we can not directly introduce such a martingale because we may require the truncations of $R_{i}(t)$ for $i \in K$. For this, we recall the definitions \eqn{limit point 2} for $\theta_{i}$ for $i \in K$ and $\alpha$ of \eqn{alpha 1}, and let, for $\theta \in \dd{R}$,
\begin{align}
\label{eqn:K alpha}
  K_{\theta} = \{i \in K; \theta_{i} \le \theta \}.
\end{align}
Obviously, $K_{\theta} = \emptyset$ for $\theta < 0$, and $K_{\alpha} = \emptyset$ for (a) of \thr{tail asymptotic 1}.

\subsection{Martingale and change of measure}
\label{sect:martingale}

We will truncate the service time distributions whose indexes are in $A \supset K_{\theta}$ for $\theta \in \dd{R}$. Thus, similarly to $M(t)$ of \eqn{martingale 3}, we define $M_{u,A(v),\theta}(t)$ for this $A$ and $u,v > 0$ as
\begin{align}
\label{eqn:martingale v}
  M_{u,A(v),\theta}(t) = f_{u,A(v),\theta}(X(t)) - \int_{0}^{t} \sr{A} f_{u,A(v),\theta}(X(s)) ds, \qquad t \ge 0,
\end{align}
which is well defined for $\theta > \ul{\theta}_{K}$ because $A \supset K_{\theta}$. From $\sr{A} f_{v,K(v),\theta}$ of \eqn{Af 1}, we have
\begin{align}
\label{eqn:Af 2}
 & \sr{A} f_{u,A(v),\theta}(X(s)) \nonumber\\
  & = \Big( \gamma_{u,A(v)}(\theta) + \eta(u,\theta) 1(R_{0}(s) \ge u) + \sum_{i \in A} \zeta_{i}(v,\theta) 1(R_{i}(s) \ge v) \nonumber\\
  & \qquad + \sum_{i \in K} \big( \zeta_{i}(\theta) 1(i \not\in A) + \zeta_{i}(v,\theta) 1(i \in A) \big) 1(X(s) \in \Gamma_{i}) \Big) f_{u,A(v),\theta}(X(s)),
\end{align}
where $\gamma_{u,A(v)}(\theta) = - \big(\eta(u,\theta) + \sum_{i \in K} \big(\zeta_{i}(\theta) 1(i \not\in A) + \zeta_{i}(v,\theta) 1(i \in A) \big)\big)$.

By \lem{terminal 1}, $f_{u,A(v),\theta}$ satisfies the terminal condition \eqn{terminal 1}, and therefore we can apply \lem{martingale 2} if $\dd{E}_{\nu}(|M_{u,A(v),\theta}(t)|) < \infty$ for each $t \ge 0$ and each $\theta \ge 0$ because $\ul{\theta}_{K} < 0$, where $\dd{E}_{\nu}$ stands for the expectation given that $X(0)$ has a distribution $\nu$. If nothing is mentioned, this $\nu$ is assumed to have a compact support concerning $L(0)$. To verify the finiteness of $\dd{E}_{\nu}(|M_{u,A(v),\theta}(t)|)$, we will use the following fact.

\begin{lemma}
\label{lem:super 1}
For each fixed $t > 0$ and for any $a > 0$, we have
\begin{align}
\label{eqn:super 1}
  \limsup_{n \to \infty} \frac 1{n} \log \dd{P}_{\nu}(N_{0}(t) \ge n) \le - a.
\end{align}
\end{lemma}

This lemma is intuitively clear, but we prove in \app{super 1} for making mathematical arguments transparent. The following lemma is a key for our arguments. For this, let $\widehat{\pi}$ be the stationary distribution of the embedded process of $X(\cdot)$ just before arrival instants, that is, $X(t_{0,n}-)$, which is extended for $n \le 0$ with $t_{0,0} = 0$. Note that $\widehat{\pi}$ is well defined because $X(\cdot)$ has the stationary distribution $\pi$. We define its conditional distribution $\nu_{k}$ on $(L(0-),L(0),\vc{R}(0-))$ as $\nu_{k}(\{\ell\} \times \{\ell'\} \times \vc{B}) = 1(\ell +1 = \ell' = k) \dd{P}_{\widehat{\pi}}(\vc{R}(0-) \in \vc{B}|L(0-) + 1 = L(0) = k)$ for $\ell, \ell' \in \dd{Z}_{+}$ and $\vc{B} \in \sr{B}(\dd{R}_{+}^{k+1})$.

\begin{lemma}
\label{lem:martingale basic}
Let $\delta_{0} = \infty$ if $K_{\alpha} = K$ while $\delta_{0} = \min \{\theta_{i}; i \in K \setminus K_{\alpha}\} - \alpha > 0$ otherwise. For $\theta \in [0, \alpha + \delta_{0})$, we have the following facts.\\
(a) For $t, v > 0$ and $v=\vartriangle$, $\sup_{s \in [0,t]} \big(\dd{E}_{\nu_{k}}\big(f_{K_{\alpha}(v),\theta}(X(s))\big) + \dd{E}_{\nu_{k}}(|M_{K_{\alpha}(v),\theta}(s)|)\Big) < \infty$.\\
(b) $\{M_{K_{\alpha}(v),\theta}(t); t \ge 0\}$ is an $\sr{F}_{t}$-martingale under $\dd{P}_{\nu_{k}}$.\\
(c) $M_{K_{\alpha}(v),\theta}(t)$ almost surely converges to $M_{K_{\alpha}(\vartriangle),\theta}(t)$ as $v \uparrow \infty$, and
\begin{align}
\label{eqn:martingale infty}
  M_{K_{\alpha}(\vartriangle),\theta}(t) = f_{K_{\alpha}(\vartriangle),\theta}(X(t)) - \int_{0}^{t} \sr{A} f_{K_{\alpha}(\vartriangle),\theta}(X(s)) ds, \qquad t \ge 0,
\end{align}
(d) $\{M_{K_{\alpha}(\vartriangle),\theta}(t); t \ge 0\}$ is an $\sr{F}_{t}$-supermartingale under $\dd{P}_{\nu_{k}}$.
\end{lemma}

\begin{remark}
\label{rem:martingale basic}
This lemma also holds true for $f_{u,K(v),\theta}$ and $M_{u,K(v),\theta}(t)$ for $u, v >0$ and $\theta > \ul{\theta}_{K}$ because all $R_{i}(t)$'s are truncated for $i \in \ol{K}$.
\end{remark}

\lem{martingale basic} enables us to apply an exponential change of measure following Palmowski and Rolski \cite{PalmRols2002}. Taking \eqn{martingale v} into account, we let
\begin{align*}
  E^{f_{K_{\alpha}(v),\theta}}(t) = \frac {f_{K_{\alpha}(v),\theta}(X(t))} {f_{K_{\alpha}(v),\theta}(X(0))} \exp\left( - \int_{0}^{t} \frac {\sr{A} f_{K_{\alpha}(v),\theta}(X(s))} {f_{K_{\alpha}(v),\theta}(X(s))} ds\right),
\end{align*}
then $E^{f_{K_{\alpha}(v),\theta}}(t)$ is an $\sr{F}_{t}$-martingale by Lemma 3.1 of \cite{PalmRols2002}. It follows from \eqn{Af 2} that
\begin{align}
\label{eqn:Ef 1}
    E^{f_{K_{\alpha}(v),\theta}}(t) = & \frac {f_{K_{\alpha}(v),\theta}(X(t))} {f_{K_{\alpha}(v),\theta}(X(0))} \exp\Big( - \gamma_{K_{\alpha}(v)}(\theta) t - \sum_{i \in K_{\alpha}} \zeta_{i}(v,\theta) \int_{0}^{t} \big(1(R_{i}(s) \ge v) ds \nonumber\\
    & - \sum_{i \in K} \big(\zeta_{i}(\theta) 1(i \not\in K_{\alpha}) + \zeta_{i}(v,\theta) 1(i \in K_{\alpha}) \big) \int_{0}^{t} 1(X(s) \in \Gamma_{i}) ds \Big).
\end{align}

Recall the conditional distribution $\nu_{k}$. Since $E^{f_{K_{\alpha}(v),\theta}}(t) > 0$ and $E^{f_{K_{\alpha}(v),\theta}}(0) = 1$, we can define a probability measure $\dd{P}^{(K_{\alpha}(v),\theta)}_{\nu_{k}}$ by
\begin{align}
\label{eqn:change 1a}
  \frac {d\dd{P}^{(K_{\alpha}(v),\theta)}_{\nu_{k}}} {d\dd{P}_{\nu_{k}}}\Big|_{\sr{F}_{t}} = E^{f_{K_{\alpha}(v),\theta}}(t),\qquad t \ge 0.
\end{align}
We denote expectations under $\dd{P}^{(K_{\alpha}(v),\theta)}_{\nu_{k}}$ by $\dd{E}^{(K_{\alpha}(v),\theta)}_{\nu_{k}}$. By Theorem 5.3 of \cite{PalmRols2002}, $X(\cdot)$ under $\dd{P}^{(K_{\alpha}(v),\theta)}_{\nu_{k}}$ is again a piecewise deterministic Markov process with generator $\sr{A}$ starting with the initial distribution $\nu_{k}$. To find the jump kernel $Q$ under the new measure $\dd{P}^{(K_{\alpha}(v),\theta)}_{\nu_{k}}$, we note the following fact, which is immediate from the definition \eqn{change 1a} and the truncation arguments.

\begin{lemma}[(a) of Proposition 3.8 in Chapter III of \cite{JacoShir2003}]
\label{lem:change 2}
Let $Z$ be a nonnegative random variable with finite expectation with respect to $\dd{P}^{(K_{\alpha}(v),\theta)}_{\nu_{k}}$, then
\begin{align}
\label{eqn:change 2}
  \dd{E}^{(K_{\alpha}(v),\theta)}_{\nu_{k}}(Z|\sr{F}_{s}) = \frac {1} {E^{f_{K_{\alpha}(v),\theta}}(s)} \dd{E}_{\nu_{k}}(E^{f_{K_{\alpha}(v),\theta}}(t) Z|\sr{F}_{s}), \qquad 0 \le s < t.
\end{align}
\end{lemma}

Letting $s \uparrow t$ in \eqn{change 2}, we have
\begin{align}
\label{eqn:change 3}
  \dd{E}^{(K_{\alpha}(v),\theta)}_{\nu_{k}}(Z|\sr{F}_{t-}) = \frac {1} {f_{K_{\alpha}(v),\theta}(X(t-))} \dd{E}_{\nu_{k}}( f_{K_{\alpha}(v),\theta}(X(t)) Z|\sr{F}_{t-}).
\end{align}
Thus, we define the transition kernel ${Q}^{(K_{\alpha}(v),\theta)}$ as, for $g \in C^{(1,p)}_{b, k+1}(S)$, 
\begin{align*}
  {Q}^{(K_{\alpha}(v),\theta)}g(X(t-)) = \frac {1} {f_{K_{\alpha}(v),\theta}(X(t-))} \dd{E}_{\nu_{k}}\big( f_{K_{\alpha}(v),\theta}(X(t)) g(X(t)) \big|X(t-)\big), \quad X(t-) \in \Gamma,
\end{align*}
then this is the jump transition kernel $Q$ under $\dd{P}^{(K_{\alpha}(v),\theta)}_{\nu_{k}}$. One can check that this is identical with (5.14) of \cite{PalmRols2002} because $f_{K_{\alpha}(v),\theta}$ satisfies the terminal condition \eqn{terminal 1}. Similarly, we define $\sr{F}_{t}$-supermartingale $E^{f_{K_{\alpha}(\vartriangle),\theta}}(t)$, possibly defective probability measure ${\dd{P}}_{\nu_{k}}^{(K_{\alpha}(\vartriangle),\theta)}$ and jump kernel ${Q}^{(K_{\alpha}(\vartriangle),\theta)}$ (see \supp{S4}). Obviously, those jump kernels do not depend on $\nu_{k}$.

We next consider the distributions of $T_{i}$ at the jump instant under $\dd{P}^{(K_{\alpha}(v),\theta)}_{\nu_{k}}$ for $v > 0$ and $v = \vartriangle$. They are assumed to be independent of the initial distribution of $(L(0-),\vc{R}(0-))$ under $\nu_{k}$, and may be written as $T_{i}^{(K_{\alpha}(v),\theta)}$. However, this notation is quite cumbersome, so we keep them as they are, which causes no confusion as long as measures for their expectations are specified.

\begin{lemma}
\label{lem:new T 1}
For each $v > 0$ and $\theta \in [0, \alpha + \delta_{0})$,
\begin{align}
\label{eqn:new T 1}
 & \dd{E}^{(K_{\alpha}(v),\theta)}_{\nu_{k}}(e^{{s} T_{0}}) = e^{\theta} \widehat{F}_{0}(\eta(\theta)+{s}), \qquad {s} \le \beta_{0} - \eta(\theta),\\
\label{eqn:new T 2}
 & \dd{E}^{(K_{\alpha}(v),\theta)}_{\nu_{k}}(e^{{s} T_{i}}) = e^{-\theta} \widehat{F}_{i}(\zeta_{i}(\theta)+{s}), \qquad {s} \le \beta_{i} - \zeta_{i}(\theta), i \in K \setminus K_{\alpha},\\
\label{eqn:new T 3}
 & \dd{E}^{(K_{\alpha}(v),\theta)}_{\nu_{k}}(e^{{s} (T_{i} \wedge v)}) = e^{-\theta} \widehat{F}_{i}(v, \zeta_{i}(v,\theta)+{s}), \qquad s \le \beta_{i} - \zeta_{i}(v,\theta), i \in K_{\alpha}.
\end{align}
The same formulas also hold true for $v=\vartriangle$.
\end{lemma}

\begin{remark}
\label{rem:new T 1}
(a) The right-hand terms of the inequalities for $s$ in \eqn{new T 1} and \eqn{new T 2} are the convergence exponents of the corresponding distributions because of \eqn{limit point 2}. \\
(b) Since $\eta(\alpha) \le 0$ and $\zeta_{i}(\alpha) < \beta_{i}$ for $i \in K \setminus K_{\alpha}$, ${\dd{E}}_{\nu_{k}}^{(K_{\alpha}(\vartriangle),\alpha)}(T_{i})$ is finite for $i=0$ and $i \in K \setminus K_{\alpha}$. Otherwise, consider $i \in K_{\alpha}$. Since $\beta_{i} - \zeta_{i}(\vartriangle,\alpha) = 0$ in this case, we have
\begin{align*}
  {\dd{E}}_{\nu_{k}}^{(K_{\alpha}(\vartriangle),\alpha)}(T_{i}) = e^{-\alpha} \dd{E}_{\nu_{k}}(T_{i} e^{\beta_{i} T_{i} }),
\end{align*}
which may be infinite as noted in \rem{concave 2}. Hence, even if $X(t)$ is a PDMP under ${\dd{P}}_{\nu_{k}}^{(K_{\alpha}(\vartriangle),\theta)}$, the corresponding service time of the $i$-th server may have an infinite mean. 
\end{remark}

\lem{new T 1} will be proved in \app{new T 1}. Since the generator $\sr{A}$ is unchanged under $\dd{P}^{(K_{\alpha}(v),\theta)}_{\nu_{k}}$ and ${\dd{P}}_{\nu_{k}}^{(K_{\alpha}(\vartriangle),\theta)}$, the distributions in \lem{new T 1} uniquely determine jump kernels ${Q}^{(K_{\alpha}(v),\theta)}$ and ${Q}^{(K_{\alpha}(\vartriangle),\theta)}$. Note that the distributions of $T_{0}$ and $T_{i}$ for $i \in K \setminus K_{\alpha}$ do not depend on $v$ under those change of measures.

Let us consider the mean drift of $L(t)$ in the off-boundary states under $\dd{P}^{(K_{\alpha}(v),\theta)}_{\nu_{k}}$ for $v >0$ and $v=\vartriangle$. For this, we compute the corresponding characteristics to $\gamma_{K_{\alpha}(v)}(\theta)$. It is easy to see that $\eta$, $\zeta_{i}$ for $i \in K \setminus K_{\alpha}$ and $\zeta_{i}(v,\cdot)$ for $i \in K_{\alpha}$ under $\dd{P}^{(K_{\alpha}(v),\theta)}_{\nu_{k}}$ are given by
\begin{align}
\label{eqn:new eta 1}
 & \eta^{(K_{\alpha}(v),\theta)}(\widetilde{\theta}) = \eta(\theta + \widetilde{\theta}) - \eta(\theta), \qquad \theta, \widetilde{\theta}+\theta \ge 0,\\
\label{eqn:new zeta 1}
 & \zeta^{(K_{\alpha}(v),\theta)}_{i}(\widetilde{\theta}) = \zeta_{i}(\theta + \widetilde{\theta}) - \zeta_{i}(\theta), \qquad \theta, \widetilde{\theta}+\theta \in [0, \alpha+\delta_{0}), i \in K \setminus K_{\alpha},\\
\label{eqn:new xi 2}
 & \zeta^{(K_{\alpha}(v),\theta)}_{i}(v,\widetilde{\theta}) = \zeta_{i}(v,\theta + \widetilde{\theta}) - \zeta_{i}(v,\theta), \qquad \theta, \widetilde{\theta}+\theta > \ul{\theta}_{K},  i \in K_{\alpha}.
\end{align}

From \eqn{new eta 1}--\eqn{new xi 2}, we have, for $v > 0$ and $v=\vartriangle$,
\begin{align}
\label{eqn:gamma n}
  \gamma^{(K_{\alpha}(v),\theta)}_{K_{\alpha}(v)}(\widetilde{\theta}) & \equiv - \Big(\eta^{(K_{\alpha}(v),\theta)}(\widetilde{\theta}) + \sum_{i \in K \setminus K_{\alpha}} \zeta^{(K_{\alpha}(v),\theta)}_{i}(\widetilde{\theta}) + \sum_{i \in K_{\alpha}} \zeta^{(K_{\alpha}(v),\theta)}_{i}(v,\widetilde{\theta}) \Big) \nonumber\\
  & = \gamma_{K_{\alpha}(v)}(\theta + \widetilde{\theta}) - \gamma_{K_{\alpha}(v)}(\theta).
\end{align}
Since $\zeta_{i}(\alpha) = \zeta_{i}(\vartriangle,\alpha)$ for $i \in K \setminus K_{\alpha}$, we have
\begin{align}
\label{eqn:alpha 2}
  \gamma_{K_{\alpha}(\vartriangle)}(\alpha) = \gamma_{K(\vartriangle)}(\alpha) = 0
\end{align}
Hence, letting $v = \vartriangle$ and $\theta = \alpha$ in \eqn{gamma n},
\begin{align}
\label{eqn:gamma alpha 1}
  \gamma^{(K_{\alpha}(\vartriangle),\alpha)}_{K_{\alpha}(\vartriangle)}(\widetilde{\theta}) = \gamma_{K_{\alpha}(\vartriangle)}(\alpha + \widetilde{\theta}) > 0, \qquad \widetilde{\theta} > 0,
\end{align}
because $\gamma_{K_{\alpha}(\vartriangle)}(\alpha + \cdot)$ is a convex function and $\gamma_{K_{\alpha}(\vartriangle)}'(0-) < 0$. They also implies that
\begin{align}
\label{eqn:gamma alpha 2}
  (\gamma^{(K_{\alpha}(\vartriangle),\alpha)}_{K_{\alpha}(\vartriangle)})'(0-) = \gamma'_{K_{\alpha}(\vartriangle)}(\alpha-) > 0.
\end{align}
This means that $L(t)$ has a positive mean drift at off-boundary states under ${\dd{P}}_{\nu_{k}}^{(K_{\alpha}(\vartriangle),\alpha)}$. Since $\gamma_{K_{\alpha}(v)}(\theta) \uparrow \gamma_{K_{\alpha}(\vartriangle)}(\theta)$ as $v \uparrow \infty$ by \lem{concave 1}, we have the following key facts.
\begin{lemma}
\label{lem:unstable 1}
(a) For all $v > 0$ and $\theta \in [0, \alpha]$, $\gamma_{K_{\alpha}(v)}(\theta) \le 0$.\\
(b) There are a $\theta > 0$ and $v_{0} > 0$ such that
\begin{align}
\label{eqn:gamma v 1}
  (\gamma^{(K_{\alpha}(v),\theta)}_{K_{\alpha}(v)})'(0-) > 0, \qquad \forall v \ge v_{0}.
\end{align}
(c) Either for $(v, \theta) = (\vartriangle, \alpha)$ or for each $v, \theta > 0$ satisfying \eqn{gamma v 1}, $L(t) \to \infty$ almost surely as $t \to \infty$ under ${\dd{P}}_{\nu_{k}}^{(K_{\alpha}(v),\theta)}$.
\end{lemma}

We can get back $\dd{P}$ from ${\dd{P}}_{\nu_{k}}^{(K_{\alpha}(v),\theta)}$. For this, let $Z(t)$ be a process which is $\sr{F}_{t}$ adapted and has a finite expectation for each $t \ge 0$. From \eqn{Ef 1} and \eqn{change 1a}, we have
\begin{align}
\label{eqn:get back 1}
  \dd{E}_{\nu_{k}} & \big(Z(t)\big) = {\dd{E}}_{\nu_{k}}^{(K_{\alpha}(v),\theta)} \big( Z(t) \nonumber\\
  & \times e^{-\theta (L(t) - L(0)) + \gamma_{K_{\alpha}(v)}(\theta) t - \eta(\theta) (R_{0}(t) - R_{0}(0)) - W_{K_{\alpha}(v),\theta}(\vc{R}(t)) + W_{K_{\alpha}(v),\theta}(\vc{R}(0))} \nonumber \\
  & \times e^{\sum_{i \in K_{\alpha}} \zeta_{i}(v,\theta) \int_{0}^{t} (1(R_{i}(s) > v) + 1(X(s) \in \Gamma_{i})) ds + \sum_{i \in K \setminus K_{\alpha}} \zeta_{i}(\theta) \int_{0}^{t} 1(X(s) \in \Gamma_{i}) ds} \big),
\end{align}
for $v > 0$ and $v = \vartriangle$, where
\begin{align*}
  W_{K_{\alpha}(v),\theta}(\vc{y}) = \sum_{i \in K \setminus K_{\alpha}} \zeta_{i}(\theta) y_{i} + \sum_{i \in K_{\alpha}} \zeta_{i}(v,\theta) (y_{i} \wedge v), \quad \vc{y} = (y_{0}, y_{1}, \ldots,y_{k}) \in \dd{R}_{+}^{k+1}.
\end{align*}

\subsection{The proof of \thr{tail asymptotic 1}}
\label{sect:tail asymptotic}

The spirit of our proof is similar to that of Theorem 2.3 of \cite{SadoSzpa1995}, which is based on the change of measure technique. However, the time index as well as the formulations are quite different. So, we will detail our proof. In this subsection, we present basic facts, then we prove the cases (a) and (b) in Sections \sect{exact} and \sect{logarithmic}, respectively.

Define random times $\sigma^{+}_{\ell}$ and $\sigma^{-}_{k}$ as
\begin{align*}
  \sigma^{+}_{\ell} = \inf \{ t > 0; L(t-) <  L(t) = \ell \}, \qquad  \sigma^{-}_{k} = \inf \{ t > 0; L(t-) = k  > L(t) \},
\end{align*}
then they are stopping times. By this definition, $\sigma^{+}_{\ell}$ is $\sr{F}_{\sigma^{+}_{\ell}-}$ measurable. Recall that $\nu_{k}$ specifies the conditional distribution of $\vc{R}(t_{0,0}-) \equiv \vc{R}(0-)$ under the stationary distribution of the embedded process $\{X(t_{0,n}-); n \in \dd{Z}\}$ given that $L(0-) + 1 = L(0-) = k$. Then, the stationary tail probability of $L$ is given by the so called cycle formula (e.g., see Corollary 2.1 of \cite{MiyaZwar2012}):
\begin{align}
\label{eqn:cycle 1}
  \dd{P}( L \ge \ell|L \ge k) = \frac {1} {\dd{E}_{\nu_{k}} (\sigma^{-}_{k})} \dd{E}_{\nu_{k}}\Big( \int_{0}^{\sigma^{-}_{k}} 1(L(s) \ge \ell) ds \Big), \qquad \ell \ge k.
\end{align}

From \eqn{get back 1} with $t=\sigma^{+}_{\ell}-$ (see, e.g., Theorem 3.4 in Chapter III of \cite{JacoShir2003}) and \lem{unstable 1}, which implies that ${\dd{P}}^{(K_{\alpha}(\vartriangle),\alpha)}_{\nu_{k}}(\sigma^{+}_{\ell} < \infty) = 1$, we have, for $\ell \ge k$,
\begin{align}
\label{eqn:L(t) 1}
 & \dd{E}_{\nu_{k}} \Big( \int_{0}^{\sigma^{-}_{k}} 1(L(s) \ge \ell) ds \Big) = \dd{E}_{\nu_{k}} \Big( \int_{\sigma^{+}_{\ell} }^{\sigma^{-}_{k}} 1(L(s) \ge \ell) ds 1(\sigma^{+}_{\ell} < \sigma^{-}_{k}) \Big) \nonumber\\
  & = \dd{E}_{\nu_{k}} \Big( \dd{E}_{\nu_{k}}\Big( \int_{\sigma^{+}_{\ell} \wedge \sigma^{-}_{k}}^{\sigma^{-}_{k}} 1(L(s) \ge \ell) ds \Big| \sr{F}_{\sigma^{+}_{\ell}-} \Big) 1(\sigma^{+}_{\ell} < \infty) \Big) \nonumber \\
    & = {\dd{E}}^{(K_{\alpha}(\vartriangle),\alpha)}_{\nu_{k}} \Big( \dd{E}_{\nu_{k}}\Big( \int_{\sigma^{+}_{\ell} \wedge \sigma^{-}_{k}}^{\sigma^{-}_{k}} 1(L(s) \ge \ell) ds \Big| \sr{F}_{\sigma^{+}_{\ell}-} \Big) 1(\sigma^{+}_{\ell} < \infty) \nonumber\\
  & \hspace{8ex} \times e^{- \alpha (L(\sigma^{+}_{\ell}-) - L(0)) - \eta(\alpha) (R_{0}(\sigma^{+}_{\ell}-) - R_{0}(0)) - W_{K_{\alpha}(\vartriangle),\alpha}(\vc{R}(\sigma^{+}_{\ell}-)) + W_{K_{\alpha}(\vartriangle),\alpha}(\vc{R}(0))} \Big) \nonumber\\
  & = e^{-\alpha (\ell-k-1)} {\dd{E}}^{(K_{\alpha}(\vartriangle),\alpha)}_{\nu_{k}} \Big( \dd{E}_{\nu_{k}}\Big( \int_{\sigma^{+}_{\ell} \wedge \sigma^{-}_{k}}^{\sigma^{-}_{k}} 1(L(s) \ge \ell) ds \Big| \sr{F}_{\sigma^{+}_{\ell}-} \Big) \nonumber\\
  & \hspace{25ex} \times e^{- W_{K_{\alpha}(\vartriangle),\alpha}(\vc{R}(\sigma^{+}_{\ell}-)) + W_{K_{\alpha}(\vartriangle),\alpha}(\vc{R}(0-))} \Big) ,
\end{align}
where we have used the fact that $\gamma_{K_{\alpha}(\vartriangle)}(\alpha) = 0$ and $R_{0}(0-)=0$.

Since $L(t)$ in \eqn{L(t) 1} is never below $k$, we can replace it by $H(t)$ defined as
\begin{align*}
  H(t) = L(0) + \widetilde{N}_{0}(t) - \sum_{i \in K} \widetilde{N}_{i}(t),
\end{align*}
where we recall that $\widetilde{N}_{0}$ and $\widetilde{N}_{i}$ are the renewal counting processes with the interval times $T_{0}$ and $T_{i}$, respectively. We denote the residual time to the next point of $N_{i}$ at time $t$ by $\widetilde{R}_{i}(t)$, which is identical with $R_{i}(t)$ for $t < \sigma^{-}_{k}$ when $L(0) = k$. Since the increments $\Delta H(t)$ only depend on $H(t-)$ and $\widetilde{\vc{R}}(t-) \equiv \{\widetilde{R}_{i}(t-); i \in \ol{K}\}$,
\begin{align*}
  Y(t) \equiv (H(t),\widetilde{\vc{R}}(t))
\end{align*}
is a continuous time Markov additive process with additive component $H(t)$ and background process $\{\widetilde{\vc{R}}(t); t \ge 0\}$. We assume that $Y(\cdot) \equiv \{Y(0-), Y(t); t \ge 0\}$ starts with $Y(0-)$ subject to $\nu_{k}$. In particular, $H(0)=k$ under this distribution.

Taking \rem{martingale basic} into account, we introduce the following $\sr{F}_{t}$-supermartingale for $Y(\cdot)$ with the initial distribution $\nu_{k}$ and for $\theta \in [0,\alpha+\delta_{0})$. 
\begin{align}
\label{eqn:Ef 2}
  \widetilde{E}^{\widetilde{f}_{u,K_{\alpha}(\vartriangle),\theta}}(t) = & \frac {\widetilde{f}_{u,K_{\alpha}(\vartriangle),\theta}(Y(t))} {\widetilde{f}_{u,K_{\alpha}(\vartriangle),\theta}(Y(0))} \nonumber\\
  & \times \exp\Big( -\gamma_{u,K_{\alpha}(\vartriangle)}(\theta) t - \eta(u,\theta) \int_{0}^{t} 1(\widetilde{R}_{0}(s) \ge u) ds \Big),
\end{align}
where $\widetilde{f}_{u,K_{\alpha}(\vartriangle),\theta}(Y(t))$ is defined as
\begin{align*}
  \widetilde{f}_{u,K_{\alpha}(\vartriangle),\theta}(Y(t)) = e^{\theta H(t) + \eta(u,\theta) \widetilde{R}_{0}(t) \wedge u + W_{K_{\alpha}(\vartriangle),\theta}(\widetilde{\vc{R}}(t))}.
\end{align*}

Similar to $\dd{P}^{(K_{\alpha}(\vartriangle),\theta)}_{\nu_{k}}$, we define  the probability measure $\widetilde{\dd{P}}^{(u,K_{\alpha}(\vartriangle),\theta)}_{\nu_{k}}$ on $(\Omega,\sr{F})$ by
\begin{align}
\label{eqn:change 1b}
  \frac {d\widetilde{\dd{P}}^{(u,K_{\alpha}(\vartriangle),\theta)}_{\nu_{k}}} {d{\dd{P}}_{\nu_{k}}}\Big|_{\sr{F}_{t}} = \widetilde{E}^{\widetilde{f}_{u,K_{\alpha}(\vartriangle),\theta}}(t), \qquad t \ge 0.
\end{align}
We denote expectation under $\widetilde{\dd{P}}^{(u,K_{\alpha}(\vartriangle),\theta)}_{\nu_{k}}$ by ${\widetilde{\dd{E}}}^{(u, K_{\alpha}(\vartriangle),\theta)}_{\nu_{k}}$. We prove the following lemma in \app{LD 1}.

\begin{lemma}
\label{lem:LD 1}
For $u > 0$ and $\theta \in [0,\alpha+\delta_{0})$ satisfying that $\gamma_{u,K_{\alpha}(\vartriangle)}(\theta) < 0$ and $\gamma_{u,K_{\alpha}(\vartriangle)}'(\theta) > 0$, there is a positive constant $C$ such that, for $\ell \ge k+1$,
\begin{align}
\label{eqn:h 1}
  \dd{E}_{\nu_{k}} \Big( e^{- W_{K_{\alpha}(\vartriangle),\theta}(\widetilde{\vc{R}}(\sigma^{+}_{\ell}-))} \int_{\sigma^{+}_{\ell}}^{\infty} 1(H(s) \ge \ell) ds \Big| \sr{F}_{\sigma^{+}_{\ell}-} \Big) < C.
\end{align}
\end{lemma}

Define $h_{+}(\ell)$ and $h_{-}(\ell)$ for $t \ge 0$ as
\begin{align*}
 & h_{+}(\ell) = \dd{E}_{\nu_{k}}\Big( \int_{\sigma^{+}_{\ell}}^{\infty} 1(H(s) \ge \ell) ds \Big| \sr{F}_{\sigma^{+}_{\ell}-} \Big) e^{- W_{K_{\alpha}(\vartriangle),\theta}(\widetilde{\vc{R}}(\sigma^{+}_{\ell}-))},\\
 & h_{-}(\ell) = \dd{E}_{\nu_{k}}\Big( \int_{\sigma^{-}_{k}}^{\infty} 1(H(s) \ge \ell) ds  1(\sigma^{+}_{\ell} < \sigma^{-}_{k}) \Big| \sr{F}_{\sigma^{+}_{\ell}-} \Big) e^{- W_{K_{\alpha}(\vartriangle),\theta}(\widetilde{\vc{R}}(\sigma^{+}_{\ell}-))}.
\end{align*}
By \lem{LD 1}, $h_{+}(\ell)$ and $h_{-}(\ell)$ are uniformly and almost surely bounded by $C$ over all $\ell \ge 0$. Then, \eqn{L(t) 1} can be written as,
\begin{align}
\label{eqn:L(t) 2}
 & e^{\alpha (\ell-k-1)} \dd{E}_{\nu_{k}} \Big( \int_{0}^{\sigma^{-}_{k}} 1(L(s) \ge \ell) ds \Big) \nonumber\\
 & \quad = \widetilde{\dd{E}}_{\nu_{k}}^{(K_{\alpha}(\vartriangle),\alpha)} \big( [h_{+}(\ell) - h_{-}(\ell)] \nonumber \\
 & \hspace{10ex} \times e^{ - (W_{K_{\alpha}(\vartriangle),\alpha}(\widetilde{\vc{R}}(\sigma^{+}_{\ell}-) - W_{K_{\alpha}(\vartriangle),\theta}(\widetilde{\vc{R}}(\sigma^{+}_{\ell}-)))} e^{W_{K_{\alpha}(\vartriangle),\alpha}(\vc{R}(0-))} \big).
\end{align}

We need two more lemmas. First one is about the Markov additive process $Y(\cdot)$ under $\widetilde{\dd{P}}_{\nu_{k}}^{(K_{\alpha}(\vartriangle),\alpha)}$. Since $\widetilde{R}_{i}(t)$ for $i \in \ol{K}$ are the residual times for the next counting times of independent renewal processes, $\{\widetilde{\vc{R}}(t)\}$ is Harris positive recurrent by the spread-out condition (iii), and therefore the embedded Markov process $\{\widetilde{\vc{R}}(t_{0,n}-)\}$ also is Harris positive recurrent under $\widetilde{\dd{P}}_{\nu_{k}}^{(K_{\alpha}(\vartriangle),\alpha)}$. Since this process can be considered as a background process of $\{Y(t_{n,0}-)\}$ and $H(t_{0,n}-)$ has a positive drift under $\widetilde{\dd{P}}_{\nu_{k}}^{(K_{\alpha}(\vartriangle),\alpha)}$, Theorem 1 of \cite{Alsm1996} (see also Theorem 4.1 of \cite{Alsm1997}) yields the following fact.
\begin{lemma}
\label{lem:Harris 1}
Under the assumptions of \thr{tail asymptotic 1}, $\{\widetilde{\vc{R}}(\sigma^{+}_{\ell}-); \ell=1,2,\ldots\}$ is Harris positive recurrent under $\widetilde{\dd{P}}_{\nu_{k}}^{(K_{\alpha}(\vartriangle),\alpha)}$, so it has the stationary distribution, which is denoted by $\nu_{\infty}(Y)$.
\end{lemma} 

Second one is about the distribution $\nu_{k}$. Recall that the embedded process $\{X(t_{0,n})\}$ has the stationary distribution because $X(\cdot)$ does so. Then, by the spread-out condition (ii), $\{\vc{R}(t_{0,n}-); n=1,2,\ldots\}$ is Harris positive recurrent. Thus, the following fact is obtained in Lemma 4.3 of \cite{SadoSzpa1995}, which assumes Harris positive recurrence for $\{X(t_{0,n}-)\}$, but its proof only requires it for $\{\vc{R}(t_{0,n}-)\}$ (see \supp{S2} for the outline of a proof for the lemma below).

\begin{lemma}
\label{lem:finite 2}
Under the assumptions of \thr{tail asymptotic 1}, $\dd{E}_{\nu_{k}}(e^{W_{K_{\alpha}(v),\theta}(\vc{R}(0-))}) < \infty$ for $v > 0$ and $v=\vartriangle$ for $\theta \in [0,\alpha + \delta_{0})$.
\end{lemma}

\subsubsection{Proof of (a) for exact decay rate}
\label{sect:exact}

By the assumptions in (a) of \thr{tail asymptotic 1}, $K_{\alpha} = \emptyset$, so $\dd{P}^{(K_{\alpha}(\vartriangle),\alpha)}_{\nu_{k}}$ is simply denoted by $\dd{P}^{(\alpha)}$. This also implies that there exist $u > 0$ and $\theta \in (0,\alpha)$ satisfying the conditions of \lem{LD 1}. Hence, from \lem{LD 1} and the fact that $\zeta_{i}(\alpha) - \zeta_{i}(\theta) > 0$ for $\alpha > \theta$, we have
\begin{align}
\label{eqn:h+ 1}
  \lim_{\ell \to \infty} {\dd{E}}_{\nu_{k}}^{(\alpha)} \big(h_{+}(\ell) & e^{ - \sum_{i \in K} (\zeta_{i}(\alpha) - \zeta_{i}(\theta)) R_{i}(\sigma^{+}_{\ell}-)} \big) \nonumber\\
  & = {\dd{E}}_{\nu_{\infty}(Y)}^{(\alpha)} \big( h_{+}(0) e^{- \sum_{i \in K} (\zeta_{i}(\alpha) - \zeta_{i}(\theta)) R_{i}(0-)} \big).
\end{align}
Note that $0 \le h_{-}(\ell) \le h_{+}(\ell)$. Furthermore, $1(H(s) \ge \ell) 1(\sigma^{+}_{\ell} < \sigma^{-}_{k})$ almost surely converges to $0$ as $\ell \to \infty$ because $\dd{P}_{\nu_{k}}(\sigma^{-}_{k} < \infty | \sr{F}_{\sigma^{+}_{\ell}-}) = 1$ for each $\ell \ge k+1$. Hence, $h_{-}(\ell)$ also converges to $0$ as $\ell \to \infty$. Thus, by \lem{finite 2}, it follows from \eqn{L(t) 2}, \eqn{h+ 1} and the bounded convergence theorem that
\begin{align}
\label{eqn:h+ 2}
  & \lim_{\ell \to \infty} e^{\alpha (\ell-k-1)} \dd{E}_{\nu_{k}}\Big( \int_{0}^{\sigma^{-}_{k}} 1(L(s) \ge \ell) ds \Big) \nonumber\\
  & \quad = {\dd{E}}_{\nu_{\infty}(Y)}^{(\alpha)} \big( h_{+}(0) e^{- \sum_{i \in K} (\zeta_{i}(\alpha) - \zeta_{i}(\theta)) R_{i}(0-)} \big) \dd{E}_{\nu_{k}}(e^{\sum_{i \in K} \zeta_{i}(\alpha) R_{i}(0-)}) < \infty.
\end{align}
Thus, \eqn{exact asymptotic 1} is obtained from \eqn{cycle 1}, which proves (a).

\subsubsection{Proof of (b) for logarithmic decay rate}
\label{sect:logarithmic}

To prove (b), we first consider an upper bound of the limit in \eqn{decay rate 1}. For this, we apply \lem{LD 1} to \eqn{L(t) 2}, then we have
\begin{align}
\label{eqn:L(t) 4}
 e^{\alpha (\ell-k-1)} & \dd{E}_{\nu_{k}}\Big( \int_{0}^{\sigma^{-}_{k}} 1(L(s) \ge \ell) ds \Big) \le C \widetilde{\dd{E}}_{\nu_{k}}^{(K_{\alpha}(\vartriangle),\alpha)} \big( e^{W_{K_{\alpha}(\vartriangle),\alpha}(\vc{R}(0-))} \big).
\end{align}
Hence, \eqn{cycle 1} and \lem{finite 2} conclude that
\begin{align}
\label{eqn:}
  \limsup_{\ell \to \infty} \frac 1{\ell} \log \dd{P}(L > \ell) \le -\alpha. 
\end{align}

We next consider a lower bound. For this, we first choose $\theta$ satisfying the conditions that, for some $v > 0, \theta > \alpha$,
\begin{align*}
  \gamma_{K_{\alpha}(v)}(\theta) < 0, \qquad \gamma_{K_{\alpha}(v)}'(\theta) > 0,
\end{align*}
which is possible because $\gamma_{K_{\alpha}(v)}(\theta) \uparrow \gamma_{K_{\alpha}(\vartriangle)}(\theta)$ as $v \uparrow \vartriangle$ and $\gamma_{K_{\alpha}(\vartriangle)}(\alpha) = 0$ by \eqn{alpha 2}. Then, it is easy to see that we can find a sufficiently large $u > 0$ by which the conditions of \lem{LD 1} are satisfied. We now use \eqn{get back 1} for $Z(t) = \dd{E}_{\nu_{k}}(\int_{t}^{\sigma^{-}_{k}} 1(L(s) \ge \ell) ds |\sr{F}_{t-}) 1(\sigma^{-}_{k} > t)$. Since $L(\sigma^{+}_{\ell}-) = \ell - 1$, $L(0) = k$, $\widetilde{R}_{0}(\sigma^{+}_{\ell}-) = \widetilde{R}_{0}(0-) = 0$, $\zeta_{i}(v,\theta) > 0$ and $\widetilde{\dd{P}}_{\nu_{k}}^{(K_{\alpha}(v),\theta)}(\sigma^{+}_{\ell} < \infty) = 1$, we have, similar to \eqn{L(t) 1},
\begin{align}
\label{eqn:L(t) 5}
  e^{(\ell - k-1) \theta} & \dd{E}_{\nu_{k}} \Big( \int_{0}^{\sigma^{-}_{k}} 1(L(s) \ge \ell) ds \Big) \nonumber \\
 & \ge \widetilde{\dd{E}}_{\nu_{k}}^{(K_{\alpha}(v),\theta)} \Big( \dd{E}_{\nu_{k}}\Big( \int_{\sigma^{+}_{\ell}}^{\sigma^{-}_{k}} 1(H(s) \ge \ell) ds \Big| \sr{F}_{\sigma^{+}_{\ell}-} \Big) e^{-W_{K_{\alpha}(v),\theta}(\vc{R}(\sigma^{+}_{\ell}-))} \Big) \nonumber\\
 & = \widetilde{\dd{E}}_{\nu_{k}}^{(K_{\alpha}(v),\theta)} \big( (h_{+}(\ell) - h_{-}(\ell)) e^{W_{K_{\alpha}(\vartriangle),\theta}(\vc{R}(\sigma^{+}_{\ell}-)) -W_{K_{\alpha}(v),\theta}(\vc{R}(\sigma^{+}_{\ell}-))}\big).
\end{align}
The last term of this inequality converges to
\begin{align*}
  \widetilde{\dd{E}}_{\nu_{\infty}(Y)}^{(K_{\alpha}(v),\theta)} \big( h_{+}(0) e^{W_{K_{\alpha}(\vartriangle),\theta}(\vc{R}(0-)) -W_{K_{\alpha}(v),\theta}(\vc{R}(0-))}\big) > 0
\end{align*}
as $\ell \to \infty$ by the following facts. First, $e^{W_{K_{\alpha}(\vartriangle),\theta}(\vc{R}(0)) -W_{K_{\alpha}(v),\theta}(\vc{R}(0))}$ is bounded by $1$ because $W_{K_{\alpha}(\vartriangle),\theta}(\vc{R}(0-)) -W_{K_{\alpha}(v),\theta}(\vc{R}(0-)) \le 0$. Second, $h_{-}(\ell)$ is not greater than $h_{+}(\ell)$ and vanishes as $\ell \to \infty$. Thus, we can apply \lem{LD 1}, the $\dd{P}_{\nu_{\infty}(Y)}^{(K_{\alpha}(v),\theta)}$-version of \lem{Harris 1} and the bounded convergence theorem. Hence, from \eqn{L(t) 5}, we have
\begin{align*}
  \liminf_{\ell \to \infty} \frac 1{\ell} \log \dd{P}(L \ge \ell) \ge -\theta.
\end{align*}
Since $\theta$ can be chosen to go to $\alpha$ as $v \uparrow \infty$, we arrive at the lower bound, which proves (b).

\subsection{Proof of \thr{heavy 1}}
\label{sect:heavy}

We consider the stationary equation of the $n$-th system corresponding to \eqn{SE 1} for $\theta \le 0$. Similar to $f_{v,K(v),\theta}$ of \eqn{test function 3}, we define $f^{(n)}_{v,K(v),\theta}$ for each $n \ge 1$ as
\begin{align*}
  f^{(n)}_{v,K(v),\theta}(\ell,\vc{y}) = e^{\theta \max(\ell,k)} e^{\eta^{(n)}(v,\theta) (y_{0} \wedge v) + \sum_{i \in K}\zeta^{(n)}_{i}(v,\theta) (y_{i} \wedge v)}, \qquad (\ell,\vc{y}) \in S.
\end{align*}
Define $\Phi^{(n)}$ as
\begin{align*}
  \Phi^{(n)}(v,\theta) = \dd{E}(f^{(n)}_{v,K(v),\theta}(X^{(n)})), \qquad v > 0, \theta \le 0,
\end{align*}
where $X^{(n)}$ is a random variable subject to the stationary distribution of the PDMP $X^{(n)}(\cdot) \equiv \{X^{(n)}(t); t \ge 0\}$, which describes the $n$-th system. Obviously, $\Phi^{(n)}(v,\theta)$ is finite for $\theta \le 0$. Similarly, we define, for $v > 0$ and $\theta \le 0$, 
\begin{align*}
 & \Phi^{(n)}_{i}(v,\theta) = \dd{E}(f^{(n)}_{v,K(v),\theta}(X^{(n)}) 1(R^{(n)}_{i} < v)), \qquad i \in \ol{K},\\
 & \Phi^{(n)}_{j,0}(v,\theta) = \dd{E}(f^{(n)}_{v,K(v),\theta}(X^{(n)}) 1(R^{(n)}_{j} = 0)), \qquad j \in K.
\end{align*}
Then, the stationary equation of the $n$-th system corresponding to \eqn{SE 1} is given by
\begin{align}
\label{eqn:SE n1}
  \eta^{(n)}(v,\theta) \Phi^{(n)}_{0}(v,\theta) + \sum_{i \in K} \zeta^{(n)}_{i}(v,\theta) \Phi^{(n)}_{i}(v,\theta) - \sum_{j \in K} \zeta^{(n)}_{i}(v,\theta) \Phi^{(n)}_{j,0}(v,\theta) = 0,
\end{align}
for $v > 0$ and $\theta \le 0$.

Let $r_{n} = 1 - \rho^{(n)}$. By the heavy traffic condition \eqn{heavy 1}, $r_{n}$ vanishes as $n \to \infty$. We scale $L^{(n)}(t)$ as $r_{n} L^{(n)}(t)$. To realize this, we let 
\begin{align*}
  \widetilde{\Phi}^{(n)}(\theta) = \Phi^{(n)}(1/r_{n}, r_{n}\theta), \quad \widetilde{\Phi}^{(n)}_{i}(\theta) = \Phi^{(n)}_{i}(1/r_{n}, r_{n}\theta), \quad \widetilde{\Phi}^{(n)}_{j,0}(\theta) = \Phi^{(n)}_{j,0}(1/r_{n}, r_{n}\theta).
\end{align*}
We correspondingly define $\widetilde{\eta}^{(n)}(\theta)$ and $\widetilde{\zeta}^{(n)}_{i}(\theta)$ as
\begin{align}
\label{eqn:tilde 1}
  \widetilde{\eta}^{(n)}(\theta) = \eta^{(n)}(1/r_{n}, r_{n} \theta), \qquad \widetilde{\zeta}^{(n)}_{i}(\theta) = \zeta^{(n)}_{i}(1/r_{n}, r_{n}\theta), \quad i \in K.
\end{align}
Then, \eqn{SE n1} can be written as
\begin{align}
\label{eqn:SE n2}
  \widetilde{\eta}^{(n)}(\theta) \widetilde{\Phi}^{(n)}_{0}(\theta) + \sum_{i \in K}  \widetilde{\zeta}^{(n)}_{i}(\theta) \widetilde{\Phi}^{(n)}_{i}(\theta) - \sum_{j \in K} \widetilde{\zeta}^{(n)}_{j}(\theta) \widetilde{\Phi}^{(n)}_{j,0}(\theta) = 0.
\end{align}

We claim that,
\begin{align}
\label{eqn:L heavy 1}
  \lim_{n \to \infty} \widetilde{\Phi}^{(n)}(\theta)  =  \frac {2 \lambda_{0}} {2 \lambda_{0} - \sum_{j \in \ol{K}} \lambda^{3}_{j} \sigma_{j}^{2} \theta}, \qquad \theta \le 0.
\end{align}
Suppose this claim is true, then the proof of \thr{heavy 1} is completed if we show that
\begin{align}
 \label{eqn:L heavy 2}
 \lim_{n \to \infty} |\widetilde{\Phi}^{(n)}(\theta) - \dd{E}(e^{r_{n} L^{(n)} \theta})| = 0.
\end{align}
Let $\widetilde{f}^{(n)}_{\theta}(X^{(n)}) = f^{(n)}_{1/r_{n},K(1/r_{n}),r_{n} \theta}(X^{(n)})$ for simplicity, then this is obtained by recalling the definition:
\begin{align*}
  \widetilde{\Phi}^{(n)}(\theta) = \dd{E} \big(\widetilde{f}^{(n)}_{\theta}(X^{(n)}) \big),
\end{align*}
and applying the bounded convergence theorem to
\begin{align*}
  |\widetilde{\Phi}^{(n)}(\theta) & - \dd{E}(e^{r_{n} L^{(n)} \theta}) |\\
  & \le \dd{E} \Big( \widetilde{f}^{(n)}_{\theta}(X^{(n)}) \Big| 1 - e^{-\theta r_{n} (k-|J^{(n)}|)} e^{-\widetilde{\eta}^{(n)}(\theta) (R^{(n)}_{0} \wedge 1/r_{n}) - \sum_{i \in K}\widetilde{\zeta}^{(n)}_{i}(\theta) (R^{(n)}_{i} \wedge v)} \Big| \Big)
\end{align*}
because $\widetilde{\eta}^{(n)}(\theta) (R^{(n)}_{0} \wedge 1/r_{n})$ and $\widetilde{\zeta}^{(n)}_{i}(\theta) (R^{(n)}_{i} \wedge 1/r_{n})$ vanish in probability as $n \to \infty$ by \eqn{bound 2} and $|r_{n} (R^{(n)}_{i} \wedge 1/r_{n})| \le 1$.

Thus, it remains to prove \eqn{L heavy 1}. For this, we will use \eqn{SE n1} and \lem{Taylor 1}. We first specialize this lemma in the following way.

\begin{lemma}
\label{lem:Taylor 3}
 Assume \eqn{moment 1} and the heavy traffic conditions \eqn{heavy 1}--\eqn{heavy 3}, and let, for each $n \ge 1$ and $r_{n} \theta \in (\ul{\theta}_{K}, \ol{\theta}_{0})$,
\begin{align}
\label{eqn:Taylor 3a}
 & \widetilde{\epsilon}^{(n)}_{0}(\theta) = \widetilde{\eta}^{(n)} \big(\theta\big) + \lambda^{(n)}_{0} r_{n} \theta + \frac 12 (\lambda^{(n)}_{0})^{3} (\sigma^{(n)}_{0})^{2} r_{n}^{2} \theta^{2}, \\
\label{eqn:Taylor 3b}
 & \widetilde{\epsilon}^{(n)}_{i}(\theta) = \widetilde{\zeta}^{(n)}_{i}\big(\theta\big) - \lambda^{(n)}_{i} r_{n} \theta + \frac 12 (\lambda^{(n)}_{i})^{3} (\sigma^{(n)}_{i})^{2} r_{n}^{2} \theta^{2},
\end{align}
then, for each $i \in \ol{K}$ and $c > 0$, there is a sequence $\{\widetilde{a}_{i, n}(c)\}$ vanishing as $n$ goes to infinity such that
\begin{align}
\label{eqn:Taylor 3c}
  \frac 1{r^{2}_{n}} |\widetilde{\epsilon}^{(n)}_{i}(\theta)| \le \widetilde{a}_{i,n}(c) |\theta|, \qquad n \ge 1, |\theta| < c.
\end{align}
\end{lemma}

 We defer the proof of this lemma to \app{Taylor 3}. We rewrite \eqn{SE n2} as
\begin{align}
\label{eqn:SE n3}
  \Big(\widetilde{\eta}^{(n)}(\theta) &+ \sum_{i \in K}  \widetilde{\zeta}^{(n)}_{i}(\theta)\Big) \widetilde{\Phi}^{(n)}(\theta) - \sum_{i \in K} \widetilde{\zeta}^{(n)}_{i}(\theta) \widetilde{\Phi}^{(n)}_{i,0}(\theta) \nonumber\\
  & = \widetilde{\eta}^{(n)}(\theta) \big(\widetilde{\Phi}^{(n)}(\theta) - \widetilde{\Phi}^{(n)}_{0}(\theta) \big) + \sum_{i \in K} \widetilde{\zeta}^{(n)}_{i}(\theta) \big(\widetilde{\Phi}^{(n)}(\theta) - \widetilde{\Phi}^{(n)}_{i}(\theta) \big).
\end{align}
We will divide both sides of this formula by $r_{n}^{2} \theta$, and let $n$ go to infinity. For this computation, we will use the following facts, which are obtained from \lem{Taylor 3}, \eqn{moment 1}, \eqn{heavy 1} and \eqn{heavy 2}.
\begin{align}
\label{eqn:heavy app 1}
 & \lim_{n \to \infty} \frac 1{r_{n}^{2} \theta} \Big(\widetilde{\eta}^{(n)}(\theta) + \sum_{i \in K}  \widetilde{\zeta}^{(n)}_{i}(\theta)\Big) = \lambda_{0}  - \frac 12 \sum_{i \in \ol{K}} \lambda_{i}^{3} \sigma_{i}^{2} \theta,\\
 \label{eqn:heavy app 2}
 & \lim_{n \to \infty} \frac 1{r_{n} \theta} \widetilde{\eta}^{(n)}(\theta) = - \lambda_{0},\\
 \label{eqn:heavy app 3}
 & \lim_{n \to \infty} \frac 1{r_{n} \theta} \widetilde{\zeta}^{(n)}_{i}(\theta) = \lambda_{i}.
\end{align}
As is well known, it follows from the renewal theorem that
\begin{align*}
  \dd{P}(R^{(n)}_{i} > v) = \lambda^{(n)}_{i} \dd{P}(R^{(n)}_{i} > 0) \int_{v}^{\infty} \dd{P}(T^{(n)}_{i} > x) dx, \qquad i \in \ol{K}.
\end{align*}
This and the uniformly integrable condition \eqn{heavy 3} yield that
\begin{align}
\label{eqn:R 1}
  \limsup_{n \to \infty} \frac 1{r_{n}} \dd{P}(R^{(n)}_{i} > 1/r_{n}) &\le \limsup_{n \to \infty} \frac 1{r_{n}} \lambda^{(n)}_{i} \int_{1/r_{n}}^{\infty} \dd{P}(T^{(n)}_{i} > x) dx \nonumber\\
  & \le \limsup_{n \to \infty} \lambda^{(n)}_{i} \dd{E}\Big(\int_{1/r_{n}}^{T^{(n)}_{i}} x 1(T^{(n)}_{i} > 1/r_{n}) dx \Big) \nonumber\\
  & \le \limsup_{n \to \infty} \frac 12 \lambda^{(n)}_{i} \dd{E}((T^{(n)}_{i})^{2} 1(T^{(n)}_{i} > 1/r_{n})) = 0,
\end{align}
Similarly, we have
\begin{align}
\label{eqn:R 2}
  \limsup_{n \to \infty} \dd{E}\big(R^{(n)}_{i} 1(R^{(n)}_{i} > 1/r_{n}) \big) \le \limsup_{n \to \infty} \frac 12 \lambda^{(n)}_{i} \dd{E}((T^{(n)}_{i})^{2} 1(T^{(n)}_{i} > 1/r_{n})) = 0.
\end{align}

Since $\widetilde{f}^{(n)}_{\theta}(X^{(n)})$ is uniformly bounded in $n$, \eqn{R 1} implies that
\begin{align*}
 0 \le \limsup_{n \to \infty} \frac 1{r_{n}} \big( \widetilde{\Phi}^{(n)}(\theta) - \widetilde{\Phi}^{(n)}_{i}(\theta) \big) \le \limsup_{n \to \infty} \dd{E}\Big( \widetilde{f}^{(n)}_{\theta}(X^{(n)}) \frac 1{r_{n}} 1(R^{(n)}_{i} > 1/r_{n}) \Big) = 0.
\end{align*}
This, \eqn{heavy app 2} and \eqn{heavy app 3} yield that the right-hand side of \eqn{SE n3} divided by $r_{n}^{2} \theta$ converges to zero as $n \to \infty$. Hence, if we can show that
\begin{align}
\label{eqn:boundary heavy 1}
  \lim_{n \to \infty} \frac 1{r_{n}^{2} \theta} \sum_{i \in K} \widetilde{\zeta}^{(n)}_{i}(\theta) \widetilde{\Phi}^{(n)}_{i,0}(\theta)  = \lambda_{0},
\end{align}
then dividing both sides of \eqn{SE n3} by $r_{n}^{2} \theta$ and letting $n \to \infty$ yield \eqn{L heavy 1} by \eqn{heavy app 1}. This is equivalent to \eqn{L moment 1}, and therefore the proof is completed. Thus, it remains to prove \eqn{boundary heavy 1}. By \eqn{heavy app 3}, \eqn{boundary heavy 1} is immediate if we have
\begin{align}
\label{eqn:boundary heavy 2}
  \lim_{n \to \infty} \frac 1{r_{n}} \sum_{i \in K} \lambda_{i} \widetilde{\Phi}^{(n)}_{i,0}(\theta)  = \lambda_{0}.
\end{align}
To prove this, we first note that dividing \eqn{SE n3} by $r_{n}^{2} \theta$ and letting $\theta \uparrow 0$ yield
\begin{align}
\label{eqn:boundary heavy 3}
  \lim_{n \to \infty} \frac 1{r_{n}} \sum_{i \in K} \lambda_{i} \widetilde{\Phi}^{(n)}_{i,0}(0) = \lambda_{0}.
\end{align}
Remark that this also follows from Little's law for customers in service at each server and the heavy traffic condition \eqn{heavy 1}. Thus, \eqn{boundary heavy 1} is obtained if we show that
\begin{align}
\label{eqn:boundary heavy 4}
  \lim_{n \to \infty} \frac 1{r_{n}} \big| \widetilde{\Phi}^{(n)}_{i,0}(\theta) - \widetilde{\Phi}^{(n)}_{i,0}(0) \big| = 0, \qquad \theta \le 0.
\end{align}
To prove this, let
\begin{align*}
  h_{n}(\theta) = \theta L^{(n)} + r_{n}^{-1} \widetilde{\eta}^{(n)}(\theta) (R^{(n)}_{0} \wedge 1/r_{n}) + \sum_{j \in K} r_{n}^{-1} \widetilde{\zeta}^{(n)}_{i}(\theta) (R^{(n)}_{j} \wedge 1/r_{n}),
\end{align*}
then $\widetilde{\Phi}^{(n)}_{i,0}(\theta) = \dd{E}(e^{r_{n} h_{n}(\theta)}1(R^{(n)}_{i} = 0))$. Thus, for $\theta < 0$, we have
\begin{align*}
  \frac 1{r_{n}} \big| \widetilde{\Phi}^{(n)}_{i,0}(\theta) - \widetilde{\Phi}^{(n)}_{i,0}(0) \big| & \le \dd{E} \Big( \Big| \frac {e^{r_{n} h_{n}(\theta)} - 1}{r_{n}} \Big|1(R^{(n)}_{i} = 0) \Big)\\
  & = \dd{E} \Big( \Big| \frac {e^{r_{n} h_{n}(\theta)} - 1}{r_{n}h_{n}(\theta)} h_{n}(\theta) \Big| 1(R^{(n)}_{i} = 0) \Big).
\end{align*}
Since $\lambda_{i} > 0$ for all $i \in \ol{K}$, it follows from \eqn{boundary heavy 3} that
\begin{align}
\label{eqn:boundary heavy 5}
  \lim_{n \to \infty} \frac {\dd{P}(R^{(n)}_{i} = 0)}{\sqrt{r_{n}}} = 0, \qquad i \in K.
\end{align}
On the other hand,
\begin{align*}
  \dd{E} & \big((R^{(n)}_{j} \wedge 1/r_{n})1(R^{(n)}_{i} = 0) \big) \\
  & = \dd{E} \big((R^{(n)}_{j} \wedge 1/r_{n})1(R^{(n)}_{i} = 0) (1(R^{(n)}_{j} > 1/\sqrt{r_{n}}) + 1(R^{(n)}_{j} \le 1/\sqrt{r_{n}}) \big) \\
  & \le \dd{E} \big(R^{(n)}_{j} (1(R^{(n)}_{j} > 1/\sqrt{r_{n}}) \big) + \frac 1{\sqrt{r_{n}}} \dd{P}( R^{(n)}_{i} = 0 ).
\end{align*}
Hence, \eqn{R 2} and \eqn{boundary heavy 5} imply that $\dd{E}(h_{n}(\theta) 1(R^{(n)}_{i} = 0))$ converges to 0 as $n \to \infty$. Thus, \eqn{boundary heavy 4} is obtained if $\frac {e^{r_{n} h_{n}(\theta)} - 1}{r_{n}h_{n}(\theta)}$ is uniformly bounded in $n$ for each $\theta < 0$, which follows from that $r_{n} h_{n}(\theta)$ is uniformly bounded in $n$. This uniform boundedness is implied by the facts that $L^{(n)} \le k-1$ on $\{R^{(n)}_{i} = 0\}$, $r_{n}^{-1} \widetilde{\eta}^{(n)}(\theta)$ and $r_{n}^{-1} \widetilde{\zeta}^{(n)}_{i}(\theta)$ are uniformly bounded in $n$ for $r_{n} \theta \in (\ul{\theta}_{K}, 0)$ by \lem{moment 1}. Consequently, \eqn{boundary heavy 4} is obtained, and we complete the proof. 

The above inequality is one of the key techniques in the present approach, and has been used in the proof of  \lem{Taylor 1} in \app{Taylor 1}. It also can be found in the last part of the proof of Lemma 4.7 of \cite{BravDaiMiya2015}.

\subsection{Proof for \thr{variance 1}}
\label{sect:variance}

Up to \eqn{SE n2}, all the computations in \sectn{heavy} are also valid for this case if we replace $r_{n}$ by $r_{n} s_{n}$ with $r_{n} = 1-\rho_{n}$. To distinguish the notation in \sectn{heavy}, we change $\widetilde{\eta}^{(n)}$, $\widetilde{\zeta}^{(n)}_{i}$ and $\widetilde{\epsilon}^{(n)}_{i}$ to $\ul{\widetilde{\eta}}^{(n)}$, $\ul{\widetilde{\zeta}}^{(n)}_{i}$ and $\ul{\widetilde{\epsilon}}^{(n)}_{i}$, respectively. Similar to \lem{Taylor 3}, we have

\begin{lemma}
\label{lem:Taylor 4}
 Assume \eqn{moment 1}, \eqn{heavy 1} and the large variance condition \eqn{variance 1}, and let, for each $n \ge 1$ and $\theta \in \dd{R}$,
\begin{align}
\label{eqn:Taylor 4a}
 & \ul{\widetilde{\epsilon}}^{(n)}_{0}(\theta) = \ul{\widetilde{\eta}}^{(n)} \big(\theta\big) + \lambda^{(n)}_{0} r_{n} s_{n} \theta + \frac 12 (\lambda^{(n)}_{0})^{3} r_{n}^{2} s_{n} (\sigma^{(n)}_{0})^{2} s_{n} \theta^{2}, \\
\label{eqn:Taylor 4b}
 & \ul{\widetilde{\epsilon}}^{(n)}_{i}(\theta) = \ul{\widetilde{\zeta}}^{(n)}_{i}\big(\theta\big) - \lambda^{(n)}_{i} r_{n} s_{n} \theta + \frac 12 (\lambda^{(n)}_{i})^{3} r_{n}^{2} s_{n} (\sigma^{(n)}_{i})^{2} s_{n} \theta^{2},
\end{align}
then, for each $i \in \ol{K}$ and $c > 0$, there is a sequence $\{\ul{\widetilde{a}}_{i, n}(c)\}$ vanishing as $n$ goes to infinity such that
\begin{align}
\label{eqn:Taylor 4c}
  \frac 1{r_{n}^{2} s_{n}} |\widetilde{\epsilon}^{(n)}_{i}(\theta)| \le \ul{\widetilde{a}}_{i, n}(c) |\theta|, \qquad n \ge 1, |\theta| < c.
\end{align}
\end{lemma}

This lemma yields similar formulas to \eqn{heavy app 1}--\eqn{heavy app 3}, in which $r_{n}$'s are replaced by $r_{n} s_{n}$. Thus, we have \eqn{variance 5} in the same way as \thr{heavy 1}.

\section{Concluding remarks}
\label{sect:concluding}
\setnewcounter

We here consider potential
of the present approach for more general models and some other asymptotics.

\subsection{More general arrival processes}
\label{sect:on arrival}

We have assumed that the arrival process is renewal or its independent superposition. However, the renewal process is not essential. As long as it is independent of the state of system, a more general class of arrival processes can be used in our framework. In this case, the terminal condition \eqn{terminal 1} needs to be changed. Suppose that $N_{0}$ be a stationary point process with finite intensity $\lambda_{0}$. Then, \eqn{eta 1} is changed to
\begin{align*}
  e^{\theta} \dd{E}_{0}(e^{\eta(\theta) T_{0}}|\sr{F}_{0-}) = 1,
\end{align*}
where $\dd{E}_{0}$ stands for the expectation by Palm distribution $\dd{P}_{0}$ concerning $N_{0}$, and $T_{0}$ is the time to the next counting instant measured from time $0$. Namely,
\begin{align*}
  \dd{E}_{0}(e^{\eta(\theta) T_{0}}) = \frac 1{\lambda_{0}} \dd{E}\Big( \sum_{i=1}^{\infty} e^{\eta(\theta) (t_{0,i+1} - t_{0,i})} 1(0 < t_{0,i} < 	1) \Big),
\end{align*}
where $t_{0,i} = \inf \{ u > t_{0,i-1}; N_{0}(\{u\}) > 0 \}$ for $i \ge 1$ with $t_{0,0} = 0$. Thus, $\eta(\theta)$ is a random variable in general, and not easy to its distribution in general. However, it may be possible to get it when $N_{0}$ is generated by a background Markov process. A typical example for this counting process is a Markov arrival process with finitely many background states. See  $\gamma^{(i,a)}(\theta)$ of (4.15) in \cite{Miya2015a}.

We next consider batch arrivals, which includes the case that $F_{0}(0) > 0$ as a special case (see \rem{tail asymptotic 2}). Denote the $n$-th arrival batch size by $A_{b,n}$. We assume that $\{A_{b,n}; n=1,2,\ldots\}$ are positive, independent and identically distributed. Furthermore they are independent of everything else. Denote their common distribution by $F_{b}$, and let $A_{b}$ be a random variable subject to $F_{b}$. We also need to verify the terminal condition \eqn{terminal 1}. This can be done similarly to the proof of \lem{terminal 1} in \app{terminal 1}, where $\ell'$ in the proof of this lemma should be changed to
\begin{align*}
  \ell' = \ell + 1(y_{0}=0) A_{b} - \sum_{i \in U \setminus U(\vc{y})} 1(y_{i} = 0),
\end{align*}
and the definition of $\eta(\theta,v)$ should be appropriately updated. For this, we need to truncate not only $T_{0}$ but also $A_{b}$. So, we define $\eta(m,v,\theta)$ for positive integer $m$ and positive real number $v$ as the solution of the following equation.
\begin{align}
\label{eqn:batch terminal 2}
  \dd{E}\big( e^{\theta (A_{b} \wedge m)} \big) \dd{E}\big( e^{\eta(m,v,\theta) (T_{0} \wedge v)} \big) = 1.
\end{align}
Similar to $\eta(\vartriangle,\theta)$, we define
\begin{align*}
  \eta(\vartriangle,\vartriangle,\theta) = \lim_{v \to \infty} \lim_{m \to \infty} \eta(m,v,\theta).
\end{align*}
It is not hard to see that these limiting operations are exchangeable. In particular, $\eta(\vartriangle,\vartriangle,\theta) = \eta(\vartriangle,\infty,\theta)$ for $\theta \ge 0$, and $\eta(\vartriangle,\vartriangle,\theta) = \eta(\infty,\vartriangle,\theta)$ for $\theta < 0$, where $\eta(m,v,\theta)$ is defined for $m = \infty$ or/and $v=\infty$ by \eqn{batch terminal 2}.

Using those notation, all the results in Sections \sect{hetero} and \sect{answers} can be extended to the batch arrival case under appropriate conditions. For example, $\alpha$ of \eqn{alpha 1} is redefined as
\begin{align}
\label{eqn:batch alpha 1}
  \alpha \equiv \sup\Big\{\theta \ge 0; - \Big(\eta(\vartriangle,\infty,\theta) + \sum_{i \in K} {\zeta}_{i}(\vartriangle,\theta) \Big) \le 0 \Big\}.
\end{align}

For the weak limit approximations, let $\eta^{(n)}(\infty,v,\theta)$ be the $\eta(\infty,v,\theta)$ of the $n$-th system, then define the error function of its 2nd order Taylor expansion for $\theta \le 0$ as
\begin{align}
\label{eqn:batch Taylor 1a}
 & \epsilon^{(n)}_{0}(\infty,v,\theta) = \eta^{(n)}(\infty,v,\theta) - \frac {\partial \eta^{(n)}}{\partial s} (\infty,v,s)\Big|_{s=0} \theta - \frac 12 \frac {\partial^{2} \eta^{(n)}}{\partial s^{2}} (\infty,v,s)\Big|_{s=0} \theta^{2},
\end{align}
since no truncation on $A_{b}$ is needed for $\theta \le 0$, where, for $\lambda^{(n)}_{0}(v) = \dd{E}(T^{(n)}_{0} \wedge v)$,
\begin{align}
\label{eqn:batch eta 1}
 & \frac {\partial \eta^{(n)}}{\partial s} (\infty,v,s)\Big|_{s=0} = - \lambda^{(n)}_{0}(v) \dd{E}(A_{b}), \\
\label{eqn:batch eta 2}
 & \frac {\partial^{2} \eta^{(n)}}{\partial s^{2}} (\infty,v,s)\Big|_{s=0} = - \lambda^{(n)}_{0}(v) \big\{ (\sigma^{(n)}_{A_{b}})^{2} + (\lambda^{(n)}_{0}(v) \dd{E}(A_{b}) \sigma^{(n)}_{0}(v))^{2}\big\},
\end{align}
where $(\sigma^{(n)}_{A_{b}})^{2}$ is the variance of $A_{b}^{(n)}$. Thus, we can apply exactly the same arguments to get the weak limits for the heavy traffic and large variance approximations.

\subsection{Some other problems and limitations}
\label{sect:some}

This paper only considers a single queue. However, the present approach is potentially applicable to systems with multiple queues and queueing networks.

For the tail asymptotic problem, this approach may relax phase-type assumptions in queues and their networks. For example, those assumptions are used for a two node generalized Jackson network in \cite{Miya2015a} and for the joint shortest queue in \cite{Saku2011}. They may be replaced by generally distributed assumptions for inter-arrival and service times.

For the weak limit approximations, a similar approach is studied by Braverman et al./ \cite{BravDaiMiya2015}. They use a slightly different test function, and consider the heavy traffic approximation of the stationary distribution for a generalized Jackson network, whose process limit is originally obtained by Reiman \cite{Reim1984}. As discussed in \cite{BravDaiMiya2015}, there are several issues to be overcome in the network case.

The tools of the present approach are based on the terminal condition \eqn{terminal 1}. However, this condition may cause serious drawbacks. For example, it will be hard to apply for a multi-class queueing network with first-come-first-served service discipline because we need to have information on the order of customer classes in a waiting line. Another hard problem is the heavy traffic approximation for the Halfin-Whitt regime (e.g., see \cite{HalfWhit1981,Reed2009}). We can describe the Halfin-Whitt regime by the present framework, but the stationary equations \eqn{SE 1} and \eqn{SE 2} may not be so useful because they do not well capture the number of busy servers.

Yet another problem of the present approach is about a process limit under the heavy traffic condition. We have only considered the stationary distribution of the PDMP. However, our starting formula \eqn{martingale 2} describes sample path evolutions. This suggests that the present approach may be applicable for a process limit.

Thus, there remain many challenging problems for future study.

\subsection*{Acknowledgements}
The author is grateful to an anonymous referee for helpful comments and suggestions. This paper is related to the recent network project \cite{BravDaiMiya2015} with Anton Braverman and Jim Dai. The author greatly appreciate their invaluable comments. This research is supported in part by JSPS KAKENHI Grant Number JP26540008.

\noindent {\it (Additional note)}: 
The author also is grateful to Toshiyuki Katsuda of Kwansei Gakuin University, who gave many helpful comments at the final stage of this paper.

\appendix

\section*{Appendix}

\section{Lemmas in \sectn{hetero}}
\label{app:section 2}
\setnewcounter

\subsection{Proof of \lem{martingale 1}}
\label{app:martingale 1}

Recall the definition \eqn{martingale 1} of $M_{0}(t)$. For a stopping time $\tau$, let $\sr{F}_{\tau-} = \sigma(\sr{F}_{0} \cup \{A \cap \{ t < \tau \}; A \in \sr{F}_{t} \} )$, then $\tau$ is $\sr{F}_{\tau-}$-measurable (see  1.1.14 of \cite{JacoShir2003}). Note that $t_{i}$ is $\sr{F}_{t_{i}-}$-measurable since $t_{i}$ is a stopping time, and $X(t_{i}-)$ is $\sr{F}_{t_{i}-}$-measurable. Hence, for $0 < s < t$,
\begin{align*}
  \dd{E}( M_{0}(t) | \sr{F}_{s}) & = M_{0}(s) + \sum_{i=1}^{\infty}  \dd{E} \Big( \dd{E}\big( f(X(t_{i})) - Qf(X(t_{i}-)) \big| \sr{F}_{t_{i}-} \big) 1(s < t_{i} \le t) \Big| \sr{F}_{s} \Big)\\
  & = M_{0}(s) + \sum_{i=1}^{\infty}  \dd{E} \Big( \big\{\dd{E}\big( f(X(t_{i})|\sr{F}_{t_{i}-} \big) - Qf(X(t_{i}-))\big\} 1(s < t_{i} \le t) \Big| \sr{F}_{s} \Big)\\
  & = M_{0}(s) + \sum_{i=1}^{\infty}  \dd{E} \Big( \big\{\dd{E}\big( f(X(t_{i})|X(t_{i}-)\big) - Qf(X(t_{i}-))\big\} 1(s < t_{i} \le t) \Big| \sr{F}_{s} \Big)\\
  & = M_{0}(s),
\end{align*}
where the third equality follows from the strong Markov property of $X(\cdot)$, and the last equality is obtained from the definition \eqn{Q 1} of $Q$. This concludes that $M_{0}(\cdot) \equiv \{M_{0}(t); t \ge 0\}$ is an $\sr{F}_{t}$-martingale since $\dd{E}(|M_{0}(t)|) < \infty$ by \eqn{finite 1}.

\subsection{Proof of \lem{terminal 1}}
\label{app:terminal 1}

We only prove \eqn{terminal 1} for $f = f_{u,K(v),\theta}$ with finite $v > 0$ since the proof is similar for $u,v=\infty$ if $\eta(\theta)$ and $\zeta_{i}(\theta)$ are finite. For this, we recall that
\begin{align*}
  \Gamma = \big\{(\ell, U, \vc{y}) \in S; y_{0} = 0 \mbox{ or } U \ne K \big\}.
\end{align*}
and $X(t-) \in \Gamma$ if and only if $\Delta N^{*}(t) = 1$. For $\vc{x} \equiv (\ell, U, \vc{y}) \in \Gamma$, let
\begin{align*}
  V(\vc{x}) = \{ i \in \{0\} \cup U; y_{i} = 0 \}.
\end{align*}
This $V(\vc{x})$ specifies an arrival or/and departures which occur at time $t$ when $X(s)$ goes to $\vc{x} \in \Gamma$ as $s \uparrow t$.

Suppose that $X(t-) = \vc{x}$ at time $t$. Denote this $X(t-)$ by $X(t,0)$. Let $X(t,1)$ be the state which is changed from $X(t,0)$ caused by the arrivals and/or departures of $V(\vc{x})$, where the server selection rule is applied at this stage independent of everything else. If $V(X(t,1)) = \emptyset$, then $X(t) = X(t,1)$, and the state transition at time $t$ is completed. Otherwise, we continue to define $X(t,2)$ using $V(X(t,1))$ and the server selection rule independent of everything else. In this way, we have the sequence of the states $X(t,n)$ for $n=0,1,2,\ldots$. Since the numbers of external arrivals and service completions at the same time instants are finite almost surely, there is a finite number $n$ for each sample path such that $V(X(t,n)) = \emptyset$ and $X(t,n) = X(t)$ almost surely. Denote the minimum of such $n$ by $\tau$, then $\dd{P}(\tau < \infty) = 1$. To describe these state transitions, we define $X(t,m) = X(t,\tau)$ for $m \ge \tau+1$. Then, we have, $X(t,\tau) = X(t)$, and therefore
\begin{align}
\label{eqn:decomposed limit}
  X(t) = \lim_{n \to \infty} X(t,n).
\end{align}

We next consider the transition kernel $Q_{0}$ defined as
\begin{align*}
  Q_{0}f(\vc{x}) = \dd{E}(f(X(t,1))| X(t,0) = \vc{x}), \qquad \vc{x} \in \Gamma, f \in C^{(1,p)}_{b,k+1}(S).
\end{align*}
Clearly, we have, for all $n =1,2,\ldots$,
\begin{align}
\label{eqn:Q 2}
  Q_{0}f(\vc{x}) = \dd{E}(f(X(t,n))| X(t,n-1) = \vc{x}), \qquad \vc{x} \in \Gamma, f \in C^{(1,p)}_{b,k+1}(S),
\end{align}
However, this $Q_{0}$ can not fully capture the transition from $X(t,n-1)$ to $X(t,n)$ because $X(t,n-1)$ may not be in $\Gamma$. To capture it, we extend $Q_{0}$ to the transition kernel $\widehat{Q}_{0}$ as
\begin{align*}
  \widehat{Q}_{0}f(\vc{x}) = \left\{
\begin{array}{ll}
 Q_{0}f(\vc{x}) \quad & \vc{x} \in \Gamma,\\
 f(\vc{x}) & \vc{x} \in S \setminus \Gamma.
\end{array}
\right.
\end{align*}
Then, the domain of $\widehat{Q}_{0}$ is $S$, and we have, for $n=1,2,\ldots$,
\begin{align*}
 & \widehat{Q}_{0}^{n-1}f(\vc{x}) = \dd{E}(f(X(t,n))| X(t,1) = \vc{x}), \qquad \vc{x} \in S, f \in C^{(1,p)}_{b,k+1}(S),\\
 & \widehat{Q}_{0}^{n}f(\vc{x}) = \dd{E}(f(X(t,n))| X(t,0) = \vc{x}), \qquad \vc{x} \in \Gamma, f \in C^{(1,p)}_{b,k+1}(S).
\end{align*}
Since these $f$'s are bounded functions, it follows from \eqn{decomposed limit} and the dominated convergence theorem that, for $\vc{x} \in \Gamma$,
\begin{align*}
  \dd{E}(f(X(t))| X(t-) = \vc{x}) & = \dd{E}(\lim_{n \to \infty} f(X(t,n))| X(t,0) = \vc{x})\\
  & = \lim_{n \to \infty} \dd{E}(f(X(t,n))| X(t,0) = \vc{x})\\
  & = \lim_{n \to \infty} \widehat{Q}_{0}^{n}f(\vc{x}).
\end{align*}
We now choose $f_{u,K(v),\theta}$ for $f$, and assume $\eta(u,\theta)$ and $\zeta_{i}(v,\theta)$ for $i \in K$ are well defined. If we show that
\begin{align}
\label{eqn:Q 3}
  Q_{0}f_{u,K(v),\theta}(\vc{x}) = f_{u,K(v),\theta}(\vc{x}), \qquad \vc{x} \in \Gamma,
\end{align}
then we have $\widehat{Q}_{0}f_{u,K(v),\theta}(\vc{x}) = f_{u,K(v),\theta}(\vc{x})$ for $\vc{x} \in S$, and therefore
\begin{align*}
  \dd{E}(f_{u,K(v),\theta}(X(t))| X(t-) = \vc{x})  = \lim_{n \to \infty} \widehat{Q}_{0}^{n}f_{u,K(v),\theta}(\vc{x}) = f_{u,K(v),\theta}(\vc{x}), \qquad \vc{x} \in \Gamma,
\end{align*}
which proves \eqn{terminal 1} for $f = f_{u,K(v),\theta}$.

Thus, it remains to prove \eqn{Q 3} to complete the proof. Recall that $U$ specifies busy servers. Suppose that, for given $(\ell, U, \vc{y}) \in \Gamma$, $X(t,0) \equiv (\ell, U, \vc{y})$ is changed to $X(t,1) \equiv (\ell', U', \vc{y}')$, where $\ell', U', \vc{y}'$ may be random. Let $U(\vc{y}) = \{ i \in U; y_{i} > 0\}$, which is the index set of servers who continue service at this transition. Since $U$ is the index set of busy servers and the residual inter-arrival and/or service times which are going to vanish at step $0$, we have
\begin{align}
\label{eqn:L to L'}
  \ell' = \ell + 1(y_{0}=0) - \sum_{i \in U \setminus U(\vc{y})} 1(y_{i} = 0) = \ell + 1(y_{0}=0) - |U| + |U(\vc{y})|.
\end{align}
Hence, the server index set which starts new service is $D(K \setminus U(\vc{y}), \ell' - |U(\vc{y}|)$, where we recall that $D(A,j)$ is the set of servers selected from $A$ for starting service of $\min (j, |A|)$ customers. Thus, $U'$ and $\vc{y}' \equiv \{y'_{i}; i \in \ol{K}\}$ are given by
\begin{align*}
 & U' = U(\vc{y}) \cup D(K \setminus U(\vc{y}), \ell' - |U(\vc{y}|),\\
 & y'_{0} = y_{0} + 1\big(y_{0}=0\big) T_{0}, \\
 & y'_{i} = y_{i} + 1\big(y_{i}=0, i \in U' \setminus U(\vc{y}) \big) T_{i}, \qquad i \in K,
\end{align*}
where $T_{i}$'s are independently chosen subject to $F_{i}$'s. Note that $y'_{i} = y_{i} > 0$ and $i \in U$ if and only if $i \in U(\vc{y})$. Since $(k - \ell') \vee 0$ servers are idle in the state $X(t,1)$ and $|U' \setminus U(\vc{y})| = (k - |U(\vc{y})|) \wedge (\ell' - |U(\vc{y})|)$, the definitions of $\eta(u,\theta)$ and $\zeta_{i}(v,\theta)$ yield
\begin{align*}
  {Q_{0}} & f_{u,K(v),\theta} (\vc{x}) = \dd{E}\big( f_{u,K(v),\theta}(X(t,1)) \big| X(t,0) = (\ell,U,\vc{y}) \big) \nonumber\\
  & = \dd{E}\big(e^{\theta (\ell' \vee k) + \eta(u,\theta) (T_{0} \wedge u) 1(y_{0} = 0) + \sum_{i \in U' \setminus U(\vc{y})} \zeta_{i}(v,\theta) (T_{i} \wedge v)}\big) e^{\eta(u,\theta) (y_{0} \wedge u) + \sum_{j \in U(\vc{y})} \zeta_{j}(v,\theta) (y_{j} \wedge v)} \nonumber\\
  & = \dd{E} \big( e^{\theta (\ell' \vee k)} \big) e^{-1(y_{0}=0) \theta + ((k - |U(\vc{y})|) \wedge (\ell' - |U(\vc{y})|)) \theta} e^{\eta(u,\theta) (y_{0} \wedge u) + \sum_{j \in U(\vc{y})} \zeta_{j}(v,\theta) (y_{j} \wedge v)} .
\end{align*}
Thus, we have \eqn{Q 3} if we can show that
\begin{align}
\label{eqn:exponent 1}
  \ell' \vee k -1(y_{0}=0) + ((k - |U(\vc{y})|) \wedge (\ell' - |U(\vc{y})|)) = \ell \vee k.
\end{align}
Since $|U' \setminus U(\vc{y})| = k \wedge \ell' - |U(\vc{y}|$ implies that
\begin{align*}
  \ell' \vee k + ((k - |U(\vc{y})|) \wedge (\ell' - |U(\vc{y})|)) = k+\ell' - |U(\vc{y})|
\end{align*}
and $|U| = \ell \wedge k$, it follows from \eqn{L to L'} that
\begin{align*}
  \mbox{The left-hand side of \eqn{exponent 1} } &= \ell' - 1(y_{0}=0) + k - |U(\vc{y})|\\
  &= \ell - |U| + k \\
  &= \ell - (\ell \wedge k) + k = \ell \vee k.
\end{align*}
Thus, we have proved \lem{terminal 1}.

\subsection{Proof of \lem{concave 1}}
\label{app:concave 1}

(a) Fix $v > 0$. Obviously, $\xi(v,\theta)$ is decreasing in $\theta$ by the definition of \eqn{xi 1}. We note that $\log \widehat{F}(v,s)$ is convex in $s \in \dd{R}$ (e.g., see Lemma 2.2.5 of \cite{DembZeit1998}). Hence, we have, for $p_{1}, p_{2} > 0$ such that $p_{1}+p_{2}=1$ and $\theta_{1}, \theta_{2} \, < \ol{\theta}$, 
\begin{align*}
  \log \widehat{F}(v,p_{1} \xi(v,\theta_{1}) + p_{2} \xi(v,\theta_{2})) & \le p_{1} \log \widehat{F}(v,\xi(v,\theta_{1})) + p_{2} \log \widehat{F}(v,\xi(v,\theta_{2}))\\
  & = - p_{1} \theta_{1} - p_{2} \theta_{2}\\
  & = \log \widehat{F}(v,\xi(v,p_{1} \theta_{1} + p_{2} \theta_{2})),
\end{align*}
where the last two equalities are obtained from \eqn{xi 1}. Since $\log \widehat{F}(v,s)$ is increasing in $s \in \dd{R}$, the above inequality implies that
\begin{align*}
  p_{1} \xi(v,\theta_{1}) + p_{2} \xi(v,\theta_{2}) \le \xi(v,p_{1} \theta_{1} + p_{2} \theta_{2}).
\end{align*}
This proves the concavity of $\xi(v,\theta)$ for $\theta \in \dd{R}$. It remains to prove for $\xi(v,\theta)$ to be infinitely differentiable. For this claim, we first consider the case that $\dd{P}(T > v) = 1$. In this case, $\xi(v,\theta) = - \theta/v$, and the claim is obviously true. Otherwise, $\widehat{F}(v,s)$ is analytic for all complex number $s$, and its derivative $\widehat{F}'(v,s) \ne 0$ for all real number $s$ because $\dd{E}(T) > 0$. Hence, $\xi(v,s)$ is analytic in a neighborhood around the real axis by the implicit function theorem. This proves the claim.\\
(b) These facts are easily inspected through \eqn{xi 1} because $\widehat{F}(v,\theta)$ is increasing in $v$ for $\theta > 0$ and decreasing in $v$ for $\theta<0$.\\
(c) By (a) and (b), we have, for $\theta \le 0$,
\begin{align*}
  0 \le \xi(v,\theta) \le \frac {\partial} {\partial x} \xi(v,x)\Big|_{x=0} \theta = - \frac 1{\dd{E}(T \wedge v)} \theta, \qquad v \ge 1,
\end{align*}
while, noting the fact that $\xi(v,\delta) < 0$ for $\delta > 0$, for $0 < \theta \le \delta$,
\begin{align*}
  \frac 1\delta \xi(v,\delta) \theta \le \xi(v,\theta) \le 0, \qquad v \ge 1.
\end{align*}
Thus, combining these two inequalities, we have \eqn{xi bound 1}.

\section{Lemmas in \sectn{answers}}
\label{app:section 3}
\setnewcounter

\subsection{Proof of \lem{moment 1}}
\label{app:moment 1}

The proof of \eqn{bound 3} is similar to that of \eqn{bound 2}, so we only show \eqn{bound 2}. We first note that
\begin{align*}
 & \eta^{(n)}(1/q_{n}, q_{n} \theta) \le \eta^{(n)}(\vartriangle, q_{n} \theta) < 0, \qquad 0 < \theta < \ol{\theta}_{0},\\
 & 0 \le \eta^{(n)}(\vartriangle, q_{n} \theta) \le \eta^{(n)}(1/q_{n}, q_{n} \theta), \qquad \theta \le 0,
\end{align*}
by \lem{concave 1}. Hence, we only need to prove that
\begin{align*}
   & \limsup_{n \to \infty} \Big|\frac {\eta^{(n)}(1/q_{n}, q_{n} \theta)} {q_{n}} \Big| \le c_{0}(\delta) |\theta|, \qquad |\theta| < \delta.
\end{align*}
However, this is already obtained in Lemma 4.4 of \cite{BravDaiMiya2015} assuming the uniformly integrable condition \eqn{heavy 3}, tacitly assuming that $\ol{\theta}^{(n)}_{0} > a$ for all $n \ge 1$ for some constant $a > 0$, which implies that $q_{n} \delta < a$ for sufficiently large $n$. This tacit assumption is easily checked because $\dd{E}(T^{(n)}_{0})$ converges to positive constant $1/\lambda_{0}$. Thus, we here show that the proof of Lemma 4.4 of \cite{BravDaiMiya2015} is valid without \eqn{heavy 3}.

It is easy to see that, in \cite{BravDaiMiya2015}, the condition \eqn{heavy 3} is only used to show that, for a given $y_{1}$, there exists a $y_{2} > y_{1}$ for any $\epsilon > 0$ such that
\begin{align*}
  \sup_{n \ge 1} \dd{P}(T^{(n)}_{0} \ge y_{2}) < \epsilon/2,
\end{align*}
which is (A.10) of \cite{BravDaiMiya2015}, where $T^{(n)}_{e,i}$ is replaced by $T^{(n)}_{0}$ of our notation. Suppose that this claim is not true, then, for some $\epsilon > 0$ and any $y_{2} > y_{1}$, there is a $n_{0} \ge 1$ such that $\dd{P}(T^{(n_{0})}_{i} \ge y_{2}) \ge \epsilon/2$, which implies that
\begin{align*}
  \dd{E}(T^{(n_{0})}_{i}) \ge y_{2} \dd{P}(T^{(n_{0})}_{i} \ge y_{2}) \ge y_{2} \epsilon/2.
\end{align*}
Hence, we can choose $y_{2}$ and $n_{0}$ such that $y_{2} \epsilon$ is arbitrarily large. This contradicts the condition \eqn{moment 1}. Thus, \eqn{moment 1} implies the claim, and therefore (A.10) of \cite{BravDaiMiya2015} is obtained without \eqn{heavy 3}.  Similarly, \eqn{moment 1} and the uniformly integrability of $\{T_{i}^{(n)}; n \ge 1\}$ imply (A.9). See \supp{S3} for its proof.

\subsection{Proof of \lem{Taylor 1}}
\label{app:Taylor 1}

We only prove \eqn{Taylor 1c} for $i=0$ since the proof is similar for $i \in K$. Fix $n \ge 1$ and $v > 0$. For simplicity, we denote $\eta^{(n)}(v,\theta)$ by $f$ in this subsection. By Taylor expansion of $f(\theta)$ around $\theta = 0$, we have, for some $\delta \in [0,1]$, 
\begin{align}
\label{eqn:Taylor 1d}
  \eta^{(n)}(v,\theta) = f'(0) \theta + \frac 12 f''(0) \theta^{2} + \frac 12( f''(\delta \theta)  - f''(0))\theta^{2},
\end{align}
since $f(0) = 0$. Taking derivatives of \eqn{xi 1} with $F = F^{(n,v)}$ and $\xi = f$, we have that
\begin{align}
\label{eqn:f 3}
  f''(\theta) & = - \frac {(\widehat{F}_{0}^{(n,v)})''(f(\theta)) (\widehat{F}_{0}^{(n,v)}(f(\theta)))^{2}} {((\widehat{F}_{0}^{(n,v)})'(f(\theta)))^{3}} + \frac {\widehat{F}_{0}^{(n,v)}(f(\theta))} {(\widehat{F}_{0}^{(n,v)})'(f(\theta))},\\
\label{eqn:f 4}
  f''(0) & = - (\lambda^{(n)}_{0}(v))^{3} (\sigma^{(n)}_{0}(v))^{2}.
\end{align}
Hence, from \eqn{Taylor 1a}, we have
\begin{align*}
  \epsilon^{(n)}_{0}(v,\theta) = \frac 12 ( f''(\delta \theta)  - f''(0))\theta^{2}, \qquad v >0, \exists \delta \in [0,1].
\end{align*}

For $\theta = 0$, obviously $\epsilon^{(n)}_{i}(v,\theta) = 0$, and \eqn{Taylor 1c} trivially holds. In what follows, we assume that $0 < \theta < \ol{\theta}_{0}$ since the case $\theta < 0$ is similarly proved. Note that
\begin{align}
\label{eqn:e function 1}
  |\epsilon^{(n)}_{0}(v,\theta)| \le \frac 12 \theta^{2} \max_{0 \le x \le \theta} |f''(x)  - f''(0)|.
\end{align}
By letting $v \uparrow \infty$, we can see that this inequality is valid for $v=\vartriangle$ because $\eta^{(n)}(\vartriangle,\theta) = \eta^{(n)}(\theta) < 0$ for $\theta > 0$ by \lem{concave 2}. In computations below, we assume that $0 < v < \infty$, but they are easily adapted for $v=\vartriangle$ by letting $v \uparrow \infty$, and lead the same conclusion.

To simply notation, let $g(s) = \widehat{F}^{(n,v)}_{0}(s)$. Then, $g'(0) = \dd{E}(T^{(n)}_{0} \wedge v)$, and \eqn{f 3} can be written as
\begin{align*}
  f''(x) = - \frac{g''(f(x))(g(f(x)))^{2}} {(g'(f(x)))^{3}} + \frac{g(f(x))} {g'(f(x))},
\end{align*}
and therefore, using the fact that $f(0) = 0$,
\begin{align*}
  f''(x) - f''(0) = \frac {g(f(x)) g'(0) - g'(f(x)) g(0)} {g'(0) g'(f(x))} + \frac {g''(0) (g'(f(x)))^{3} - g''(f(x)) (g'(0))^{3}(g(f(x)))^{2}} {(g'(0) g'(f(x)))^{3}}.
\end{align*}
We evaluate the numerators as
\begin{align*}
 & g(f(x)) g'(0) - g'(f(x)) g(0) = g(0) \big(g'(0) - g'(f(x))\big) + g'(0) \big(g(f(x)) - g(0)\big),\\
 & g''(0) (g'(f(x)))^{3} - g''(f(x)) (g'(0))^{3} (g(f(x)))^{2}\\
 & = g''(0) \big((g'(f(x)))^{3} - (g'(0))^{3}\big) + (g'(0))^{3} \big(g''(0) - g''(f(x))\big) + (g'(0))^{3} g''(f(x)) (1-(g(f(x))^{2}).
\end{align*}
Hence, it is important to evaluate the differences $g(f(x)) - g(0)$, $g'(f(x)) - g'(0)$ and $g''(f(x)) - g''(0)$. Since $f(x)$ is decreasing in $x$ and $g(s), g'(s), g''(s)$ are all increasing in $s$, we have, using Taylor expansions again,
\begin{align*}
 & \max_{0 \le x \le \theta} |g(f(x)) - g(f(0))| = |g(f(0)) - g(f(\theta))| \le |f(\theta)| g'(f(0)) = |f(\theta)| \dd{E}(T^{(n)}_{0}),\\
 & \max_{0 \le x \le \theta} |g'(f(x)) - g'(f(0))| = |g'(f(0)) - g'(f(\theta))| \le |f(\theta)| g''(f(0)) = |f(\theta)| \dd{E}((T^{(n)}_{0})^{2}),\\
 & \max_{0 \le x \le \theta} |g''(f(x)) - g''(f(0))| = g''(f(0)) - g''(f(\theta))\\
 & \hspace{10ex} = |\dd{E}\big( (T^{(n)}_{0} \wedge v)^{2} (1 - e^{(T^{(n)}_{0} \wedge v) f(\theta)}) \big)|\\
 & \hspace{10ex} \le \dd{E}\big( (T^{(n)}_{0} \wedge v)^{2} e^{(T^{(n)}_{0} \wedge v) |f(\theta)|} (1-e^{-(T^{(n)}_{0} \wedge v) |f(\theta)|})\big)\\
 & \hspace{10ex} \le |f(\theta)| e^{v |f(\theta)|} \dd{E}\big((T^{(n)}_{0})^{2} v (v^{-1} T^{(n)}_{0} \wedge 1))\big)\\
 & \hspace{10ex} \le v |f(\theta)| e^{v |f(\theta)|} \dd{E}\big((T^{(n)}_{0})^{2} (v^{-1} T^{(n)}_{0} \wedge 1) (1(T^{(n)}_{0} > \sqrt{v}) + 1(T^{(n)}_{0} \le \sqrt{v}))\big)\\
 & \hspace{10ex} \le v |f(\theta)| e^{v |f(\theta)|} \dd{E}\big((T^{(n)}_{0})^{2} \big(1(T^{(n)}_{0} > \sqrt{v}) + \sqrt{v}^{-1} \big)\big).
\end{align*}
Similarly, we have
\begin{align*}
  \max_{0 \le x \le \theta} \frac {1} {g'(0) g'(f(x))} \le \frac 1{(g'(0))^{2}} e^{|v f(\theta)|}.
\end{align*}
Since $g'(0) = \dd{E}(T^{(n)}_{0} \wedge v)$ and $g''(0) = \dd{E}((T^{(n)}_{0} \wedge v)^{2})$, we have
\begin{align*}
 & \dd{E}(T^{(n)}_{0}) - \dd{E}(T^{(n)}_{0} 1(T^{(n)}_{0} > v)) \le g'(0) \le \dd{E}(T^{(n)}_{0}),\\
 & \dd{E}((T^{(n)}_{0})^{2}) - \dd{E}((T^{(n)}_{0})^{2} 1(T^{(n)}_{0} > v)) \le g''(0) \le \dd{E}((T^{(n)}_{0})^{2}).
\end{align*}

We now let $v = 1/r_{n}$ and replace $\theta$ by $r_{n} \theta$, then $\frac 1{r_{n}} f(r_{n} \theta)$ is uniformly bounded by the linear function of $\lambda_{0} |\theta|$ as $n \to \infty$ by \lem{moment 1}. Hence, recalling $f(\theta) =\eta^{(n)}(v,\theta)$ and letting $f_{n}(x) = \eta^{(n)}(1/r_{n},r_{n} x)$,
\begin{align*}
  \limsup_{n \to \infty} \max_{0 \le x \le \theta} |f_{n}''(x)  - f_{n}''(0)| =0.
\end{align*}
Thus, in \eqn{e function 1}, let either $v=1/r_{n}$ or letting $v \uparrow \infty$ and substitute $r_{n} \theta$ to $\theta$, then applying the above limit completes the proof.

\subsection{Proof of \lem{Taylor 2}}
\label{app:Taylor 2}

The proof is essentially the same as \lem{Taylor 1}. Because of the conditions \eqn{heavy 1} and \eqn{heavy 2}, we need to divide the error functions by $s_{n} r_{n}^{2}$ instead of $r_{n}^{2}$. Except this manipulation, all the arguments in the proof of \lem{Taylor 1} go through.

\section{Lemmas in \sectn{proofs}}
\label{app:section 4}
\setnewcounter

\subsection{Proof of \lem{super 1}}
\label{app:super 1}

Recall that $F_{0}(0) = 0$ is assumed. It implies that $\dd{P}(N_{0}(t) \ge n) = \dd{P}(T_{0}(1)+T_{0}(2)+\ldots + T_{0}(n) \le t)$, so there is a sufficiently large $n \ge 1$ for any $\delta > 0$ and $t \ge 0$ such that $t \le n \delta$, which implies that
\begin{align*}
  \dd{P}(N_{0}(t) \ge n) \le \dd{P}(T_{0}(1)+T_{0}(2)+\ldots + T_{0}(n) \le n \delta).
\end{align*}
Hence, \eqn{super 1} is obtained if we can find $\delta > 0$ for any $a > 0$ such that
\begin{align}
\label{eqn:Cramer 1}
  \limsup_{n \to \infty} \frac 1n \log \dd{P}\Big(\frac 1n (T_{0}(1)+T_{0}(2)+\ldots + T_{0}(n)) \le \delta \Big) \le - a.
\end{align}
To prove this inequality, we use the Cram\'er theorem on large deviations (e.g., see Theorem 2.2.3 of \cite{DembZeit1998}) that, for any $\delta > 0$, 
\begin{align}
\label{eqn:Cramer 2}
  \limsup_{n \to \infty} \frac 1n \log \dd{P}\Big(\frac 1n (T_{0}(1)+T_{0}(2)+\ldots + T_{0}(n)) \le \delta \Big) \le - \inf_{x \le \delta} \Lambda^{*}(x),
\end{align}
where $\Lambda^{*}(x) = \sup_{\eta \in \dd{R}} (\eta x - \log \widehat{F}_{0}(\eta))$. Since $\dd{E}(T_{0}) > 0$ and $\widehat{F}_{0}(\eta)$ is finite for $\eta \le 0$ by \lem{concave 1}, it follows from Lemma 2.2.5 of \cite{DembZeit1998} that, for $\delta < \dd{E}(T_{0})$,
\begin{align}
\label{eqn:Cramer 3}
  \inf_{x \le \delta} \Lambda^{*}(x) = \Lambda^{*}(\delta) = \sup_{\eta < 0} (\eta \delta - \log \widehat{F}_{0}(\eta)).
\end{align}
Hence, the proof is completed if $\Lambda^{*}(\delta)$ can be arbitrarily large for some $\delta \in (0,\dd{E}(T_{0}))$. Let $\theta = - \log \widehat{F}_{0}(\eta)$, then $\eta = \xi(\theta)$ for $\xi$ of \eqn{xi 1} with $v=\infty$ and $T=T_{0}$. Then, by \lem{concave 1}, $\xi(\theta) < 0$ if and only if $\theta > 0$, and therefore $\Lambda^{*}(\delta)$ can be represented as
\begin{align*}
  \Lambda^{*}(\delta) = \sup_{\theta > 0} (\xi(\theta) \delta + \theta).
\end{align*}
For each $\theta > 0$, we can choose a sufficiently small $\delta \in (0,\dd{E}(T_{0}))$ such that $\Lambda^{*}(\delta) \ge \frac 12 \theta$. Hence, letting $a = 2 \theta$, \eqn{Cramer 1} follows from \eqn{Cramer 2} and \eqn{Cramer 3}. 

\subsection{Proof of \lem{martingale basic}}
\label{app:f infty finite}

(a) Since $\theta \ge 0$,  $\eta(\theta) R_{0}(t)$ can be omitted in the exponent in $f_{K_{\alpha}(v),\theta}(X(t))$ to prove the finiteness of its expectation because ${\eta}(\theta) \le 0$ for $\theta \ge 0$ by \lem{concave 1}. Since
\begin{align}
\label{eqn:L bound}
  L(t) \le L(0) + N_{0}(t), \qquad t \ge 0,
\end{align}
\lem{super 1} and $\dd{P}_{\nu_{k}}(L(0)=k) = 1$ lead to
\begin{align}
\label{eqn:E N0 finite}
  \dd{E}_{\nu_{k}}(e^{\theta L(t)}) \le e^{\theta k} \dd{E}_{\nu_{k}}(e^{\theta N_{0}(t)}) < \infty, \qquad \theta \ge 0.
\end{align}
We next check the finiteness of $\dd{E}_{\nu_{k}}(e^{W_{K_{\alpha}(v),\theta}(\vc{R}(t))})$ for $\theta \in [0,\alpha+\delta_{0}]$. Let us consider $R_{i}(0-)$ on $\nu_{k}$. Let $U_{t}(0) = \{i \in K; R_{i}(0-) \le t\}$. On $\nu_{k}$, one of $R_{i}(0-)$ vanishes for $i \in K$, so $U_{t}(0) \in 2^{K} \setminus \{\emptyset\}$. Hence, it is enough to verify that $\dd{E}_{\nu_{k}}(e^{W_{K_{\alpha}(v),\theta}(\vc{R}(t))}1(U_{t}(0) = A)) < \infty$ for each $A \in 2^{K} \setminus \{\emptyset\}$.

We now fix $A \in 2^{K} \setminus \{\emptyset\}$, and consider $\vc{R}(t)$ on $U_{t}(0) = A$ under $\nu_{k}$. If $i \in K \setminus A$, then $R_{i}(t) = R_{i}(0-) - t > 0$. Otherwise, if $i \in A$, then server $i$ may start service in the time interval $[R_{i}(0-),t]$. We apply similar arguments in the proof of Lemma 4.3 of \cite{SadoSzpa1995}. For this, let $t_{s}(m)$ be the $m$-th time when service is started and let $t_{s,i}(n)$ be the $n$-th time when server $i$ starts service, where if more than one server simultaneously start service, we count them separately. Then, there is a unique $n$ for each $m \ge 1$ such that $t_{s,i}(n) \le t_{s}(m) < t_{s,i}(n+1)$. Denote this $n$ by $J_{i}(m)$. For $\vc{y} \in \dd{R}_{+}^{k}$, let 
\begin{align*}
  g_{A_{\alpha}(v),\theta}(\vc{y}) = 1(U_{t}(0) = A, y_{i} \ge 0, i \in A) & \prod_{i \in A \setminus K_{\alpha}} \dd{E} \big( e^{\zeta_{i}(\theta) (T_{s,i} - y_{i})} 1(T_{s,i} > y_{i}) \big)\\
  & \times \prod_{i \in A \cap K_{\alpha}} \dd{E} \big( e^{\zeta_{i}(v,\theta) (v \wedge (T_{s,i} - y_{i}))} 1(T_{s,i} > y_{i}) \big),
\end{align*}
then, from the fact that $N_{0}(t)+1$ customers arrive in $[0,t]$, we have
\begin{align}
\label{eqn:bound 1}
  \dd{E}_{\nu_{k}}\big( e^{W_{K_{\alpha}(v),\theta}(\vc{R}(t))} 1(U_{t}(0) & = A) \big) \le \dd{E}_{\nu_{k}} \Big( \sum_{m=1}^{N_{0}(t)+1} e^{W_{A^{c}_{\alpha}(v),\theta}(\vc{R}(0-) - t \vcn{1})} \nonumber\\
 & \times g_{A_{\alpha}(v),\theta}(t - t_{s,1}(J_{1}(m)), \ldots,t - t_{s,k}(J_{k}(m)))\Big),
\end{align}
where $W_{A^{c}_{\alpha}(v),\theta}(\vc{y}) = \sum_{i \in K \setminus (A \cup K_{\alpha})} \zeta_{i}(\theta) y_{i} + \sum_{i \in K_{\alpha} \cap (K \setminus A)} \zeta_{i}(v,\theta) (y_{i} \wedge v)$. Hence, recalling the definitions of $\zeta_{i}(\theta)$ and $\zeta_{i}(v,\theta)$, \lem{finite 2} and \eqn{bound 1} yield
\begin{align}
\label{eqn:finite 2}
  \dd{E}_{\nu_{k}}\big( e^{W_{K_{\alpha}(v),\theta}(\vc{R}(t))} 1(U_{t}(0) & = A) \big) \le \dd{E}_{\nu_{k}}\big( e^{W_{A^{c}_{\alpha}(v),\theta}(\vc{R}(0-))}\big) e^{|A|\theta} \dd{E}_{\nu_{k}}(N_{0}(t)+1) < \infty.
\end{align}
Since $|2^{A} - \{\emptyset\}| = 2^{k} - 1 < 2^{k}$, \eqn{E N0 finite} implies that
\begin{align}
\label{eqn:finite 3}
  \dd{E}_{\nu_{k}}(f_{K_{\alpha}(v),\theta}&(X(t))) \nonumber\\
  & < e^{\theta k} \dd{E}_{\nu_{k}}(e^{\theta N_{0}(t)}) \dd{E}_{\nu_{k}}(e^{W_{K_{\alpha}(v),\theta}(\vc{R}(0-))}) 2^{k} e^{(k-1)\theta} \dd{E}_{\nu_{k}}(N_{0}(t)+1),
\end{align}
which is finite for each $t,v > 0$ and $v = \vartriangle$ by \lem{finite 2}. Since it follows from \eqn{martingale v} and \eqn{Af 2} for $u=\infty$ and $A = K_{\alpha}(v)$ that, for $\theta \ge 0$,
\begin{align}
\label{eqn:finite 4}
  |M_{K_{\alpha}(v),\theta}(t)| \le f_{K_{\alpha}(v),\theta}(X(t)) + (|\eta(\theta)|+\sum_{i \in K} \zeta_{i}(v,\theta)) \int_{0}^{t} f_{K_{\alpha}(v),\theta}(X(s)) ds,
\end{align}
$\dd{E}_{\nu_{k}}(|M_{K_{\alpha}(v),\theta}(t)|)$ also is finite for $t,v > 0$, $v=\vartriangle$ and $\theta \in [0,\alpha+\delta_{0})$. This proves (a), and therefore Lemmas \lemt{martingale 2} and \lemt{terminal 1} imply (b).\\
(c) Since $e^{{\zeta}_{i}(v,\theta) R_{i}(t) \wedge v}$ almost surely converges to $e^{{\zeta}_{i}(\vartriangle,\theta) R_{i}(t)}$ as $v \to \infty$ by the definition of ${\zeta}_{i}(\vartriangle,\theta)$, $f_{K_{\alpha}(v),\theta}(X(u))$ for $u \in [0,t]$ (see \supp{S4}), and therefore $M_{K_{\alpha}(v),\theta}(t)$ almost surely converges to $f_{K_{\alpha}(\vartriangle),\theta}(X(u))$ and $M_{K_{\alpha}(\vartriangle),\theta}(t)$, respectively, as $v \to \infty$.\\
(d) By (b), $M_{K_{\alpha}(v),\theta}(t)$ is an $\sr{F}_{t}$-martingale under $\dd{P}_{\nu_{k}}$ for each $v > 0$, which implies that
\begin{align}
\label{eqn:martingale K}
  \dd{E}_{\nu_{k}}( M_{K_{\alpha}(v),\theta}(t) 1_{D}) = \dd{E}_{\nu_{k}}( M_{K_{\alpha}(v),\theta}(s) 1_{D}), \qquad D \in \sr{F}_{s}, 0 \le s < t.
\end{align}
However, this may not imply (d) because of limiting operation although $M_{K_{\alpha}(v),\theta}(s)$ almost surely converges to $M_{K_{\alpha}(\vartriangle),\theta}(s)$ as $v \uparrow \infty$ for each $s \ge 0$. Instead of using \eqn{martingale K}. we can directly prove (d) by substituting $f = f_{K_{\alpha}(\vartriangle),\theta}$ into \eqn{martingale 2} and \lem{martingale 1}. See \supp{S4} for its details.

\subsection{Proof of \lem{new T 1}}
\label{app:new T 1}
We only prove \eqn{new T 1} because \eqn{new T 2} and \eqn{new T 3} are similarly proved. Since, at arrival instant $t$, for $L(t-) \ge k$,
\begin{align*}
  \frac {f_{K_{\alpha}(v),\theta}(X(t))} {f_{K_{\alpha}(v),\theta}(X(t-))} &= e^{\theta} e^{\eta(\theta) T_{0}},
\end{align*}
and, for $L(t-) \le k-1$,
\begin{align*}
  \frac {f_{K_{\alpha}(v),\theta}(X(t))} {f_{K_{\alpha}(v),\theta}(X(t-))} &= e^{\theta} e^{\eta(\theta) T_{0}} \prod_{i \in \Delta J(t) \setminus K_{\alpha}} e^{-\theta} e^{\zeta_{i}(\theta) T_{i}} \prod_{i \in \Delta J(t) \cap K_{\alpha}} e^{-\theta} e^{\zeta_{i}(v,\theta) (T_{i} \wedge v)}
\end{align*}
where $\Delta J(t) = J(t) - J(t-)$ and either $\Delta J(t) \setminus K_{\alpha}$ or $\Delta J(t) \cap K_{\alpha}$ is an empty set, it follows from \eqn{change 3} and the definitions of $\zeta_{i}$ and $\zeta^{(n)}_{i}$ that
\begin{align*}
  \dd{E}^{(K_{\alpha}(v),\theta)}_{\nu_{k}}(e^{{u} T_{0}}|R_{0}(t-)=0) & = \dd{E}_{\nu_{k}}\Big( \frac {f_{K_{\alpha}(v),\theta}(X(t))} {f_{K_{\alpha}(v),\theta}(X(t-))} e^{uT_{0} } \Big| R_{0}(t-)=0 \Big)\\
  & = e^{\theta} \dd{E}_{\nu_{k}}(e^{(\eta(\theta)+{u}) T_{0} }) = e^{\theta} \widehat{F}_{0}(\eta(\theta)+{u}), \qquad \eta(\theta) + {u} \le \beta_{0}.
\end{align*}

\subsection{Proof of \lem{LD 1}}
\label{app:LD 1}

This lemma is partly a continuous time version of Lemma 4.5 of \cite{SadoSzpa1995}, and we will use the same idea for our proof. Let $\widetilde{\tau}_{0}(t) = \inf \big\{ s \ge t; H(s) \ge \ell \big\}$, then $\int_{0}^{\infty} 1(H(s) \ge \ell) ds \ge t$ implies that $\widetilde{\tau}_{0}(t) < \infty$. Let $A_{t} = \{\widetilde{\tau}_{0}(t) < \infty\}$. To evaluate $\dd{P}_{\nu_{k}}\big(A_{t} | \sr{F}_{\sigma^{+}_{\ell}-} \big)$, we will use the probability measure ${\widetilde{\dd{P}}}^{(u,K_{\alpha}(\vartriangle),\theta)}_{\nu_{k}}$ of \eqn{change 1b}. Since $\gamma'_{u,K_{\alpha}(\vartriangle)}(\theta) > 0$, we have $\widetilde{\dd{P}}^{(u,K_{\alpha}(\vartriangle),\theta)}_{\nu_{k}}(A_{t}| \sr{F}_{\sigma^{+}_{\ell}-} ) = 1$. Furthermore, $\eta(u,\theta) \le 0$ for $\theta \ge 0$. Thus, we have, from \eqn{Ef 2},
\begin{align*}
 & \dd{E}_{\nu_{k}} (e^{- W_{K_{\alpha},\theta}(\widetilde{\vc{R}}(\sigma^{+}_{\ell}-))} 1_{A_{t}}| \sr{F}_{\sigma^{+}_{\ell}-} ) = \widetilde{\dd{E}}^{(u,K_{\alpha}(\vartriangle),\theta)}_{\nu_{k}} \Big( e^{- W_{K_{\alpha},\theta}(\widetilde{\vc{R}}(\sigma^{+}_{\ell}-))} \frac {\widetilde{f}_{u,K_{\alpha}(\vartriangle),\theta}(Y(\sigma^{+}_{\ell}-))} {\widetilde{f}_{u,K_{\alpha}(\vartriangle),\theta}(Y(\widetilde{\tau}_{0}(t)))} \nonumber\\
  & \hspace{30ex} \times e^{\gamma_{u,K_{\alpha}(\vartriangle)}(\theta) \widetilde{\tau}_{0}(t) + \eta(u,\theta) \int_{0}^{\widetilde{\tau}_{0}(t)} 1(\widetilde{R}_{0}(s) \ge u) ds} \Big| \sr{F}_{\sigma^{+}_{\ell}-} \Big)\\
  & \le \widetilde{\dd{E}}^{(u,K_{\alpha}(\vartriangle),\theta)}_{\nu_{k}} \big(e^{-\theta (H(\widetilde{\tau}_{0}(t)) - (\ell-1))- \eta(u,\theta) (\widetilde{R}_{0}(\widetilde{\tau}_{0}(t)) \wedge u) - W_{K_{\alpha},\theta}(\widetilde{\vc{R}}(\widetilde{\tau}_{0}(t))) + \gamma_{u,K_{\alpha}(\vartriangle)}(\theta) t} \big| \sr{F}_{\sigma^{+}_{\ell}-} \big)\\
  & \le \widetilde{\dd{E}}^{(u,K_{\alpha}(\vartriangle),\theta)}_{\nu_{k}} \big(e^{ - \eta(u,\theta) u + \gamma_{u,K_{\alpha}(\vartriangle)}(\theta) t} \big| \sr{F}_{\sigma^{+}_{\ell}-} \big),
\end{align*}
where the last inequality uses the fact that $H(\widetilde{\tau}_{0}(t)) \ge \ell$ and $\zeta_{i}(\vartriangle, \theta) \ge 0$ for $\theta \ge 0$. Integrating both sides of this inequality for $t \in [0,\infty)$, we have \eqn{h 1} for $C = (-\gamma_{u,K_{\alpha}(\vartriangle)}(\theta))^{-1} e^{ - \eta(u,\theta) u}$ because $\int_{0}^{\infty} 1(H(s) \ge \ell) ds \le \int_{0}^{\infty} 1_{A_{s}} ds$.

\subsection{Proof of \lem{Taylor 3}}
\label{app:Taylor 3}

We only need to prove \eqn{Taylor 3c} for $i=0$ because the proof is similar for $i \in K$. We introduce the following notation for the $n$-th system.
\begin{align*}
 & \widetilde{\lambda}^{(n)}_{0} = \big(\dd{E}(T^{(n)}_{0} 1(T^{(n)}_{0} \le 1/r_{n}))\big)^{-1}, \\
 & \widetilde{\sigma}^{(n)}_{0} = \dd{E}((T^{(n)}_{0})^{2} 1(T^{(n)}_{0} \le 1/r_{n})) - (\widetilde{\lambda}^{(n)}_{0})^{-2}.
\end{align*}
Comparing \eqn{Taylor 1a} with \eqn{Taylor 3a}, we can see that \lem{Taylor 3} follows from \lem{Taylor 1} if we have,
\begin{align}
\label{eqn:moment 1a}
 & \limsup_{n \to \infty} \frac 1{r_{n}} |\lambda^{(n)}_{0} - \widetilde{\lambda}^{(n)}_{0}| = 0,\\
\label{eqn:moment 2a}
 & \limsup_{n \to \infty} |(\sigma^{(n)}_{0} )^{2}- (\widetilde{\sigma}^{(n)}_{0})^{2}| = 0.
\end{align}
We will use the uniform integrability condition \eqn{heavy 3} to verify them. We first verify \eqn{moment 1a} in the following way.
\begin{align*}
   \limsup_{n \to \infty} \frac 1{r_{n}} |\lambda^{(n)}_{0} - \widetilde{\lambda}^{(n)}_{0}| &= \limsup_{n \to \infty} \frac 1{r_{n}}  \lambda^{(n)}_{0} \widetilde{\lambda}^{(n)}_{0} \dd{E}\big(T_{0}^{(n)} (1 - 1(T_{0}^{(n)} \le 1/r_{n})) \big) \\
   &= \limsup_{n \to \infty} \frac 1{r_{n}}  \lambda^{(n)}_{0} \widetilde{\lambda}^{(n)}_{0} \dd{E}\big(T_{0}^{(n)}  1(T_{0}^{(n)} > 1/r_{n}) \big) \\
   &\le \limsup_{n \to \infty} \lambda^{(n)}_{0} \widetilde{\lambda}^{(n)}_{0} \dd{E}\big((T_{0}^{(n)})^{2} 1(T_{0}^{(n)} > 1/r_{n}) \big) = 0,
\end{align*}
where the last equality is obtained by \eqn{moment 1} and \eqn{heavy 3}. Similarly, \eqn{moment 2a} is verified by
\begin{align*}
   & \limsup_{n \to \infty} |(\sigma^{(n)}_{0})^{2} - (\widetilde{\sigma}^{(n)}_{0})^{2}|\\
   & \quad \le \limsup_{n \to \infty} \big( \dd{E}\big((T_{0}^{(n)})^{2} 1(T_{0}^{(n)} > 1/r_{n}) \big) + (\lambda^{(n)}_{0} \widetilde{\lambda}^{(n)}_{0})^{-2} |(\lambda^{(n)}_{0} - \widetilde{\lambda}^{(n)}_{0}) (\lambda^{(n)}_{0} + \widetilde{\lambda}^{(n)}_{0})| \big) = 0.
\end{align*}
This completes the proof of \lem{Taylor 3}.

\newpage

\section*{Supplement for some proofs}
\setcounter{equation}{0}
\setcounter{lemma}{0}
\setcounter{remark}{0}

This paper will be published in Advances in Applied Probability, 2017. In this published version, we have omitted proofs of Lemmas \lemt{concave 2} and \lemt{finite 2} and some arguments in the proof of Lemmas \lemt{moment 1} and \lemt{martingale basic} because their proofs are not so hard. However, those lemmas are important in our arguments, and their detailed proofs may help the reader. Thus, we here include them as supplement to the paper.

\setcounter{section}{19}
\setcounter{subsection}{0}
\subsection{Proof of \lem{concave 2}}
\label{supp:S1}

We only prove (b) of Lemma 2.5 for $\theta < 0$ since others are obvious. Note that $\xi(v,\theta)$ is decreasing in $v$ as $v \uparrow \infty$ for this $\theta$. Hence, we always have
\begin{align*}
  0 \le \xi(\vartriangle,\theta) \le \xi(v,\theta), \qquad v > 0,
\end{align*}
and therefore $\xi(\vartriangle,\theta)$ and $\widehat{F}(\xi(\vartriangle,\theta))$ are finite for all $\theta < 0$.

Let us consider the two cases $\widehat{F}(\beta_{*}) = \infty$ and $\widehat{F}(\beta_{*}) < \infty$ separately. First, assume that $\widehat{F}(\beta_{*}) = \infty$. In this case, $\theta_{*} = - \infty$, and $\xi(\theta)$ exists and is finite for all $\theta < 0$. Hence,
\begin{align*}
  0 \le \xi(\theta) \le \xi(\vartriangle,\theta) \le \xi(v,\theta), \qquad v > 0,
\end{align*}
which implies that
\begin{align*}
  \dd{E}(e^{\xi(\theta) (T \wedge v)}) \le \dd{E}(e^{\xi(\vartriangle,\theta) (T \wedge v)}) \le \dd{E}(e^{\xi(v,\theta) (T \wedge v)}) = e^{-\theta}.
\end{align*}
Then, letting $v \uparrow \infty$, the monotone convergence theorem and the definition of $\xi(\theta)$ yield that
\begin{align*}
  e^{-\theta} = \dd{E}(e^{\xi(\theta) T}) \le \dd{E}(e^{\xi(\vartriangle,\theta) T}) \le e^{-\theta}
\end{align*}
Thus, we have $\dd{E}(e^{\xi(\theta) T}) = \dd{E}(e^{\xi(\vartriangle,\theta) T}) = e^{-\theta}$, which concludes that $\xi(\theta)=\xi(\vartriangle,\theta)$.

We next assume that $\widehat{F}(\beta_{*}) < \infty$. In this case, $\theta_{*} > -\infty$, and $e^{\theta_{*}} \widehat{F}(\beta_{*}) = 1$, equivalently, $\xi(\theta_{*}) = \beta_{*}$. If $\theta \in (\theta_{*}, 0)$, then $\xi(\theta)$ exists and is finite, and we have $\xi(\theta)=\xi(\vartriangle,\theta)$ similarly to the case $\widehat{F}(\beta_{*}) = \infty$. If $\theta = \theta_{*}$, then $\xi(\vartriangle,\theta_{*}) = \xi(\theta_{*})=\beta_{*}$. Finally, if $\theta < \theta_{*}$, then $\xi(\vartriangle,\theta_{*}) \le \xi(\vartriangle,\theta)$, and threfore
\begin{align*}
  0 \le \dd{E}(e^{\xi(\vartriangle,\theta_{*}) (T \wedge v)}) \le  \dd{E}(e^{\xi(\vartriangle,\theta) (T \wedge v)}), \qquad v > 0.
\end{align*}
Letting $v \uparrow \infty$, this implies
\begin{align*}
   e^{-\theta_{*}} = \widehat{F}(\beta_{*}) = \dd{E}(e^{\xi(\vartriangle,\theta_{*}) T}) \le \dd{E}(e^{\xi(\vartriangle,\theta) T}) = \widehat{F}(\xi(\vartriangle,\theta)).
\end{align*}
Hence, if $\beta_{*} < \xi(\vartriangle,\theta)$, then $\widehat{F}(\xi(\vartriangle,\theta)) = \infty$ because $\widehat{F}(s) = \infty$ for $s > \beta_{*}$. However, this contradicts the finiteness of $\widehat{F}(\xi(\vartriangle,\theta))$. Thus, $\beta_{*} = \xi(\vartriangle,\theta)$ for $\theta < \theta_{*}$, and we have completed a proof for (b) for $\theta < 0$.
\pend

\subsection{Outline of a proof for \lem{finite 2}}
\label{supp:S2}

Note that this lemma obviously holds if $\theta = 0$ because $\zeta_{i}(0) = \zeta_{i}(\vartriangle,0) = 0$, and therefore we assume that $\theta > 0$. We first consider the case that $0 < v < \infty$. If $K_{\alpha} = K$, then all the remaining service times are truncated in $W_{K_{\alpha},\theta}(\vc{R}(0-))$, and therefore the lemma is immediate. Thus, we can assume that $K_{\alpha} \not= K$. As in the proof of Lemma 4.3 of \cite{SadoSzpa1995}, it is sufficient to prove
\begin{align}
\label{eqn:sup 1}
  \sup_{t > 0} \dd{E}_{\vc{0}}(e^{W_{K_{\alpha}(v),\theta}(\vc{R}^{(s)}(t))}) < \infty,
\end{align}
where $\vc{R}^{(s)}(t) = \{R_{i}(t); i \in K\}$, and $\dd{E}_{0}$ stands for the conditional expectation given that $\vc{R}_{s}(0) = \vc{0}$. Similar to $W_{K_{\alpha}(v),\theta}$, let, for $A \in 2^{K} \setminus \{\emptyset\}$,
\begin{align*}
  W_{A_{\alpha}(v),\theta}(\vc{y}) = \sum_{i \in A \setminus K_{\alpha}}\zeta_{i}(\theta) y_{i} + \sum_{i \in A \cap K_{\alpha}} \zeta_{i}(v,\theta) (y_{i} \wedge v), \qquad \vc{y} \in \dd{R}_{+}^{k}.
\end{align*}
A key step to prove \eqn{sup 1} is the first inequality on page 553 of \cite{SadoSzpa1995}, which can be obtained in the present case as
\begin{align*}
  \dd{E}\big(& e^{W_{A_{\alpha}(v),\theta}(\vc{T}^{(s)} - \vc{u})} 1(\vc{T}^{(s)} > \vc{u}) \big)\\
  & = \dd{E}\big(e^{\sum_{i \in A \setminus K_{\alpha}}\zeta_{i}(\theta) (T_{i} - u_{i})  + \sum_{i \in A \cap K_{\alpha}} \zeta_{i}(v,\theta) ((T_{i} - u_{i}) \wedge v)} 1(\vc{T}^{(s)} > \vc{u}) \big)\\
  & \le e^{- \wedge_{i \in A \setminus K_{\alpha}} u_{i} \sum_{i \in A \setminus K_{\alpha}}\zeta_{i}(\theta)} \dd{E}\big(e^{\sum_{i \in A \setminus K_{\alpha}}\zeta_{i}(\theta) T_{i}  + \sum_{i \in A \cap K_{\alpha}} \zeta_{i}(v,\theta) (T_{i} \wedge v)} \big)\\
  & = e^{- \wedge_{i \in A \setminus K_{\alpha}} u_{i} \sum_{i \in A \setminus K_{\alpha}}\zeta_{i}(\theta)} e^{|A \setminus K_{\alpha}|\theta + |A \cap K_{\alpha}| v}, \qquad \vc{u} \in \dd{R}_{+}^{k},
\end{align*}
because $\dd{E}(e^{\zeta_{i}(\theta) T_{i}}) = e^{\theta}$ for $i \in A \setminus K_{\alpha}$, where $\vc{T}^{(s)} = \{T_{i}; i \in K\}$. Since $K_{\alpha} \not= K$, there is an $A$ such that $A \setminus K_{\alpha} \ne \emptyset$. Furthermore, $\sum_{i \in A \setminus K_{\alpha}}\zeta_{i}(\theta) > 0$ for such $A$ and $\theta \in (0,\alpha+\delta_{0})$. Hence, we can apply the exactly same argument to prove \eqn{sup 1} as in \cite{SadoSzpa1995}.

We next consider the case that $v = \vartriangle$. In this case, there is no truncation for the remaining service times. As we have shown in \supp{S1}, we have
\begin{align*}
  \dd{E}(e^{\zeta_{i}(\vartriangle,\theta) T_{i}}) = e^{\theta}, \qquad \theta \ge 0, i \in K_{\alpha}.
\end{align*}
Hence, applying the same argument as the case that $0 < v < \infty$, we have
\begin{align*}
  \dd{E}\big(& e^{W_{A_{\alpha}(v),\theta}(\vc{T}^{(s)} - \vc{u})} 1(\vc{T}^{(s)} > \vc{u}) \big) \le e^{- \wedge_{i \in A} u_{i} (\sum_{i \in A \setminus K_{\alpha}}\zeta_{i}(\theta) + \sum_{i \in A \cap K_{\alpha}} \zeta_{i}(\vartriangle,\theta))} e^{|A|\theta},
\end{align*}
which proves \eqn{sup 1} in the same way as in the case that $0 < v < \infty$.

\subsection{(A.9) in the proof of \lem{moment 1} in \app{moment 1}}
\label{supp:S3}

Recall that (A.9) of \cite{BravDaiMiya2015} is
\begin{align}
\label{eqn:S1}
  \exists y >0, \exists \epsilon \in (0,1), \qquad \liminf_{n \to \infty} \dd{P}(T_{i}^{(n)} > y) \ge \epsilon.
\end{align}
We prove that \eqn{S1} holds under the assumptions that \eqn{moment 1} holds for positive and finite $\lambda_{i}$ and $\{T_{i}^{(n)}; n \ge 1\}$ is uniformly integrable, that is,
\begin{align}
\label{eqn:S2}
  \lim_{a \to \infty} \sup_{n \ge 1} \dd{E}(T_{i}^{(n)}1(T_{i}^{(n)} > a)) = 0.
\end{align}
Suppose that \eqn{S1} does not hold. Then, 
\begin{align}
\label{eqn:S3}
  \forall y >0, \forall \epsilon \in (0,1), \qquad \liminf_{n \to \infty} \dd{P}(T_{i}^{(n)} > y) < \epsilon.
\end{align}
On the other hand, it follows from the uniform integrability condition \eqn{S2} that, for any $\epsilon' > 0$, there exists $a_{0} > 0$ such that, for any $a \ge a_{0}$,
\begin{align}
\label{eqn:S4}
  \forall n \ge 1, \qquad \dd{E}(T_{i}^{(n)} 1(T_{i}^{(n)} > a)) < \epsilon'.
\end{align}
Hence, from \eqn{S3} and \eqn{S4}, we have
\begin{align*}
  \dd{E}(T_{i}^{(n)}) &= \dd{E}(T_{i}^{(n)} 1(T_{i}^{(n)} \le y)) + \dd{E}(T_{i}^{(n)} 1(y < T_{i}^{(n)} \le a)) + \dd{E}(T_{i}^{(n)} 1(T_{i}^{(n)} > a))\\
  & \le y + a \dd{P}(y < T_{i}^{(n)}) + \epsilon',
\end{align*}
and therefore
\begin{align*}
  \liminf_{n \to \infty} \dd{E}(T_{i}^{(n)}) \le y + a \epsilon + \epsilon'.
\end{align*}
In this inequality, we first let $y, \epsilon \downarrow 0$, then let $\epsilon' \downarrow 0$. This concludes that $\liminf_{n \to \infty} \dd{E}(T_{i}^{(n)}) = 0$, which contradicts with $\lambda_{i} < \infty$. Hence, we have \eqn{S1}.

\subsection{(c) and (d) in the proof of \lem{martingale basic} in \app{f infty finite} and the supermartingale $E^{f_{K_{\alpha}(\vartriangle),\theta}}(t)$}
\label{supp:S4}

In the proof for (c), we claim without proof that $e^{{\zeta}_{i}(v,\theta) R_{i}(t) \wedge v}$ almost surely converges to $e^{{\zeta}_{i}(\vartriangle,\theta) R_{i}(t)}$ as $v \to \infty$. This is easily verified by the following inequalities.
\begin{align*}
  \zeta_{i}(\vartriangle,\theta) R_{i}(t) \wedge v \le \zeta_{i}(v,\theta) R_{i}(t) \wedge v \le \zeta_{i}(v,\theta) R_{i}(t), \qquad v > 0, \theta \ge 0.
\end{align*}

We next detail a proof of (d). If $K_{\alpha}(\vartriangle) = \emptyset$, $\{M(t)\}$ is obviously a martingale, so (d) holds. Thus, we assume that $K_{\alpha}(\vartriangle) \not= \emptyset$. We first note that \eqn{terminal 1} may not hold for $f = f_{K_{\alpha}(\vartriangle),\theta}$, but we have
\begin{align}
\label{eqn:terminal 2}
  Q f_{K_{\alpha}(\vartriangle),\theta}(\vc{x}) \le f_{K_{\alpha}(\vartriangle),\theta}(\vc{x}), \qquad \vc{x} \in \Gamma,
\end{align}
because $\ell$ in $\vc{x} \equiv (\ell,U,\vc{y})$ is decreased by 1 when service is completed and, for $i \in K_{\alpha}(\vartriangle)$,
\begin{align*}
  \dd{E}(e^{\zeta_{i}(\vartriangle,\theta)T_{i}}) = \left\{
\begin{array}{ll}
 \dd{E}(e^{\beta_{i}T_{i}}) < e^{\theta}, \quad & \theta_{i} < \theta,\\
 e^{\theta} &  \theta \le \theta_{i}, 
\end{array}
\right.
\end{align*}
where $\theta_{i} \le \alpha$ for $i \in K_{\alpha}(\vartriangle)$. Define 
\begin{align*}
  A(t) = \int_{0}^{t} (f_{K_{\alpha}(\vartriangle),\theta}(X(s-)) - Qf_{K_{\alpha}(\vartriangle),\theta}(X(s-)) dN^{*}(s), \qquad t \ge 0,
\end{align*}
then $A(t)$ is predictable and nondecreasing by \eqn{terminal 2}. Using $M_{0}(t)$ of \eqn{martingale 1} for $f = f_{K_{\alpha}(\vartriangle),\theta}$, we rewrite \eqn{martingale 2} for $f = f_{K_{\alpha}(\vartriangle),\theta}$ as
\begin{align*}
  M_{0}(t) + f_{K_{\alpha}(\vartriangle),\theta}(X(0)) = M_{K_{\alpha}(\vartriangle),\theta}(t) + A(t).
\end{align*}
Hence, 
\begin{align}
\label{eqn:supermartingale 1}
  M_{K_{\alpha}(\vartriangle),\theta}(t) - f_{K_{\alpha}(\vartriangle)\theta}(X(0)) = M_{0}(t) - A(t)
\end{align}
is an $\sr{F}_{t}$-supermartingale. This proves (d) of \lem{martingale basic}.

Let $M(t) = M_{K_{\alpha}(\vartriangle),\theta}(t) - f_{K_{\alpha}(\vartriangle)\theta}(X(0))$. We show that $E^{f_{K_{\alpha}(\vartriangle),\theta}}(t)$ is an $\sr{F}_{t}$-supermartingale. Let
\begin{align*}
  Y(t) = \frac 1{f_{K_{\alpha}(\vartriangle),\theta}(X(0))} \exp\left(- \int_{0}^{t} \frac {\sr{A} f_{K_{\alpha}(\vartriangle),\theta}(X(s))} {f_{K_{\alpha}(\vartriangle),\theta}(X(s))} ds \right).
\end{align*}
By \eqn{martingale infty}, we can see that
\begin{align*}
  \int_{0}^{t} Y(u) dM(u) & =  \int_{0}^{t} Y(u) d f_{K_{\alpha}(\vartriangle),\theta}(X(u)) \\
  & \quad -  \int_{0}^{t} \frac {\sr{A} f_{K_{\alpha}(\vartriangle),\theta}(X(u))} {f_{K_{\alpha}(\vartriangle),\theta}(X(0))} \exp\left(- \int_{0}^{u} \frac {\sr{A} f_{K_{\alpha}(\vartriangle),\theta}(X(s))} {f_{K_{\alpha}(\vartriangle),\theta}(X(s))} ds \right) du\\
  & = \int_{0}^{t} Y(u) d f_{K_{\alpha}(\vartriangle),\theta}(X(u)) - \int_{0}^{t} f_{K_{\alpha}(\vartriangle),\theta}(X(u)) dY(u)\\
  & = \left[ f_{K_{\alpha}(\vartriangle),\theta}(X(u)) Y(u) \right]_{0}^{t}\\
  &  = f_{K_{\alpha}(\vartriangle),\theta}(X(t)) Y(t) - f_{K_{\alpha}(\vartriangle),\theta}(X(0)) Y(0)\\
  &= E^{f_{K_{\alpha}(\vartriangle),\theta}}(t) - 1.
\end{align*}
Hence, by \eqn{supermartingale 1}, we have
\begin{align*}
  E^{f_{K_{\alpha}(\vartriangle),\theta}}(t) = 1 + \int_{0}^{t} Y(u) dM_{0}(u) - \int_{0}^{t} Y(u) dA(u).
\end{align*}
This implies that $E^{f_{K_{\alpha}(\vartriangle),\theta}}(t)$ is an $\sr{F}_{t}$-supermartingale since $M_{0}(t)$ is an $\sr{F}_{t}$-martingale, $Y(t)$ is continuous and $\int_{0}^{t} Y(u) dA(u)$ is nondecreasing in $t$ due to $Y(u) > 0$.

It is easy to see that $E^{f_{K_{\alpha}(\vartriangle),\theta}}(t)$ is a multiplicative functional of \cite{KuniWata1963} because it is an $\sr{F}_{t}$-supermartingale and has the exponential form. Hence, we can use it for change of measure as we have done by \eqn{change 1a}.

\renewcommand{\refname}{\large Additional references}.


\end{document}